\documentclass[11pt]{article}
\usepackage[OT1]{fontenc} 
\usepackage{amssymb}
\usepackage{amsmath}
\usepackage{amsthm}
\usepackage{amsfonts}
\usepackage{xcolor}
\usepackage{mathbbol}
\usepackage{bm}
\usepackage{hyperref}

\usepackage{ dsfont }
\usepackage{times}
\usepackage{graphicx}
\usepackage{bm}

\usepackage{dsfont}
\RequirePackage[algo2e]{algorithm2e}
\RequirePackage{enumitem}
\RequirePackage{algorithm}
\usepackage{titling}
\usepackage{tabularx}
\usepackage[numbers]{natbib}
\usepackage{cleveref}
\usepackage[utf8]{inputenc}
\usepackage{setspace}
\usepackage[margin=1in]{geometry}

\DeclareMathOperator*{\argmin}{argmin}

\title{Online Estimation with Rolling Validation: \\
Adaptive Nonparametric Estimation with Streaming Data}
\author{
    Tianyu Zhang\textsuperscript{1} and Jing Lei\textsuperscript{2} \vspace{5mm}\\
\textsuperscript{1}Department of Statistics and Applied Probability, University of California, Santa Barbara \\
\textsuperscript{2}Department of Statistics \& Data Science, Carnegie Mellon University
}
\date{}

\usepackage[page]{appendix}
\newtheorem{theorem}{Theorem}[section]
\newtheorem{assumption}{Assumption}[section]
\newtheorem{lemma}[theorem]{Lemma}
\newtheorem{corollary}[theorem]{Corollary}
\newtheorem{definition}[theorem]{Definition}
\newtheorem{remark}{Remark}[section]
\newtheorem{example}[theorem]{Example}
\usepackage{chngcntr} 
\counterwithin{figure}{section} 
\begin{document}
\fontsize{12}{22}\selectfont

\maketitle
\begin{abstract}
\fontsize{12}{18}\selectfont
Online nonparametric estimators are gaining popularity due to their efficient computation and competitive generalization abilities. An important example includes variants of stochastic gradient descent. These algorithms often take one sample point at a time and incrementally update the parameter estimate of interest. In this work, we consider model selection/hyperparameter tuning for such online algorithms. We propose a weighted rolling validation procedure, an online variant of leave-one-out cross-validation, that costs minimal extra computation for many typical stochastic gradient descent estimators and maintains their online nature. Similar to batch cross-validation, it can boost base estimators to achieve better heuristic performance and adaptive convergence rate. Our analysis is straightforward, relying mainly on some general statistical stability assumptions. The simulation study underscores the significance of diverging weights in practice and demonstrates its favorable sensitivity even when there is only a slim difference between candidate estimators. 
\end{abstract}
% \tableofcontents

\section{Introduction}
Online estimators are a collection of statistical learning methods where the estimate of the parameter of interest is sequentially updated during the reception of a stream of data (\citep[Section~3]{hoi2021online}). In contrast to traditional batch (or offline) estimators that learn from the entire training data set all at once, online estimators gradually improve themselves after more data points are processed. They have gained popularity in recent years especially due to their lower computational expense compared to traditional batch learning methods. This characteristic has made online learning particularly attractive in scenarios where data arrives continuously and computational resources are limited. Even with a complete, ``offline'' data set in hand, online methods are still routinely applied given their fast training speed and competitive prediction quality. Over the past two decades, there has been significant progress in the development of online estimators, both in terms of their implementation and our understanding of their statistical properties \cite{calandriello2017efficient, dieuleveut2016nonparametric, marteau2019globally, tarres2014online, ying2008online,langford2009sparse}. Several works \citep{bach2013non, dieuleveut2016nonparametric, zhang2022sieve} have shown that estimators trained with Stochastic Gradient Descent (SGD) can achieve certain statistical optimality.

The performance of nonparametric estimators, usually measured by predictive accuracy, often crucially depends on some hyperparameters that need to be specified in the algorithm, such as the number of basis functions in basis expansion-based algorithms and step size in gradient-based algorithms. Automated procedures that can select good hyperparameter values are both practically important and theoretically desirable. Cross-Validation (CV) is a commonly applied technique in batch learning to assess the generalization ability of a fitted model and perform model selection. The data is divided into multiple (say, $M$) equal-sized folds. The model is then trained on $M-1$ folds and evaluated on the remaining fold. This process rotates over the holdout fold and takes the average holdout validation accuracy as the overall assessment. Although CV has a long history \citep{stone1978cross, picard1984cross} and arguably the most implemented model selection methods in practice \cite{browne2000cross, rao2008dangers, arlot2010survey, roberts2017cross, ding2018model}, its practical performance and theoretical investigation leaves many questions to be answered \cite{bates2024cross}. There is strong empirical evidence that CV tends to give slightly biased model selection results, and several empirical adjustments have been proposed, such as the ``0.632 rule'' \cite{efron1997improvements} and the ``one-standard-error rule'' \citep{tibshirani2009bias}.
Theoretically, while CV is known to be risk-consistent \citep{bousquet2002stability,homrighausen2017risk} under very general conditions, when it comes to model selection, it is shown \citep{yang2007,gyorfi2002distribution} that batch CV may not consistently select the correct model under some natural scenarios.

In this work, we study the problem of online model/hyperparameter selection using a variant of CV. To begin with, we must clarify the meaning of model-selection consistency in online estimation, as it turns out to be subtly different from the batch setting. In the theory of nonparametric and high-dimensional statistics, the optimal hyperparameter is determined by three components: the complexity of the target parameter, the amount of noise in the data, and the sample size. Therefore, in a batch setting, where the sample size is known and fixed, one can simply compare different hyperparameters as the other two components are (assumed to be) constant. 

However, in the online setting, the sample size keeps increasing and the optimal hyperparameter must change accordingly as predicted by the theory. For example, the optimal Lasso \cite{tibshirani1996regression} penalty should scale inversely as the square root of the sample size. In this work, we propose to focus on selecting among \emph{sequences} of hyperparameters rather than among individual hyperparameters. More specifically, we are given a collection of sequences of hyperparameters, each specified by a particular function form with respect to the sample size, and the formal inference task is to find out the optimal sequence that gives the best prediction accuracy when enough information is available. When the $n$-th sample is revealed to the estimation algorithm, it updates the current estimate, for example, by taking an SGD step. The specifics of this update, such as the step size and/or the model capacity of the update, are determined by the $n$-th element of the hyperparameter sequence, $\lambda_n$. When the next sample arrives, the update procedure will iterate with $\lambda_{n+1}$. More detailed discussion and motivation of the tuning sequence selection problem and its relationship to the batch tuning problem is given in \Cref{subsec:batch_vs_seq_tuning}.

To select the best hyperparameter sequence, we propose a new methodology called weighted Rolling Validation (wRV), which can be viewed as a variant of leave-one-out (LOO) CV adapted to online learning. Unlike batch learning, where data is artificially split into fitting and validation folds, in the online setting the ``next data point'' naturally serves as the validating sample for the current estimate. RV efficiently exploits the computational advantage of online learning algorithms, because the standard LOO batch CV requires re-fitting the model $n$ times, whereas RV incurs no extra fitting at all. Our theoretical and methodological framework for RV allows for different weights to different time points in calculating the cumulative validated risks, which offers remarkable practical improvements. Our main theoretical result states that, under certain stability conditions on the fitting algorithms, wRV can select the best hyperparameter sequence with probability tending to one as the sample size diverges to infinity, provided the candidate sequences are sufficiently different. This result agrees with the existing results for batch CV.  If the candidate sequences are too similar, such as both being $\sqrt{n}$-consistent, then CV-based model selection is known to be inconsistent \citep{Shao93,Zhang93,yang2007}. We also extend the discussion to the settings where the number of candidate estimator sequences diverges along with sample size, which corresponds to the settings when refined hyperparameter sequences are applied or more modeling possibilities are considered as more data points become available. 
 
 \paragraph{Related work.}
Procedures that are formally similar to RV have been implemented in real-world data analysis, especially for time series data \cite{armstrong1972comparative, vien2021machine}, but none have they been implemented nor formally well-studied in the context of online adaptive nonparametric estimation. In the field of time series analysis, these procedures are referred to as ``expanding window cross-validation'' or ``walk-forward validation''. To the best of the authors' knowledge, the earliest related work was published in 1972 \citep[page 216]{armstrong1972comparative}, which implements a similar performance metric for long-range market forecasting tasks. The idea of RV is also mentioned in an earlier review paper \cite{tashman2000out}, under the name ``rolling-origin calibration''. In a more recent study \citep{benkeser2018online}, the authors applied unweighted RV---referred to as ``online cross-validation'' in their work---to ensemble learning and present some favorable guarantees of the proposed estimators. Compared with their work, we explicitly engage with nonparametric procedures, elucidate what the hyperparameters are in these settings, and go into more detail about the theoretical analysis. Moreover, our proposed wRV statistics can much better trace the transition of estimators' generalization error than the unweighted version, which may take $10$k more samples for the latter to make the correct choice in some simple settings (Section~\ref{section: numerical}). This deficiency of the unweighted RV is also mentioned in \cite[Section~5]{wei1992predictive}.

\paragraph{Notation.} We will use $\mathbb{R}$ ($\mathbb{R}^+$) to denote (positive) real numbers and $\mathbb{N}^+$ to denote the set of positive integers. For $n \in \mathbb{N}^+$, we denote $\{1,2,...,n\}$ as $[n]$. For any vector $x$ in $\mathbb{R}^p$ or $(\mathbb{N}^+)^p$, we use $x[j]$ to denote its $j$-th element. Notation $a\vee b = \max\{a,b\}$ for $a,b\in\mathbb{R}$. 
% For a random variable $X$ and a probability distribution $P$, the notation $X\sim P$ means ``the distribution of $X$ is $P$''. 
In some conditional expectation notation such as $\mathbb{E}[f(X,Y)|F^i]$, we use $F^i$ to denote the $\sigma$-algebra generated by sample $\{(X_j,Y_j), j\in [i]\}$. The notation $\mathbb{E}_{X_0}[\hat f(X_0)]$ means taking a conditional expectation with respect to $X_0$. We use $\delta_{jk}$ to represent Kronecker delta: $\delta_{jk} = 1$ when $j=k$ and $\delta_{jk} = 0$ otherwise. The inequality $a \lesssim b, a,b\in\mathbb{R}$ means $a\leq Cb$ for some constant $C$, but the value of $C$ may be different from line to line. When condition $\mathcal{A}$ holds, $\mathbb{1}(\mathcal{A}) = 1$, otherwise it takes $0$.

\section{Stochastic Approximation and the Problem of Online Model Selection}
\label{section:SGD}

In a typical statistical learning scenario, we have sample points of a pair of random variables $(X,Y)\in\mathbb R^{p}\times\mathbb R^1$ with joint distribution $P$, and want to minimize the predictive risk over a class $\mathcal F$ of regression functions: $\min_{f\in\mathcal F}\mathbb{E}_{(X,Y)\sim P} [l(f;X,Y)]$ under a loss function $l(\cdot)$.  
 We denote the marginal distribution of $X$ as $P_X$ and assume that $\mathbb{E}[Y^2]<\infty$. For regression problems, we consider the squared loss $l(f;X,Y) = (f(X) - Y)^2$. When $\mathcal{F}$ is the function space of square-integrable functions $L_2(P_X)$, the corresponding minimizer
\begin{equation*}
    f_0 = \argmin_{f\in L_2(P_X)} \mathbb{E}[(f(X) - Y)^2]
\end{equation*}
is the conditional mean function or the true regression function.

In contrast to batch learning, where all the data is available at once, streaming samples are usually revealed to the researchers continuously. For computational convenience, it is also increasingly more common to treat large data sets as a data stream and apply online methods by incrementally processing the samples. In both cases, the estimates need to be frequently updated whenever new samples arrive. When a new sample $(X_{i+1},Y_{i+1})\sim P(X,Y)$ arrives, the algorithm uses this new piece of information to update the current estimate $\hat f_{i}$ to a new estimate $\hat f_{i+1}$. One may convert any batch estimator to an online one by storing all the past data and refitting the whole model repeatedly. But this is often computationally infeasible. For example, fitting a kernel ridge regression estimator using $i$ samples (by matrix inversion) requires $i^3$-order computation. This would result in a cumulative $n^4$ computation to process $n$ samples in an online fashion. Moreover, batch estimators often need to load the whole data set into memory, causing significant storage costs, whereas genuine online algorithms operate incrementally and maintain a much smaller memory footprint.

As many iterative algorithms, the new estimate $\hat f_{i+1}$ can be written as the output of a pre-defined $\text{Update}(\cdot)$ mapping, whose inputs are the current estimate $\hat f_i$, the new sample $(X_{i+1},Y_{i+1})$ and hyperparameters $\lambda_i\in\mathbb{R}^{\Lambda}$, for some $\Lambda \in \mathbb{N}^+$. Formally,
\begin{equation}
\label{eq: abstract update rule}
\begin{aligned}
    \hat f_0&  \leftarrow 0\\
    \hat{f}_{i+1} &\leftarrow \operatorname{Update}\left(\hat{f}_i,\left(X_{i+1}, Y_{i+1}\right), \lambda_i\right).
\end{aligned}
\end{equation}
For simplicity, we set the initial estimate to the zero function without loss of generality. Such an online scheme has the appealing \emph{single-pass property} which leads to substantial computation and storage savings.
\begin{enumerate}
    \item The estimate update only depends on the current estimate $\hat f_i$ and the current data point.
    \item Each data point is used only once. The past estimates $\{\hat f_j:j\in [i-1]\}$  and data points $\{(X_{j},Y_j):j \in [i-1]\}$ are not stored.
\end{enumerate}
Online algorithms typically require a sequence of hyperparameters $\{\lambda_i: i \in \mathbb{N}^+\}$, which allows much flexibility in the algorithm.

\subsection{Examples of Online Estimators and their Hyperparameters}\label{section: example_online_estimators}

In this section, we review the basic ideas of stochastic approximation, SGD, and the explicit forms of several hyperparameter sequences. 
Consider a linear model $\mathcal{F}_{\rm linear}$: for each $f\in\mathcal{F}_{\rm linear}$, there exists a vector $\beta\in\mathbb{R}^p$ such that $f(X;\beta) = X^\top \beta$. This brings us to a parametric stochastic approximation problem:
\begin{equation*}
\min _{\beta \in \mathbb{R}^p} \mathbb{E} [l(f(\ \cdot\ ; \beta) ; X, Y)] .
\end{equation*}

The gradient descent method, a basic numerical optimization technique, uses the following iterative estimation procedure: First, initialize $\tilde\beta_0 = 0$, then iteratively update the estimate against the direction of the local gradient:
\begin{equation*}
    \tilde \beta_{i+1} \leftarrow \tilde\beta_i - \gamma_i\nabla_\beta \mathbb{E}\left[l(f(\ \cdot\ ;\beta); X,Y)\right]\big|_{\tilde{\beta}_i},
\end{equation*}
where $\gamma_i\in\mathbb{R}^+$ is a pre-specified sequence of step sizes (or learning rates). In practice, the expectation $\mathbb{E}[\ \cdot\ ]$ cannot be evaluated and must be approximated from the data. SGD procedures replace it with a one-sample (or mini-batch) estimate. In the case of linear regression:
\begin{equation*}
    \tilde \beta_{i+1} \leftarrow \tilde\beta_i - \gamma_i\left.\nabla_\beta\left(Y_i - X_{i+1}^{\top} \beta\right)^2\right|_{\tilde{\beta}_{i}} =\tilde{\beta}_i+2 \gamma_i\left(Y_{i+1}-X_{i+1}^{\top} \tilde{\beta}_{i}\right) X_{i+1}.
\end{equation*}
% Since $\tilde f_i(X) = X^\top\tilde\beta_i$ is uniquely determined by the coefficients, 
We have our first concrete example of the stochastic approximation rule \eqref{eq: abstract update rule}. Here, the hyperparameter sequence is a real number series $\lambda_i = \gamma_i\in\mathbb{R}^+$. Alternatively, we can rewrite the above iteration in a function-updating style:
\begin{equation*}
    \tilde{f}_{i+1}(\cdot) \leftarrow \tilde{f}_i(\cdot)+2 \gamma_i\left(Y_{i+1}-\tilde{f}_i\left(X_{i+1}\right)\right) \sum_{k=1}^p \phi_k\left(X_{i+1}\right) \phi_k(\cdot).
\end{equation*}
The function $\phi_k:\mathbb{R}^p\rightarrow\mathbb{R}$ maps a $p$-dimensional vector to its $k$-th entry $\phi_k(X) = X[k]$. We present it in this more intricate format to better align with the nonparametric versions that follow. The Polyak-averaging of $\{\tilde f_i\}$:
\begin{equation}
\label{eq: polyak averaging}
    \hat f_i = i^{-1}\sum_{m=1}^i \tilde f_m = \frac{i-1}{i}\hat f_{i-1} + \frac{1}{i} \tilde f_i
\end{equation}
is shown \cite{bach2013non} to achieve the parametric minimax rate of estimating $f_0$ when the model is well-specified. 

Now we turn to nonparametric online estimation with a diverging model capacity. First, let's consider the reproducing kernel SGD estimator (kernel-SGD). Let  $\mathcal K(\cdot,\cdot):\mathbb{R}^p\times\mathbb{R}^p\rightarrow\mathbb{R}$ be a reproducing kernel with associated eigenvalue–eigenfunction pairs $\{(q_k, \phi_k):k\in \mathbb{N}^+\}$. The kernel-SGD update rule is
\begin{equation}
\label{eq: kernel sgd update}
\begin{aligned}
      \tilde{f}_{i+1}(\cdot) & \leftarrow \tilde{f}_i(\cdot)+\gamma_i\left(Y_{i+1}-\tilde{f}_i\left(X_{i+1}\right)\right) \mathcal K\left(X_{i+1}, \ \cdot\ \right) \\
&=\tilde{f}_i(\cdot)+\gamma_i\left(Y_{i+1}-\tilde{f}_i\left(X_{i+1}\right)\right) \sum_{k=1}^{\infty} q_k \phi_k\left(X_{i+1}\right) \phi_k(\cdot).
\end{aligned}
\end{equation}
This is our second example of the abstract update rule \eqref{eq: abstract update rule}. The related rate-optimal estimator \cite{dieuleveut2016nonparametric} is its Polyak-averaging $\hat f_i = i^{-1}\sum_{m=1}^i\tilde f_m$, which can also be updated as described earlier \eqref{eq: polyak averaging}.

The performance of kernel-SGD depends on the selected kernel function $\mathcal{K}$ and the learning rate $\gamma_i$. It has been shown that using a sequence of slowly decaying step sizes $\gamma_i = A i^{-\zeta}$ yields rate-optimal estimators, where $A \in \mathbb{R}^+$ and $\zeta \in (0, 1/2)$. However, the best choice of $\zeta$ depends on the expansion coefficients $\langle f_0,\phi_k\rangle_{L_2(P_X)}$ and its relative magnitude to the eigenvalues $q_k$, which is in general not available to the algorithm. Data-adaptive procedures for choosing $\gamma_i$ are needed for better performance of the kernel-SGD in practice. Our framework also allows a varying reproducing kernel, such as Gaussian kernels with a dynamic bandwidth. In the notation of \eqref{eq: abstract update rule}, the hyperparameter $\lambda_i = \gamma_i$ is one-dimensional when the kernel function is fixed, and becomes two-dimensional, $\lambda_i = (\gamma_i, \sigma_i)^\top$, when the kernel bandwidth $\sigma_i$ is allowed to vary.

The second line of the kernel-SGD update \eqref{eq: kernel sgd update} also suggests a direct way to construct the update from basis expansions without specifying a kernel.
Recently, \cite{zhang2022sieve, chen2024stochastic} consider ``sieve-type" nonparametric online estimators that explicitly use an increasing number of  orthonormal basis functions $\{\phi_k:k\ge 1\}$, with an update rule
\begin{equation}
\label{eq: sieve sgd update}
\begin{aligned}
      \tilde{f}_{i+1}(\cdot) \leftarrow \tilde{f}_i(\cdot)+\gamma_i\left(Y_{i+1}-\tilde{f}_i\left(X_{i+1}\right)\right) \sum_{k=1}^{J_i} k^{-2 \omega} \phi_k\left(X_{i+1}\right) \phi_k(\cdot), 
\end{aligned}
\end{equation}
and the averaged version $\hat f_i$ is similarly updated as in \eqref{eq: polyak averaging}. Here, $\gamma_i$ denotes the learning rate sequence as before, $J_i$ is a positive integer representing the number of basis functions used in the $(i+1)$-th update, and $\omega$ is a shrinkage parameter, which can typically be set to $1/2$ in practice. Sieve-SGD \eqref{eq: sieve sgd update} is computationally more efficient since $J_i$ is usually much smaller than $n$ for moderate-dimension problems (where we do not model $p$ increasing with $n$). Under the assumption that $f_0$ belong to some Sobolev ellipsoid $W(s)$ spanned by $\{\phi_k\}$:
\begin{equation}\label{eq: sobolev ellipsoid}
    W(s) = W(s;Q) = \left\{f=\sum_{k=1}^{\infty} \beta_k \phi_k \text{ such that } \sum_{k=1}^{\infty}\left(\beta_k k^s\right)^2 \leq Q^2\right\},
\end{equation}
it has been shown \cite{zhang2022sieve} that the sieve-SGD has minimal space expense (computer memory cost) among all statistically optimal estimators. The theoretical optimal hyperparameters are 
\begin{equation}
\label{eq: step size and basis number}
\left(\gamma_i, J_i\right)=\left(A i^{-1 /(2 s+1)}, B i^{1 /(2 s+1)}\right),
\end{equation}  
where $A,B\in\mathbb{R}^+$ and $s\geq 1/2$. The hyperparameter sequence $\{\lambda_i\}$ of sieve-SGD \eqref{eq: sieve sgd update} is typically two-dimensional, with $\lambda_i = (\gamma_i, J_i)^\top$.

\subsection{Model Selection in Batch and Online Settings}\label{subsec:batch_vs_seq_tuning}

Suppose we have a dataset of size $i$ (or have collected $i$ samples from a data stream), along with $K$ estimates $\hat{f}_i^{(k)}$ for $k \in [K]$. One of the main goals of model selection is identifying the optimal index $k_i^*$ that minimizes the expected risk:
\begin{equation}\label{eq: finite_best_model}
    k_i^*  = \argmin_{k\in[K]}\mathbb{E}[r_{i,k}]
\end{equation}
where the risk is defined as
\begin{equation*}
\begin{aligned}
r(\hat{f})&=\mathbb{E}_{X_0}\left[\left(\hat{f}\left(X_0\right)-f_0\left(X_0\right)\right)^2\right] \\
r_{i, k} & =r\big(\hat{f}_i^{(k)}\big).
\end{aligned}
\end{equation*}
Here $X_0$ is a random variable independent from the estimator $\hat f$. And $r(\hat f)$, the risk of $\hat f$, still possesses the randomness from constructing $\hat f$. The expectation in \Cref{eq: finite_best_model} is marginal, making $\mathbb{E}[r_{i,k}]$ a real number that only depends on $(i,k)$. 

For batch learning, this task can be accomplished by comparing some (CV) estimates of the risks \cite{Shao93}. Such effective risk estimates are available because data can be repeatedly visited and we can empirically evaluate the performance of $\hat f_i^{(k)}$---potentially trained with a slightly smaller sample size---on many left-out samples. However, for online learning algorithms \eqref{eq: abstract update rule} satisfying the single-pass property, we can only obtain a single realization of the loss function at each sample size $i$. Methods measuring the overall risk behavior of online algorithms are often implemented in this setting.

Another way to formulate the model selection problem in the infinite-horizon setting is to asymptotically identify an index $k^* \in [K]$---if it exists---such that
\begin{equation}\label{eq:seq_selec_obj}
\limsup_{i\rightarrow\infty}\ \mathbb{E}[r_{i,k^*}]/\mathbb{E}[r_{i,k}] < 1, \text{ for all } k\in[K]\backslash\{k^*\}\,. 
\end{equation}
While this criterion can be used to assess both batch and online model selection procedures, it is particularly natural for the online setting because the sample size is an actual varying component of the problem. 

For instance, in the sieve-SGD example discussed in \Cref{section: example_online_estimators}, suppose the only hyperparameter to be selected is ${J_i}$, the number of basis functions used in the estimator. Existing theory \cite{zhang2022sieve, chen2024stochastic} suggests that the best choice of $J_i$ has the form $J_i=B i^{1/(2s+1)}$ for a pair of constants $(B,s)$. Thus the candidate set of models corresponds to a collection of pairs $\{(B^{(1)},s^{(1)}),...,(B^{(K)},s^{(K)})\}$, each defining one of $K$ infinite sequences specifying the number of basis functions:
\begin{equation*}
\mathbf J^{(k)}=\left\{J_i^{(k)}=B^{(k)} i^{1/\left(2 s^{(k)}+1\right)}, i\in \mathbb{N}\right\}\,,~~ k\in[K].
\end{equation*}
Correspondingly, we obtain $K$ \emph{sequences} of estimators $\{\hat f^{(k)}_i, i\in\mathbb{N}\}, k\in [K]$, where each of them is trained with a different $\mathbf{J}^{(k)}$. An online model selection procedure should guide us toward the sequence $\mathbf{J}^{(k)}$ that eventually performs best, thereby identifying the optimal combination of $(B^{(k)}, s^{(k)})$ based on the observed data. In contrast, batch model selection methods focus more on a single sample size and do not suggest a good choice of model for a different sample size. For example, if batch model selection returns a choice of $J=10$ at 
sample size $100$, it is unclear whether this $J$ comes from a $(B,s)$ combination of $(1,~1/2)$ or $(\sqrt{10},~3/2)$, which would lead to different choices of $J$ for other values of sample size.

\section{Methods: Weighted Rolling Validation}
\label{section:RV}

For a stochastic approximation estimator \eqref{eq: abstract update rule}, new sample points are sequentially applied to update the current estimates. The idea of wRV is simple: we treat the new sample as a natural left-out to update the prediction accuracy metric before applying it to update the current estimates.

Consider $K$ sequences of estimators for $f_0$: $\{\hat f^{(k)}_i\}, k\in [K]$ with iterative formula
\begin{equation*}
\hat{f}_{i+1}^{(k)}=\operatorname{Update}\left(\hat{f}_i^{(k)},\left(X_{i+1}, Y_{i+1}\right), \lambda_i^{(k)}\right) .
\end{equation*}

To compare their online prediction accuracy, we calculate the wRV statistic sequence for each estimator trajectory:
\begin{equation}
\label{eq: rolling validation}
{\rm RV}_i^{(k)}={\rm RV}\left(\left\{\hat{f}_l^{(k)}, l\in [i]\right\}, \xi\right)=\sum_{l=0}^il^{\xi}\left(\hat{f}_l^{(k)}\left(X_{l+1}\right)-Y_{l+1}\right)^2,
\end{equation}
where $\xi\ge 0$ is a weighting exponent to be chosen by the user. We will discuss some rule-of-thumb choices for it soon in \Cref{section: detection delay}. A positive $\xi$ will put more weight on larger values of $l$---i.e., the samples processed later---which can improve finite-sample performance by alleviating some ``detection delay" due to the cumulative nature of wRV. Although the definition of wRV statistics used squared loss, the method is general and can be extended to other loss functions (see \Cref{section: discussion} for a quantile-regression example).

At each step $i$, we can use $\hat k^*_i=\arg\min_{k\in [K]} {\rm RV}_i^{(k)}$ as an estimate of $k^*_i$ \eqref{eq: finite_best_model} or $k^*$ \eqref{eq:seq_selec_obj}. A statistically favorable procedure is expected to achieve $\mathbb{P}(\hat k^*_i=k^*)=1 - o(1)$ when the superior sequence exists.

It is direct to see that ${\rm RV}_i^{(k)}$ is an unbiased estimator of accumulated prediction error:
\begin{equation}\label{eq: expectation of RV stat}
    \mathbb{E}\left[{\rm RV}_i^{(k)}\right]=\sum_{l=1}^i l^\xi     \mathbb{E}\left[\left(\hat{f}_l^{(k)}(X)-Y\right)^2\right].
\end{equation}
For wRV to select the correct $k^*$, we need the following two conditions:
\begin{enumerate}
	\item There is a gap between $    \mathbb{E}\left[{\rm RV}_i^{(k^*)}\right]$ and $\inf_{k\neq k^*}    \mathbb{E}\left[{\rm RV}_i^{(k)}\right]$, which is closely related to the excess risk condition \eqref{eq:seq_selec_obj}.
	\item The sample quantities ${\rm RV}_i^{(k)}$ do not deviate too far away from their expected values---specifically, their deviation should remain smaller than the aforementioned "gap." This is technically more challenging as the dependence structure in ${\rm RV}_i^{(k)}$ is more complicated than simple martingales. We resort to a stability-based argument to formally quantify its variability.
\end{enumerate}
A formal analysis will be presented in detail in Section~\ref{section: consistency}, with an extension to diverging $K$ in \Cref{section: diverging cardinality}.

\subsection{Choice of Weighting Exponent}\label{section: detection delay}

A positive weighting exponent $\xi$ can significantly enhance finite-sample model selection quality, which is also demonstrated later in the numerical examples (\Cref{section: numerical}). The optimal $\xi$ depends on the specific problem and the estimators utilized, making it challenging to determine in practice. In general, we recommend using a fixed $\xi = 1$ to balance selection sensitivity and stability. This choice is determined under certain assumptions on the excess risks. We provide some intuition below, and the technical details can be found in \Cref{app: delay ratio}. 

Assume we receive a sequence of IID samples $(X_i,Y_i) = (X_i, f_0(X_i) + \epsilon_i)$, where $\epsilon_i$ is a centered noised variable $\mathbb{E}[\epsilon_i \mid X_i] = 0$. Consider a sequences of estimators $\{f_i^{(1)}\}$ whose risk at step $i$ is
\begin{equation*}
\mathbb{E}\left[\left(f_i^{(1)}(X)-Y\right)^2\right]  =\mathbb{E}[r_{i,1}]+\mathbb{E}\left[\epsilon_i^2\right].
\end{equation*}

We similarly denote the excess risk of another sequence $\{f_i^{(2)}\}$ as $\mathbb{E}[r_{i,2}]$. After revealing $i$ samples, ideally, we want to instantly identify $k_i^*$ and use it for predictions or downstream procedures. However, wRV only provides a noisy version of cumulative risks $\mathbb{E}\left[{\rm RV}_i^{(k)}\right]$, $k=1,2$, whose relative orders may differ from that of $\mathbb{E}[r_{i,k}]$ for finite $i$.

To simplify the discussion, we assume both models have polynomial orders of excess risks: $\mathbb{E}[r_{i,1}] = Ai^{-a}$ and $\mathbb{E}[r_{i,2}] = Bi^{-b}$ with $0<A<B$ and $0\leq a<b<1$. For smaller sample sizes $i\leq(B / A)^{1/(b-a)}$, $\{f_i^{(1)}\}$ is a better estimator, but their relative performance flips afterward because $f_i^{(2)}$ converges to $f_0$ at a faster rate. However, the relative magnitudes of $\mathbb{E}\left[{\rm RV}_i^{(1)}\right]$ and $\mathbb{E}\left[{\rm RV}_i^{(2)}\right]$ do not switch until $i\ge (B/A)^{1/(b-a)} T(a,b,\xi)$ with a delay ratio $T(a,b,\xi) > 1$ (whose value interestingly does not depend on $A$ or $B$). 

A calculation (\Cref{lemma: delay ratio}) shows that increasing the weighting exponent $\xi$ can reduce the delay ratio. Specifically, for any target threshold $t > 1$, the maximum delay ratio $T(a, b, \xi)$ across all $(a, b)$ can be controlled to stay below $t$ by choosing $\xi > 1 / \log t$, as illustrated in \Cref{fig: threshold xi}. To enhance the algorithm's sensitivity, a larger $\xi$ that minimizes the delay ratio is preferred. However, this comes at the cost of reducing the concentration of ${\rm RV}_i^{(k)}$ on $\mathbb{E}[{\rm RV}_i^{(k)}]$ by a constant factor. For example, if $\xi = 3$, the $i$-th sample has a weight of $i^3$ when calculating ${\rm RV}_i$, while the $i/2$-th sample has a weight of $i^3/8$, contributing much less to the wRV criterion. We recommend selecting $\xi = 1$ to control the worst-case $T(a,b,\xi)$ at no more than $2.4$ while maintaining good stability.

We can also compare the excess risk of the selected sequence to that of the optimal one. When $i$ is slightly smaller than $I_F = (B / A)^{1 /(b-a)} T(a, b, \xi)$, the difference between $\mathbb{E}[r_{i,1}]$ and $\mathbb{E}[r_{i,2}]$ is maximized, while wRV on average still prefer the sub-optimal sequence, as discussed above. We can compute $\mathbb{E}[r_{I_F,1}]/\mathbb{E}[r_{I_F,2}]$ and use it as a measure of the worst-case efficiency loss. It takes a simple form
\begin{equation*}
\mathbb{E}[r_{I_F,1}]/\mathbb{E}[r_{I_F,2}] = 1 + \frac{b-a}{\xi +1-b} \leq 1 + \xi^{-1}.
\end{equation*}
A choice of $\xi = 1$ limits the efficiency loss to within a factor of two of the optimal.

The wRV-based decision improves dramatically when $\xi$ takes a positive value, but the benefit stabilizes after it is larger than $1$. We illustrate this phenomenon in \Cref{section: numerical} and \Cref{app: extra numerical}---the ranking plots for $\xi = 0$ differ significantly from the rest.

\subsection{Computational Advantages of Rolling Validation}

\begin{algorithm}[t]
\SetAlgoLined
\SetKwInOut{Input}{Input}
\SetKwInOut{Output}{Output}
\SetKwInOut{Initialize}{Initialize}
\vspace{3pt}

\Input{
A stream of training samples $(X_i, Y_i)$ for $i \in \mathbb{N}^+$.\\
\ An estimator update rule: $\text{Update}(f, (X, Y), \lambda)$.\\
    \ Candidate hyperparameter sequences $\{\lambda^{(k)}_i : i \in \mathbb{N}^+\}$ for $k \in [K]$.\\
    \ Weighting exponent $\xi$.
}

\Output{
    Online model selection index $\hat{k}^*_i$ at each time step $i$.
}

\Initialize{
Set estimators $\hat{f}^{(1)} = \dots = \hat{f}^{(K)} = 0$.\\
\ Set wRV statistics ${\rm RV}^{(1)} = \dots = {\rm RV}^{(K)} = 0$.\\
\ Set index $i = 0$.
}
\While{there are remaining training samples}{
    Read new sample $(X_{i+1}, Y_{i+1})$\;
    
    \For{$k \in [K]$}{
        Predict $\hat{Y}_{i+1}^{(k)} \leftarrow \hat{f}^{(k)}(X_{i+1})$\;
        
        Update wRV:
        ${\rm RV}^{(k)} \leftarrow {\rm RV}^{(k)} + i^\xi (Y_{i+1} - \hat{Y}_{i+1}^{(k)})^2$\;
        
        Update estimator: 
        $\hat{f}^{(k)} \leftarrow \text{Update}(\hat{f}^{(k)}, (X_{i+1}, Y_{i+1}), \lambda^{(k)}_{i+1})$\;
    }
    
    Output current selection: $\hat{k}_i^* \leftarrow \arg\min_k {\rm RV}^{(k)}$\;
    
    Increment index: $i \leftarrow i + 1$\;
}
\caption{Stochastic approximation with weighted rolling validation (wRV) for regression.}
\label{algorithm1}
\end{algorithm}

The wRV statistics \eqref{eq: rolling validation} is not only a natural target to examine when comparing online estimators, but also offers much computational savings when implemented efficiently. In Algorithm~\ref{algorithm1}, we summarize our recommended framework on integrating model training and model selection. In contrast to the batch-learning scenario where training, validation, and testing are typically done in separated phases, we proposed a procedure that: 1) updates multiple estimators simultaneously; 2) calculates the wRV statistics interchangeably with model training; and 3) can offer a selected model at any time during the training process.

We take the kernel-SGD estimator \eqref{eq: kernel sgd update} as an example to quantify the computational gain. Examining the update rule, we can verify that both the trajectory functions $\{\tilde f_i\}$ and the estimators $\{\hat f_i\}$ are a linear combination of $i$ kernel functions ``centered'' at the past sample covariate vectors. Formally, there are some $\tilde\beta_{mi}, \hat\beta_{mi}\in\mathbb{R}:$
\begin{equation*}
    \tilde f_i(\cdot) = \sum_{m=1}^i \tilde \beta_{mi} \mathcal K(X_m,\cdot\ ) ,\quad \hat f_i(\cdot) = \sum_{m=1}^i \hat \beta_{mi} \mathcal K(X_m,\cdot\ ).
\end{equation*}

Suppose evaluating one kernel function $\mathcal K(x,z)$ takes $O(p)$ computation---recall that $p$ is the dimension of $X_i$. Given a new sample $(X_{i+1}, Y_{i+1})$, the majority of computation for the $\text{Update}(\cdot)$ operator is calculating $\mathcal K(X_{i+1}, X_m)$ for all $m \in [i]$. This means $O((p+C)i)$ calculation when processing one sample $(X_{i+1}, Y_{i+1})$, where $C$ is a small constant corresponding to the extra computation after evaluating all $\mathcal K(X_{i+1}, X_m)$. Moreover, the above computational expense accumulates to the order of $O((p+C)n^2)$ for processing the first $n$ samples $\{(X_1,Y_1),...,(X_n,Y_n)\}$. Naively training $K$ kernel-SGD estimators one at a time with $n$ samples would require $O(K(p+C)n^2)$ computation if no intermediate quantities are saved. However, when the $K$ sequences of trajectories only differ at the learning rate $\gamma_i$, they can be simultaneously updated with only $O((p+KC)n^2)$ computation using the scheme in \Cref{algorithm1}, simply because the evaluated $\mathcal K(X_{i+1}, X_m)$ kernel value at each step can be shared between the estimators.

In addition to the computational gain for simultaneous training, these evaluated kernel values can be further shared with the wRV statistic calculation that only increases the total expense to at most $O((p+2KC)n^2)$. This effectively makes performance assessment a ``free lunch" relative to the cost of training, without a separate hold-out set for model selection. The only extra space expense is saving extra $K$ numbers $\{{\rm RV}^{(k)}\}$. A similar calculation also holds for sieve-SGD \eqref{eq: sieve sgd update}: the computational expense of training and validating $K$ estimators can be bounded by $O((p+2KC)nJ_n)$, where $J_n$ is the largest basis number among the $K$ candidate estimators. For a fixed $p$, $J_n$ is typically a diverging sequence with $J_n=o(\sqrt{n})$.

\section{Consistency of Rolling Validation, Fixed Number of Candidates}
\label{section: consistency}

In this section, we analyze the statistical properties of the proposed RV procedure in regression settings. We focus on whether wRV can select $k^*$ specified in \eqref{eq:seq_selec_obj} from a finite candidate; the increasing candidate-cardinality case is covered in \Cref{section: diverging cardinality}. We will show that under sufficient risk separation and stability of the wRV statistic, the proposed procedure can consistently select the superior model. For simplicity, we consider $K=2$ with two sequences of online estimators of $f_0$. This consistency result extends to any finite $K \geq 2$ by considering $k^*$ against each candidate and using the union bound. Below, we list the technical conditions.

First, we assume IID samples from a fixed underlying distribution.
\begin{assumption}\label{ass:iid}
    The data points $\{(X_i, Y_i)\in\mathbb{R}^{p+1}, i\in\mathbb{N}^+\}$ are IID samples from a common distribution $P(X,Y)$. The centered noise variable $\epsilon_1 = Y_1 - f_0(X_1)$ has finite variance: $\mathbb{E}[\epsilon_1^2 | X_1 =~x] \leq \sigma^2$ for any $x$ in the support of $X_1$.
\end{assumption}

We allow the noise variables $\epsilon_i$ to be associated with the covariate $X_i$ and only require them to have a finite second moment given $X_i$. This assumption is milder than those in some existing works on model selection via CV \citep{vaart2006oracle, benkeser2018online}.

Next, we formalize the assumption on the sub-optimal model’s estimation accuracy, essentially stating that its average risk is no better than $M i^{-a}$.

\begin{assumption}\label{ass: inf estimator_quality}
    Let $\{\hat f_i, i \in \mathbb{N}\}$ be a sequence of estimators of $f_0$. There exist constants $a\in[0,1]$ and $0<M<\infty$ such that
\begin{equation}
\label{eq: inf estimator quality}
    \liminf_{i \rightarrow \infty} i^a \mathbb{E}[r(\hat{f}_i)]=M.
\end{equation}
\end{assumption}

Recall $r(\hat{f}_i)=$ $\mathbb{E}_{X_0}\left[\left\{\hat{f}_i(X_0)-f_0(X_0)\right\}^2\right]$ and the outer expectation is unconditional, making the left-hand-side in \eqref{eq: inf estimator quality} non-random. We impose a similar condition for the superior sequence, flipping the limit inferior to limit superior:

\begin{assumption}\label{ass: sup estimator second moment}
   Let $\{\hat f_i, i \in \mathbb{N}\}$ be a sequence of estimators of $f_0$. There exist constants $a\in[0,1]$ and $0<M<\infty$ such that
\begin{equation}
\label{eq: sup estimator second moment}
    \limsup_{i \rightarrow \infty} i^a \mathbb{E}[r(\hat{f}_i)]=M.
\end{equation}
\end{assumption}

We also require the fourth moment of the estimators to converge to zero at a sufficient rate to ensure that the sample version of wRV concentrates closely around its mean.

\begin{assumption}\label{ass: sup estimator_quality}
   %We assume $\mathbb{E}[r(\hat{f}_i)]$ is uniformly bounded. Moreover, we require that t
   We assume $\sup_i \mathbb{E}[r(\hat{f}_i)]<\infty$  and there exist constants $a\in [0,1]$, $0<C_1<\infty$ such that
\begin{equation}
\label{eq: forth moment conditions}
    \limsup_{i \rightarrow \infty} i^{2a} \mathbb{E}\left[\left\{\hat{f}_{i}\left(X\right)-f_0\left(X\right)\right\}^4\right]=C_1.
\end{equation}
\end{assumption}

\begin{remark}
    Although \Cref{ass: sup estimator_quality} implies \Cref{ass: sup estimator second moment} with $M = C_1^{1/2}$, we state these conditions separately because, in Theorem~\ref{th: main theorem}, the superior and inferior sequences are assumed to satisfy different combinations of them. We will also compare the constants in \eqref{eq: inf estimator quality} and \eqref{eq: sup estimator second moment}, while the constant $C_1$ in \eqref{eq: forth moment conditions} is of lesser importance. Additionally, we do not require the sequence $i^a \mathbb{E}[r(\hat{f}_i)]$ itself to have a limit as $i \rightarrow \infty$.
\end{remark}
\begin{remark}
    The values of $(M,a)$ reflect the quality of the estimators, which partially reflects the problem's difficulty.  They are determined by the fitting methods and tuning parameters. For optimal parametric SGD, $a = 1$. For a minimax-optimal estimation procedure in a Sobolev ellipsoid \eqref{eq: sobolev ellipsoid}---as well as some related RKHS space---$a = 2s/(2s+1)$. The constant $M$ is less explicit but usually depends on the axis lengths of the ellipsoid, interpreted as the Sobolev/RKHS norm of $f_0$ or as a measure of $f_0$'s smoothness. Mean-squared error (MSE) convergence results for online (nonparametric) estimators are available in the literature, including \cite{dieuleveut2016nonparametric,zou2021benign, zhang2022sieve, zhang2023online,chen2024stochastic}. Specifically, \cite{zou2021benign} also provides an interesting lower bound for the MSE error.
\end{remark}

In addition to the moment conditions above, we also need the following stability conditions to control certain ``non-standard" deviation quantities that cannot be fully reduced to martingale-type concentration. Let $(X'_j,Y'_j)$ be an IID copy of $(X_j,Y_j)$ for a $j \geq 1$. The actual estimator $\hat f_i(x)$ is trained with samples $\bm{Z}_i = \{(X_l,Y_l), l \in [i]\}$ and can be explicitly expressed as $\hat f_i(x) = \hat f_i(x;\bm{Z}_i)$. We denote an imaginary dataset $\bm{Z}_i^j$, obtained by replacing $(X_j,Y_j)$ in $\bm{Z}_i$ by its IID copy $(X'_j,Y'_j)$. The imaginary estimator trained with $\bm{Z}_i^j$ is $\hat f_i(x;\bm{Z}_i^j)$. We need conditions on the difference between $\hat f_i(x;\bm{Z}_i)$ and its perturb-one version: $\nabla_j \hat{f}_{i}\left(x\right):= \hat{f}_{i}\left(x;\bm{Z}_i\right) - \hat{f}_{i}\left(x;\bm{Z}_i^j\right) $.

\begin{assumption}[Estimator stability]\label{ass:stability}
There exist constants $b, C_2> 0$ such that for any $i\in\mathbb{N}^+$ and $j\in[i]$, the following holds:
\begin{equation}\label{eq:stability_condition}
    i^{2b}\mathbb{E}[\{\nabla_j \hat f_i(X)\}^2| F^j] \leq C_2\qquad \text{almost surely},
\end{equation}	
    where $F^j$ is the $\sigma$-algebra generated by the first $j$ samples in $\bm{Z}_i$.
\end{assumption}

\begin{remark}
   The stability condition quantifies the variability of the estimator when switching one training sample. The exponent $b$ in \Cref{ass:stability} is of the most interest in our analysis. For many parametric estimators, the stability rate $b$ equals $1$, while for general nonparametric procedures, the rate is less than $1$, indicating less stability. In Section~\ref{sectoin:example stability of estimators}, we establish bounds on $b$ for estimator examples, including batch projection estimators and both parametric and non-parametric SGD estimators. A similar condition to \Cref{ass:stability} is also utilized in batch high-dimensional CV procedures \cite{kissel2022high}. 
   \end{remark}

\begin{remark}\label{remark: holder tradeoff}
   The sup-norm stability condition in \eqref{eq:stability_condition} can be simplified and relaxed to $\mathbb{E}[\{\nabla_j\hat f_{i}(X)\}^2] = O(i^{-2b})$ if we strengthen the convergence condition in \Cref{ass: sup estimator_quality} to the sup norm: $\sup_x \{|\hat f_i(x)-f_0(x)|\}= O(i^{-a})$. In fact, we require the product $\mathbb{E}[\nabla_j \hat{f}_i(X)\{\hat{f}_i(X)-f_0(X)\}\mid F^i]$ to be small enough in the proof (equation \eqref{important decomposition for stability}), which allows for H\"older's-inequality type trade-off between the stability term and the estimation error.
\end{remark}

We now present our main result, which establishes the consistency of wRV for a fixed number of candidate sequences. As the online sample size $n$ increases, it selects the better estimator sequence with probability approaching one.

\begin{theorem}
\label{th: main theorem}
Consider the IID learning setting in \Cref{ass:iid}. Let $\{h_i\}$ be an inferior sequence of estimators satisfying \Cref{ass: inf estimator_quality}~\&~\ref{ass: sup estimator_quality} with parameters $a = a_h\in [0,1)$ and $ M = M_h$, and \Cref{ass:stability} with $b = b_h > 1/2 + a_h/2$. Let $\{g_i\}$ be the superior sequence of estimators \eqref{eq:seq_selec_obj} that satisfies \Cref{ass: sup estimator second moment} \& \ref{ass: sup estimator_quality} with parameters $a_g$ and $M_g$ and \Cref{ass:stability} with $b_g> 1/2 + a_h - a_g/2$. 

We further assume either one of the following conditions holds:
\begin{enumerate}
    \item The estimators $\{g_i\}$ achieve a better rate: $a_g \in (a_h, 1]$;
    \item The nonparametric estimators $\{g_i\}$ achieves a better constant: $a_g = a_h < 1$ and $M_g < M_h$.
\end{enumerate}
Under these assumptions, we have selection consistency: 
\begin{equation*}
\lim _{n \rightarrow \infty} \mathbb{P}\left[{\rm RV}\left(\left\{g_i, i\in [n]\right\}, \xi\right)<{\rm RV}\left(\left\{h_i,i\in [n]\right\}, \xi\right)\right]
=1
\end{equation*}
where $\xi \geq 0$ is a fixed number. 
\end{theorem}

The proof of Theorem~\ref{th: main theorem} is presented in Appendix~\ref{app:proof of main theorem}. We will discuss key components of the proof in Section~\ref{section: technical results}.

\begin{remark}
   In Theorem~\ref{th: main theorem}, smaller $a_h$ and $a_g$ result in a milder stability requirement. We can plug in specific values to illustrate the interplay between convergence and stability rates:
\begin{itemize}
	\item For two sequences of inconsistent estimators, where $a_g = a_h = 0$, we only require $b_g, b_h > 1/2$.
	\item When $g_i$ converges to $f_0$ at a parametric rate ($a_g = 1$), we need $b_g > a_h$, which is stronger than the requirement of $b_h$: $b_h > (a_h+1)/2$.
	\item When both estimator sequences have the same rate $a_g = a_h = 1/2$, the stability requirements are also identical: $b_g,b_h > 3/4$.
\end{itemize}
\end{remark}

\begin{remark}
    (Any $\xi$ leading to consistency) In \Cref{th: main theorem} we stated that any fixed choice of $\xi$ would yield favorable model selection results. Intuitively, multiplying the $i^\xi$ factor enlarges both the signal and the noise by the same proportion. While the specific choice of $\xi$ has a meaningful impact on the finite-sample performance (\Cref{section: detection delay}), it does not alter the asymptotic consistency result under the current theoretical granularity.
    
    Moreover, for any choice of $\xi$, at any step $i$, approximately $i$ many samples are assigned a non-trivial weight when calculating wRV. For example, there are roughly $(1-2^{-1/\xi})i$ samples receiving a weight that is greater than $i^\xi/2$ (the weight of the most recent sample divided by $2$). However, this is not the case for the exponential weighting scheme $\exp(i\xi)$, which we believe does not yield similar results. Due to the rapid divergence of the exponential function, only a constant number of samples are assigned non-trivial weights. 
\end{remark}

\subsection{Sketch of the Analysis and Key Technical Ingredients}
\label{section: technical results}

The wRV statistics is an unbiased estimator of the estimators' cumulative population risk \eqref{eq: expectation of RV stat}. Theorem~\ref{th: main theorem} assumed a ``gap'' between the rolling population risks for $g_i$ and $h_i$. To ensure statistical consistency of wRV, we need to establish some concentration results of the sample RV statistics on their individual means. Specifically, the fluctuation of the sample quantities should be of a smaller order than the magnitude of the population risk gap. Consequently, the realized trajectories of ${\rm RV}(\{g_i\})$ and ${\rm RV}(\{h_i\})$ will eventually separate. To elucidate the concentration properties, we rewrite the difference between sample wRV and its population mean as follows. Define $u_i=i^{\xi}\left\{\hat{f}_{i-1}\left(X_i\right)-f_0\left(X_i\right)\right\}^2$:
\begin{equation}\label{concentration splitting}
    \begin{aligned}
        & \sum_{i=1}^n i^\xi \{ \hat f_{i-1}(X_i) - Y_i\}^2 - i^\xi\epsilon_i^2 - \sum_{i=1}^n i^\xi \mathbb{E}[\{\hat f_{i-1}(X) - f_0(X)\}^2]\\
    = & \sum_{i=1}^n u_i - 2i^\xi\epsilon_i \{\hat f_{i-1}(X_i) - f_0(X_i)\} - \sum_{i=1}^n \mathbb{E}[u_i]\\
    = & \left(\sum_{i=1}^n u_i - \mathbb{E}[u_i | F^{i-1}]\right) + \left(\sum_{i=1}^n \mathbb{E}[u_i|F^{i-1}] - \mathbb{E}[u_i]\right)  \\
    & \ -\left(2\sum_{i=1}^n i^\xi\epsilon_i \{\hat f_{i-1}(X_i) - f_0(X_i)\}\right)\,.
    \end{aligned}
\end{equation}

For each of the three terms, we can derive a concentration result. The first and the third terms are typical, and deriving their concentration does not need the stability condition \Cref{ass:stability}. For the third term, we have the following result.

\begin{lemma}
\label{puzzle 1}
Assume IID sampling scheme \Cref{ass:iid} and $\{\hat f_{i}\}$ satisfying \Cref{ass: sup estimator_quality} with a constant $a$, we have
\begin{equation*}
    \lim _{n \rightarrow \infty} \mathbb{P}\left(\left|\sum_{i=1}^n i^\xi\epsilon_i\left\{\hat{f}_{i-1}\left(X_i\right)-f_0\left(X_i\right)\right\}\right| \geq c_n\right)=0
\end{equation*}
for any positive $c_n$ satisfying $\lim_{n\rightarrow\infty} c_n^{-2}(n^{1 - a + 2\xi} \vee \log n)= 0$.
\end{lemma}
The proof of \Cref{puzzle 1} is based on a direct application of Chebyshev's inequality and controlling the variance. The detail is given in Appendix~\ref{app: proof of puzzle 2}. We can also derive the concentration properties of the first term using a similar argument. Very often it is a higher-order term compared with the other two.

\begin{lemma}
\label{puzzle 2}
Assume IID sampling scheme \Cref{ass:iid} and $\{\hat f_{i}\}$ satisfying \Cref{ass: sup estimator_quality} with a constant $a$, we have
\begin{equation*}
\lim _{n \rightarrow \infty} \mathbb{P}\left(\left|\sum_{i=1}^n u_i-\mathbb{E}\left[u_i \mid F^{i-1}\right]\right| \geq c_n\right)=0
\end{equation*}
for any positive $c_n$ such that $\lim_{n\rightarrow\infty} c_n^{-2}(n^{1 - 2a + 2\xi}\vee \log n\cdot \mathbb{1}(2a -2\xi = 1) \vee 1) = 0$.
\end{lemma}

The second term in \eqref{concentration splitting} is the most technically challenging part since we are no longer handling a simple martingale sequence.
We need the stability of the estimator to control its variability.
\begin{lemma}
\label{puzzle 3}
Assume IID sampling scheme \Cref{ass:iid} and $\{\hat f_{i}\}$ satisfying \Cref{ass: sup estimator_quality} \& \ref{ass:stability} with constants $a, b$, we have
\begin{equation}\label{eq: non_martingale_concentration}
\lim _{n \rightarrow \infty} \mathbb{P}\left(\left|\sum_{i=1}^n \mathbb{E}\left[u_i \mid F^{i-1}\right]-\mathbb{E}\left[u_i\right]\right| \geq c_n\right)=0
\end{equation}
for any positive sequence $c_n$ satisfying 
\begin{equation}\label{eq: non_martingale_cn}
\lim _{n \rightarrow \infty} c_n^{-2}\left(\log n \vee n^{3-a-2 b+2 \xi}\left(\log n \cdot \mathbb{1}\{3-a-2 b+2 \xi=1\} \vee 1\right)\right)=0.
\end{equation}
\end{lemma}
The proof of Lemma~\ref{puzzle 3} is presented in Appendix~\ref{app:proof of puzzle 3}. We also provide examples in which ${\rm Var}(\sum_{i=1}^n \mathbb{E}\left[u_i \mid F^{i-1}\right]-\mathbb{E}\left[u_i\right])$ can be explicitly characterized, which implies the requirement of $c_n$ in \eqref{eq: non_martingale_cn} cannot be further relaxed to achieve the desired concentration property \eqref{eq: non_martingale_concentration} in general. See \Cref{app: non_martingale_lemma_is_sharp} for the derivation and numerical verification. 

\begin{remark}
A typical $c_n$ that we are interested in is $c_n = n^{1-a+\xi}$---this is the order of $\sum_{i=1}^n\mathbb{E}[u_i]$ when $a\in[0,1)$. If we want the deviation $\left|\sum_{i=1}^n \mathbb{E}[u_i|F^{i-1}] - \mathbb{E}[u_i]\right|$ to be of a smaller order than the mean $\sum_{i=1}^n\mathbb{E}[u_i]$, we need the stability index $b$ to be greater than $(a+1)/2$ (stated in Theorem~\ref{th: main theorem}). As $a$ getting larger, the convergence rate of the considered estimator gets better, but the stability requirement becomes more stringent.   
\end{remark}

Combining the above lemmas with a union bound, we can derive the following concentration of sample wRV over its expectation, which then yields the consistency result claimed in \Cref{th: main theorem}, as detailed in Appendix~\ref{app:proof of main theorem}.
\begin{corollary}
\label{main theorem each sequence}
Let $\{\hat f_{i}\}$ be a given sequence of estimators. Assume that \Cref{ass:iid}, \ref{ass: sup estimator_quality} \& \ref{ass:stability} hold with constants $a, b$. We have
\begin{equation*}
 \lim _{n \rightarrow \infty} \mathbb{P}\left(\left|\sum_{i=1}^ni^\xi\left\{\hat{f}_{i-1}\left(X_i\right)-Y_i\right\}^2-i^\xi\epsilon_i^2-\sum_{i=1}^n i^\xi \mathbb{E}[r(\hat{f}_{i-1})]\right| \geq c_n\right)=0
\end{equation*}
for any $c_n$ satisfying 
\begin{equation*}
\label{clean cn scale}
    \lim _{n \rightarrow \infty} c_n^{-2}\left(\log n \vee n^{1-a+2\xi} \vee n^{3-a-2b+2\xi} \log n\right)=0\,.
\end{equation*}
\end{corollary}

\section{Consistency with Diverging Candidate Cardinality}\label{section: diverging cardinality}

So far, we have focused on comparing the performance of a fixed number of candidate estimator sequences specified by their hyperparameters $\{(\lambda_i^{(k)}: i \in \mathbb{N}^+),~k\in[K]\}$. It is also natural to explore a diverging number of model sequences in online learning settings. As more data is collected, researchers may refine or consider a broader range of hyperparameter sequences and explore more modeling possibilities (e.g., including higher-order interactions between interesting covariates). 

We will extend our theoretical results in \Cref{th: main theorem} to the case of diverging candidate cardinality, where the number of models, $K_i$, may increase with the sample size $i$.  It is also possible that $K_i$ decreases with $i$ in practice, corresponding to terminating certain inferior sequences. This scenario is statistically less challenging to establish consistency, as fewer models are under comparison.  At each step $i$, we denote the collection of all hyperparameters and their corresponding estimates as
\begin{equation*}
    \Lambda_i = \{\lambda_i^{(k)}, k\in[K_i]\}\quad \text{and}\quad G_i = \{\hat f_i^{(k)}, k\in[K_i]\}.
\end{equation*}

Similarly, at step $i+1$, we have $\Lambda_{i+1}$ and $G_{i+1}$, with cardinalities $K_{i+1} \geq K_i$.

In the increasing $K_i$ setting, the history trajectory of $\hat f^{(k)}_i$ may be unavailable if $k\notin [K_{i-1}]$ (assuming new models are assigned larger model indices $k$). When $K_{i+1} > K_i$, there are $K_{i+1} - K_i$ hyperparameter/estimator sequences whose definition only starts from step $i+1$. To facilitate formal discussion, we will retrospectively define their value before step $i+1$. We offer one option in \Cref{remark: history_definition}. The theoretical guarantees remain similar for other natural choices: (1) setting past estimates to $0$, (2) starting wRV calculation from the step they appear (rather than from $1$).

\begin{remark}\label{remark: history_definition}
For each element $\lambda_{i+1}^{(k)}$ in $\Lambda_{i+1}$, we define its representer in $\Lambda_{i}$ as $\lambda_{i-}^{(k)}$. Typically, the representer is the element in $\Lambda_i$ that is most ``similar" to $\lambda_{i+1}^{(k)}$. For $k \leq K_i$ the representer is simply $\lambda_{i-}^{(k)} = \lambda_{i}^{(k)}$. For new hyperparameters $\lambda_{i+1}^{(k)}$ with $k >K_i$, we find their representers in $\Lambda_i$ by minimizing a mathematical distance, or by identifying conceptual relationships between $\Lambda_i$ and $\Lambda_{i+1}$. For example, some elements in $\Lambda_{i+1}$ are a refinement of earlier hyperparameters. Once the representer mapping is established, we can map each element in $G_{i+1}$ to one in $G_i$, assigning all $\hat f_{i+1}^{(k)}, k \in [K_{i+1}]$ a previous value $\hat f_{i}^{(k)}$. This process can be traced back to $i=1$, ensuring all sequences of hyperparameters and estimators are well-defined, allowing us to continue our discussion.
\end{remark}

In the next two subsections, we outline the technical conditions required for consistency and then formally state our results.

\subsection{Technical Conditions for Consistency}

The proof structure for the general consistency result mirrors that in \Cref{section: technical results}. However, instead of finite-degree moment conditions, we require stronger light-tail conditions on the noise and scaled estimators, effectively imposing regularity conditions on all moments.

\begin{definition}(Sub-Weibull)
Let $L$ and $\theta$ be positive numbers. A random variable $X$ is defined as $(L, \theta)$-sub-Weibull if any of the following equivalent conditions holds:
\begin{enumerate}
    \item There exists a constant $C$ such that $\mathbb{P}\left(|X| \geq Lt\right) \leq C \exp(-t^\theta)$ for all $t>0$.
    \item There exists a constant $c$ such that $(\mathbb{E}[X^q])^{1/q} \leq c L q^{1 / \theta}$ for all $q \geq 1$.
\end{enumerate}
\end{definition}

The parameter $\theta$ characterizes the tail behavior of the random variable $X$. When $\theta = 2$, the sub-Weibull distribution coincides with the sub-Gaussian distribution, and $L$ corresponds to the Orlicz norm $\|X\|_{\psi_2}$. However, we consider the broader sub-Weibull family since sub-Gaussian (or exponential) random variables are not closed under multiplication: the product of two sub-Gaussian variables is not necessarily sub-Gaussian. For a more detailed discussion on sub-Weibull random variables, see \cite{vladimirova2020sub}.

We consider the following sub-Weibull versions of the moment conditions in \Cref{th: main theorem}. 

\begin{assumption}\label{ass: subW noise}
    The data points $\left\{\left(X_i, Y_i\right) \in\right.$ $\left.\mathbb{R}^{p+1}, i \in \mathbb{N}^+\right\}$ are IID samples from a common distribution $P(X, Y)$. The centered noise variable $\epsilon_1=Y_1-f_0\left(X_1\right)$ is $(K_\epsilon, \theta_\epsilon)$-sub-Weibull. 
\end{assumption}

\begin{assumption}\label{ass: subW estimator quality}
    Let $\{\hat f_i^{(k)}, k \in [K_i]\}$ be the $K_i$ estimators of $f_0$ after collecting $i$ samples. For sufficiently large $i > N_{\rm est}$, and each $k \in K_i$, the scaled error variable $i^{a_k / 2}\left(\hat{f}_i^{(k)}-f_0\right)(X)$ is $(K_{\rm est}, \theta_{\rm est})$-sub-Weibull.
%     \begin{equation}\label{eq: error subW}
% \mathbb{P}\left(i^{a_k/2}\left(\hat{f}_i^{(k)}-f_0\right)(X) \geq t K_{\rm est}\right) \leq \exp \left(-t^{\theta_{\rm est}}\right) \text{for each }k\in K_i.
%     \end{equation}
    The constants $a_k\in [0,1)$ may vary with the model index $k$. The constants $N_{\rm est}, K_{\rm est}$, $\theta_{\rm est} > 0$ are universal. 
    
    Additionally, for $i > N_{\rm est}$ and each $k\in [K_i]$
\begin{equation}\label{eq: exact ak}
\inf_{\ell \geq i} \ell^{a_k} \mathbb{E}[r_{\ell,k}] \geq M_k > 0,
\end{equation}
    with positive constants $M_k$ allowed to vary with $k$.
\end{assumption}

The condition \eqref{eq: exact ak} in \Cref{ass: subW estimator quality} guarantees that $a_k$ reflects the actual convergence rate---rather than a conservative lower bound. Without this line of restriction, one could take $a_k = 0$ and the rescaled errors $i^{a_k / 2}\left(\hat{f}_i^{(k)}-f_0\right)(X)$ are trivially sub-Weibull for consistent estimators.

\begin{assumption}\label{ass: subW stability}
    Let $\{\hat f_i^{(k)}, k \in [K_i]\}$ be the $K_i$ estimators of $f_0$ after collecting $i$ samples. For any $i$ and $j \in [i]$, the scaled stability variables $i^{b_k}
\nabla_j \hat{f}_i(X)$ are $(K_{\rm stab}, \theta_{\rm stab})$-sub-Weibull for each $k\in[K_i]$. The constants $b_k > 0$ may vary with model index $k$, while $K_{\rm stab}$, $\theta_{\rm stab} > 0$ are universal.
\end{assumption}

In Assumption~\ref{ass: subW estimator quality} and \ref{ass: subW stability}, the estimator quality parameter $a_k$ and stability parameter $b_k$ can vary between sequences. However, we assumed the tail behaviors of the scaled random variables are depicted by two sets of universal constants $(K_{\rm est},\theta_{\rm est})$ and $(K_{\rm stab},\theta_{\rm stab})$. This assumption could be relaxed to allow $k$-dependent sub-Weibull parameters plus a uniform bound over them. For presentation simplicity, we have retained the $k$-agnostic conditions. 

\subsection{Consistency with a Diverging Number of Candidate Sequences}

Under the new light-tail conditions, wRV can identify the model with superior predictability against all alternatives. We present the formal results below, with comments on the technical conditions in the following remarks.

\begin{theorem}\label{th: diverging cardinality}
    Assume the learning setting in Assumption~\ref{ass: subW noise}. Let $\left\{h_i^{(k)}, i \in[n]\right\}$, $k \in\left[K_n\right]$ be sequences of estimators that satisfy Assumption~\ref{ass: subW estimator quality} and \ref{ass: subW stability} with constants $\{(a_k, b_k)\}$.
    
    We further assume:
    \begin{itemize}
        \item (U1) The estimator quality parameters satisfy
    \begin{equation*}
        \lim_{n\rightarrow\infty} \left(\frac{\log n}{\log \log n}\right)\min _{k \in\left[K_n\right]}\left(b_k-a_k / 2-1 / 2\right) \wedge\left(1-a_k\right) = \infty.
    \end{equation*}
    Moreover, we require $\log K_n \leq C\log n$ for a positive constant $C$.
    \item (U2) There exists a positive number $ M_{h} > 0$ such that
    \begin{equation*}
        \liminf_{n\rightarrow \infty} n^{a_{h,n}} \min_{k \in K_n} \mathbb{E}[r(h^{(k)}_n)] = M_{h},
    \end{equation*}
    where $a_{h,n} = \max_{k\in K_n} a_k$.
    \end{itemize}
    Let $\left\{g_i\right\}$ be the sequence of superior estimators \eqref{eq:seq_selec_obj}, satisfying \Cref{ass: sup estimator second moment} and \ref{ass: sup estimator_quality} with parameter $a = a_g$ and $M = M_g$, as well as Assumption~\ref{ass:stability} with $b_g>1 / 2+ a_h-a_g / 2$, where $a_h := \limsup_n a_{h,n}$.
    
    Then, for either of the following scenarios:
    \begin{enumerate}
        \item The estimators $\left\{g_i\right\}$ achieve a better rate: $a_h < a_g \leq 1$;
        \item The estimators $\left\{g_i\right\}$ achieve a better constant: $a_h = a_g<1$ and $M_g < M_h$, 
    \end{enumerate}
    we have
\begin{equation*}
    \lim _{n \rightarrow \infty} \mathbb{P}\left[{\rm RV}\left(\left\{g_i, i\in[n]\right\}, \xi\right)< \inf_{k\in [K_n]}{\rm RV}\left(\left\{h_i^{(k)}, i\in[n]\right\}, \xi\right)\right]=1,
\end{equation*}
where $\xi \geq 0$ is a fixed number. 
\end{theorem}

The proof of Theorem~\ref{th: diverging cardinality} is presented in Appendix~\ref{app: infinite model}.
\begin{remark}(Assumption interpretation)
   Condition (U1) allows both $b_k-a_k / 2-1 / 2$ and $\left(1-a_k\right)$ to converge slowly to $0$. These conditions generalize the requirements $b_h-a_h / 2-1 / 2 > 0$ and $1-a_h > 0$ in \Cref{th: main theorem}. When both $a_h$ and $a_g$ equal $1$, we are comparing two parametric models, and it is known that CV would not choose the favorable model with a probability approaching $1$. Our framework accommodates $a_g = 1$ while allowing $a_h$ to converge to the parametric threshold. Condition (U2) imposes a uniform limit for all the sub-optimal estimator sequences, essentially requiring that the excess risk of all alternative models is no less than $M_hn^{-a_h}$.
\end{remark}

\begin{remark} (Regarding model number $K_n$) In \Cref{th: diverging cardinality}, we required the number of candidate models to satisfy $\log K_n \lesssim \log n$. Given the current literature works with light-tail assumptions, one may expect the number of items under comparison (here $K_n$) to scale like $\log K_n \sim n^{\alpha}$ for some $\alpha > 0$ when assuming similar exponential light-tail conditions. This is true if we further require $\liminf_n\min _{k \in\left[K_n\right]}\left(b_k-a_k / 2-1 / 2\right) \wedge\left(1-a_k\right)>0$. However, the current (U1) condition means to explore the regime where the convergence rates and stability conditions approach the boundary, which can afford at most polynomially many $K_n$.
\end{remark}

\begin{remark}
    The problem of maintaining a dynamic list of candidate models has been considered in the literature of online functional data analysis. For example, \cite{fang2023online} proposed to maintain a list of candidate bandwidths obtained as of weighted average of past candidates and the current candidates. At a high level, this is similar to our proposal of linking a new candidate to the closest one in the previous time step. However, there are two important distinctions between our framework and the setting in \cite{fang2023online}. First, the tuning parameter selection in \cite{fang2023online} is partially model-based, as it requires pilot estimates of functionals of unknown true and correctly specified model parameters. In contrast, our method is model-free and only requires a loss function to evaluate the predictive quality of the fitted model.  Second, the averaging approach taken in \cite{fang2023online} requires the candidates to be objects in the same vector space, such as different bandwidths.  In contrast, our method is much more flexible as it works for arbitrary collections of model fitting methods. For example, our method can compare local polynomial smoothers with basis expansion estimators, and can even incorporate black-box estimators.
\end{remark}

\section{Stability of Estimators}
\label{sectoin:example stability of estimators}
The consistency of wRV model selection relies on three interpretable conditions: IID sampling, MSE convergence, and estimator stability. While the convergence rates of estimators have been extensively studied in the literature, the stability condition is relatively new and deserves more discussion. 
In this section, we evaluate the stability rates for multiple prototypical estimators, providing algorithm-specific conditions for the assumptions in Theorem~\ref{th: main theorem} to hold.

\subsection{Batch Sieve Estimators}
\label{section: batch sieve stability}
The concept of stability applies to both online and batch estimators. We begin our discussion with a batch sieve estimator, which is simpler and easier to understand at first reading. 

Recall the basic nonparametric regression model for IID sample $(X_i,Y_i)$ on $[0,1]\times\mathbb R$:
\begin{equation*}
    Y_i = f_0(X_i) + \epsilon_i,\ i\in [n]\,.
\end{equation*}

Assume that there is an orthonormal basis of $L_2(P_X)$, $\{\phi_k\}$, such that $f_0 = \sum_{k=1}^\infty \beta_k\phi_k$. We assume the coefficients satisfying a Sobolev ellipsoid \eqref{eq: sobolev ellipsoid} condition that $\sum_{k=1}^\infty k^{2s}\beta_k^2 \leq Q^2$
for some $s \geq 1, Q>0$. 

Consider estimators using the first $J = J_n$ functions in $\{\phi_k\}$:
\begin{equation}\label{eq: simple sieve 1}
    \hat f = \sum_{k=1}^{J} \hat\beta_k \phi_k\,,
\end{equation}
where the regression coefficients are calculated from the sample as:
\begin{equation}\label{eq: simple sieve 2}
    \hat \beta_k = n^{-1}\sum_{i=1}^n Y_i\phi_k(X_i)\,.
\end{equation}

This estimator is also called the Projection Estimator, whose estimation quality is studied under evenly-spaced $X_i$ in \cite{introtononpara}. We consider using $J = n^{\alpha}$ many basis functions for some $\alpha \geq 0$ and quantify its stability.

Let $\hat f(\cdot) = \hat f_n(\ \cdot\ ; \bm{Z}_n)$, $\hat f^\prime = \hat f_n(\ \cdot\ ;\bm{Z}_n^j)$ be two estimates but are different at the $j$-th training sample. For any $j\in[n]$,
    \begin{align*}
        \nabla_j \hat{f}_n(X) & =\hat{f}\left(X\right)-\hat{f}^{\prime}\left(X\right)\\
        & = n^{-1}\sum_{k=1}^{J_n}\left\{Y_j \phi_k\left(X_j\right)-Y_j^{\prime} \phi_k\left(X_j^{\prime}\right)\right\} \phi_k\left(X\right)\,.
    \end{align*}
Therefore, 
\begin{equation}\label{eq: stability batch sieve}
   \begin{aligned}
       \mathbb{E}\left[\left\{\nabla_j \hat{f}_n(X)\right\}^2 \mid F^j\right] & =n^{-2} \mathbb{E}\left[\left(\sum_{k=1}^{J_n}\left\{Y_j \phi_k\left(X_j\right)-Y_j^{\prime} \phi_k\left(X_j^{\prime}\right)\right\} \phi_k(X)\right)^2 \mid F^j\right] \\
       & \stackrel{(I)}{=} n^{-2}\sum_{k=1}^{J_n}\left\{Y_j \phi_k\left(X_j\right)-Y_j^{\prime} \phi_k\left(X_j^{\prime}\right)\right\}^2 \stackrel{(II)}{\lesssim}  n^{\alpha-2}\text{ almost surely}\,,
\end{aligned} 
\end{equation}
where in step $(I)$ we used $\{\phi_k\}$ is orthonormal: $\mathbb{E}[\phi_l(X)\phi_k(X)] = \delta_{lk}$, and in step $(II)$ we assumed the outcome variable and the basis functions are uniformly bounded. We summarize the stability results established above as follows:
\begin{theorem}\label{th: stability batch sieve}
    Suppose we observe an IID sample $\{(X_i,Y_i), i\in[n]\}$ satisfying Assumption~\ref{ass:iid}. We further assume the outcome $Y_i$ and basis functions $\{\phi_k\}$ are uniformly bounded by $L>0$. Consider orthonormal basis $\{\phi_k\}$ w.r.t. $P_X$: $\mathbb{E}\left[\phi_l(X) \phi_k(X)\right]=\delta_{l k}$. Then for the estimator defined in \eqref{eq: simple sieve 1} and \eqref{eq: simple sieve 2} with $J = n^{\alpha}$, we have for any $j\in [n]$:
    \begin{equation*}
n ^{2-\alpha} \mathbb{E}\left[\left\{\nabla_j \hat{f}_{n}(X)\right\}^2 \mid F^j\right] \leq 4L^4 \text { almost surely }.
\end{equation*}
\end{theorem}

When $\phi_k$ is an orthonormal basis under $P_X$, the MSE of batch sieve estimators can be decomposed as (detail presented in \Cref{app: rigorous_discussion} )
\begin{equation}\label{eq: convergence_batch_sieve}
         \mathbb{E}[r(\hat f)]
      \lesssim J/n + J^{-2s}Q^2.
\end{equation}
In the case when $\alpha < (2s+1)^{-1}$, the bias term $J^{-2s}Q^2$ would dominate the variance term $J/n$ and the convergence rate of $\hat f$ will be $n^{-2s\alpha}$. When we use significantly more basis than needed ($\alpha > (2s+1)^{-1}$), the variance term $J/n$ will dominate and the convergence rate behaves like $n^{\alpha-1}$. To obtain rate-optimal estimators, we need to balance the $J/n$ and the $J^{-2s}Q^2$ terms so that they are of the same order. This leads to $J \sim n^{1/(2s+1)}$, with an MSE convergence rate $n^{-2s/(2s+1)}$.

At this point, both the convergence and stability rates are explicitly related to the diverging speed of the number of basis functions. Suppose we are estimating $f_0$ using repeatedly refitted batch sieve estimators (just for theoretical interest) and we are comparing the following two sequences of candidates: $\{g_i\}$ uses the first $i^{1/(2s+1)}$ basis functions and $\{h_i\}$ uses $i^{\alpha}$ of them with  $\alpha\neq 1/(2s+1)$. We know $g_i$ is the better sequence because it implements the correct order of model capacity, so a good model selection procedure is expected to choose it. In view of Theorem~\ref{th: main theorem}, we require the following stability conditions to ensure consistent selection:
\begin{equation}\label{eq: recap stability}
    \begin{aligned}
    & b_h > 1/2 + a_h/2\\
        & b_g>1 / 2+a_h-a_g / 2
    \end{aligned}
\end{equation}
where $b_h = 1-\alpha/2, a_h = (1-\alpha) \wedge 2s\alpha, b_g = 1-(4s+2)^{-1}$ and $a_g = 2s(2s+1)^{-1}$. These make \eqref{eq: recap stability} equivalent to requiring $\alpha < (2s+1)^{-1}$.

That is, we can ensure the wRV procedure consistently selects $\{g_i\}$ over $\{h_i\}$ if the latter uses fewer basis functions than the optimal order. When $\alpha > (2s+1)^{-1}$, we need the $b_h > 1 - \alpha/2$ to establish the asymptotic consistency, but according to \Cref{th: stability batch sieve}, $b_h = 1 - \alpha/2$, which lies at the boundary of the open interval. Although the simulation results show wRV with probability $1$ can rule out the inferior sequence (\Cref{fig: batch_sieve_consistency}), the current analysis based on stability conditions cannot theoretically reject this particular type of over-fitting estimator. Given the individual analysis of the three terms in \Cref{concentration splitting} is optimal and the stability bound is exact, we conjecture there are some intricate cancellations between the terms in \eqref{concentration splitting} that contribute to the numerical observation. 

In addition to the above hard setting when we compare two consistent estimators, if we instead assume $\{h_i\}$ uses $i^{\alpha}$ many basis functions but is an inconsistent estimator of $f_0$, then conditions in \eqref{eq: recap stability} are much easier to satisfy and wRV can always pick the better sequence $\{g_i\}$ so long as $\alpha < 1$. 
% We have more discussion regarding the evenly-spaced grid requirement in \Cref{app: rigorous_discussion} to make the discussion above more rigorous. 

\subsection{Online Parametric SGD}
In Section~\ref{section: batch sieve stability} we discussed the stability properties of batch sieve estimators. Now we switch back to the online setting. In this section, we consider a linear SGD estimator with the update formula
\begin{equation}
\label{eq: parametric iteration}
\begin{aligned}
       \beta_0 & = 0 \in\mathbb{R}^p,\\
       \beta_i & = \beta_{i-1} + \gamma(Y_i - X_i^\top \beta_{i-1})X_i\,,
\end{aligned}    
\end{equation}
where $(X_i,Y_i)\in\mathbb R^{p+1}$ and $\gamma > 0$ is a fixed learning rate.
We consider the estimator $\hat f_i$ using the averaged parameter:
\begin{equation*}
    \hat f_{i}(x) = x^\top \bigg(i^{-1} \sum_{k = 1}^i \beta_k\bigg)\,.
\end{equation*}
This estimator has been shown to achieve the parametric minimax rate in \cite{marteau2019globally} when the truth is linear. Stability results under strongly convex losses are known in \cite{hardt2016train}, but our treatment can handle the cases when the loss is ``on average strongly convex'' (\cite{hardt2016train} does not cover the simple, unpenalized regression setting). The following results may be of interest itself, whose proof techniques are also a primer of the nonparametric results in Section~\ref{section: sieve SGD stability}.

\begin{theorem}
\label{th: stability of parametric sgd}
Assume we observe IID samples $\{X_i,Y_i\}$ satisfying \Cref{ass:iid}. We further assume that $\|X_i\|\leq R$ and $|Y_i - X_{i-1}^\top\beta_i| \leq M$ almost surely. The minimal eigenvalue $\lambda_{min}(\mathbb{E}[XX^\top]) \geq \underline\lambda>0$. When the learning rate $\gamma \leq R^{-1}$, we have, for any $j\in [i]$,
    \begin{equation*}
i^2\mathbb{E}\left[\left\{\nabla_j \hat f_{i}\left(X\right)\right\}^2 \mid F^j\right] \leq C(\gamma, \underline\lambda, M, R)\quad\text{ almost surely },
\end{equation*}
where $C(\gamma,\underline\lambda, M, R)>0$ is a constant depending on $\gamma, \underline\lambda, M, R$.
\end{theorem}

The proof is given in Appendix~\ref{app:stability of SGD}. Theorem~\ref{th: stability of parametric sgd} states that parametric SGD has a strong mode of stability, in the language of \Cref{ass:stability}, whose rate $b = 1$. Our result is also sharp: it is direct to verify that the stability rate cannot be better than $1$ under some simple scenarios ($p=1$, $X_i = 1$).

\subsection{Nonparametric SGD: Sieve-type}
\label{section: sieve SGD stability}
In this section, we investigate the stability of a genuine online nonparametric estimator: the sieve-SGD \eqref{eq: sieve sgd update}, which combines the approaches of batch sieve estimators and online parametric SGD estimators. Recall that the estimator $\hat f_i$ is a weighted linear combination of $J_i$ pre-specified basis functions
\begin{equation}
\label{eq: recall seive sgd}
    \hat f_i(x) = \sum_{k=1}^{J_i} \hat \beta_i[k]\phi_k(x),
\end{equation}
where the weights $\hat\beta_i\in\mathbb{R}^{J_i}$ is determined using data. 

Denote 
\begin{equation*}
\boldsymbol\phi_i = \phi(X_i;J_i) = (\phi_1(X_i),...,\phi_{J_i}(X_i))^\top\in\mathbb{R}^{J_i}.
\end{equation*}
%%%%%%
Very similar to the parametric correspondence, the coefficient update rule for sieve-SGD is:
\begin{equation}
\label{eq: sieve coeff update rule}
  \begin{aligned}
& \beta_0=0 \in \mathbb{R} \\
& \beta_i = \beta_{i-1}^\wedge  + \gamma_i ( Y_i - \boldsymbol\phi_i^\top \beta_{i-1}^\wedge) D_i\boldsymbol\phi_i\,,
\end{aligned}  
\end{equation}
where the $\beta^\wedge$ notation means embedding a vector $\beta$ into a higher dimension Euclidean space---so that its length is compatible with the other components--- and padding the extra dimensions with $0$. For example, in the second line of \eqref{eq: sieve coeff update rule}, $\beta_{i-1}^\wedge = (\beta_{i-1}^\top, 0,...,0)^\top\in\mathbb{R}^{J_i}$, which is of the same dimension as $\boldsymbol\phi_i$. The matrix $D_i$ is a diagonal with elements $(1^{-2\omega},...,J_i^{-2\omega})$. The final coefficient vector in \eqref{eq: recall seive sgd} is just the average of the trajectory:
\begin{equation*}
    \hat \beta_i = i^{-1}\sum_{l=1}^i \beta_l^{\wedge} \in \mathbb{R}^{J_i}\,.
\end{equation*}

Now we are ready to state our main result regarding the stability of sieve-SGD:
\begin{theorem}
\label{th: stability sieve with shrinkage}
Assume we observe IID samples $\{X_i,Y_i\}$ satisfying \Cref{ass:iid} and consider bounded orthonormal basis $\{\phi_k\}$ w.r.t. some measure $\nu$: $\int \phi_l(x) \phi_k(x) d \nu(x)=\delta_{l k}$. We further assume that $|Y_i - \boldsymbol\phi_{i-1}^\top\beta_i| \leq M$ and the density (with respect to $\nu$) of $X$ is bounded from above and below.
When $\gamma_i = i^{-\zeta}$, $J_i = i^\zeta$, for some $\zeta \geq 0, \omega > 1/2$, we have for any $j \in [i]$
\begin{equation*}
i^{2-4\omega\zeta} \mathbb{E}\left[\left\{\nabla_j \hat{f}_{i}\left(X\right)\right\}^2 \mid F^j\right] \leq C\quad \text { almost surely},
\end{equation*}
for some constant $C > 0$.
\end{theorem}
The proof of Theorem~\ref{th: stability sieve with shrinkage} is given in Appendix~\ref{app:stability of Sieve SGD}. 

Suppose we have multiple sieve-SGD estimators that assume a different degree of smoothness $s^{(k)}$, $k\in [K]$ and take the hyperparameter $\zeta$ in Theorem~\ref{th: stability sieve with shrinkage} as $\zeta^{(k)} = (2s^{(k)}+1)^{-1}$. Denote their convergence rates as $a^{(1)},...,a^{(K)}$. If one of the models properly specifies the smoothness $s^*$---more accurately, it specifies the largest possible $s$ such that $f_0\in W(s)$, it will achieve an MSE convergence rate $2s^*/(2s^*+1)$. Theorem~\ref{th: stability sieve with shrinkage} tells us that this  estimator will have a smaller wRV value against other trajectories that satisfy (using \eqref{eq: recap stability}):
    \begin{align*}
        1 - 2\omega \zeta^{(k)}  & > 1/2 + a^{(k)}/2,\\
        1 - 2\omega (2s^* + 1)^{-1} & > 1/2 + a^{(k)}/2 - s^*/(2s^*+1).
    \end{align*}
Rearranging the terms and let $\omega$ approach $1/2$, the above condition becomes:
\begin{equation}\label{eq: clean range sieve sgd}
    \begin{aligned}
        & a^{(k)}  < \frac{2s^*-1/2}{2s^*+1},\\
        & a^{(k)} + \frac{2}{2s^{(k)} + 1} < 1\ .
    \end{aligned}
\end{equation}
According to \eqref{eq: clean range sieve sgd}, any alternative sieve-SGD will be excluded by wRV if it: (i) has a gap in convergence rate against the better value $2s^*/(2s^*+1)$ (first line of \eqref{eq: clean range sieve sgd}) and (ii) uses a large enough $s^{(k)}$ (second line of \eqref{eq: clean range sieve sgd}). This is similar to the results we derived in Section~\ref{section: batch sieve stability} for batch sieve estimators. We can theoretically rule out under-fitting sequences, and simulation studies provide complementary evidence for ruling out over-fitting sequences (\Cref{section: numerical}). We conjecture the $i^{2-4\omega\zeta}$ factor in Theorem~\ref{th: stability sieve with shrinkage} can possibly be improved using different conditions or arguments, and will pursue this in future work.

\begin{remark}
    The uniform bound on $X$ and/or $Y$ in \Cref{th: stability batch sieve}-\ref{th: stability sieve with shrinkage} cannot be dropped when establishing the uniform bounds on the stability terms. This is more evident from examining \eqref{eq: stability batch sieve} where the calculation is very straightforward. However, as we discussed in \Cref{remark: holder tradeoff}, the almost sure bound on stability is not always necessary. When the estimator converges in a stronger mode to $f_0$, such as in $\|\cdot\|_{\infty}$, it can be relaxed to $\mathbb{E}\left[\left\{\nabla_j \hat{f}_i(X)\right\}^2\right] \lesssim i^{-2b}$. To establish the bounds on unconditional expectations of the stability terms, we can follow the same treatment but under milder moment assumptions such as $\mathbb{E}[\|X\|^2] \leq R$ and $\mathbb{E}[(Y - f_0(X))^2 \mid X] \leq M$.
\end{remark}

\section{Numerical Examples}
\label{section: numerical}

We consider two settings for the simulation studies: a univariate case with a theoretically more tractable true regression function and a more realistic $10$-dimensional setting. 

\begin{example}[One-dimensional nonparametric regression]\label{exa:1}
We generate IID samples as follows:
\begin{align}
X_i & \sim \mathcal{U}([0,1]),\nonumber\\
       \epsilon_i & \sim \mathcal{N}(0, 0.5^2),\nonumber\\
       Y_i & = \sum_{k=1}^{30} k^{-2.5} \cos \left((k-1) \pi X_i\right)+\epsilon_i,\label{eq:example 1 setting}
\end{align}

where $\mathcal{U} (\mathcal{N}$) indicates a uniform (normal) distribution.
\end{example} 
We compare four sieve-SGD models trained with different smoothness parameters. It is expected that the one with the correctly specified smoothness would have the lowest estimation error and should be selected by our proposed wRV procedure. In the notation of \eqref{eq: sieve sgd update}, we set $\omega = 0.51$, $\phi_k(x) = \cos \left((k-1) \pi x \right)$, $J_i = i^{1/(2s+1)}$ and $\gamma_i = 0.1\cdot i^{-1/(2s+1)}$. The smoothness parameter $s$ varies between the four models $\in\{1,2,3,4\}$. Note that the true regression function in \eqref{eq:example 1 setting} and the basis functions of sieve-SGD are both cosine functions. Our main interest is whether wRV can pick $s = 2$ as the better model. Data is generated and processed in an online fashion: each mini-bath contains $100$ samples. When revealed to the learners, the samples are processed one by one (i.e. taking $100$ stochastic gradient descent steps within each mini-batch).

\begin{figure}[htbp!]
    \centering
    \includegraphics[width = \textwidth]{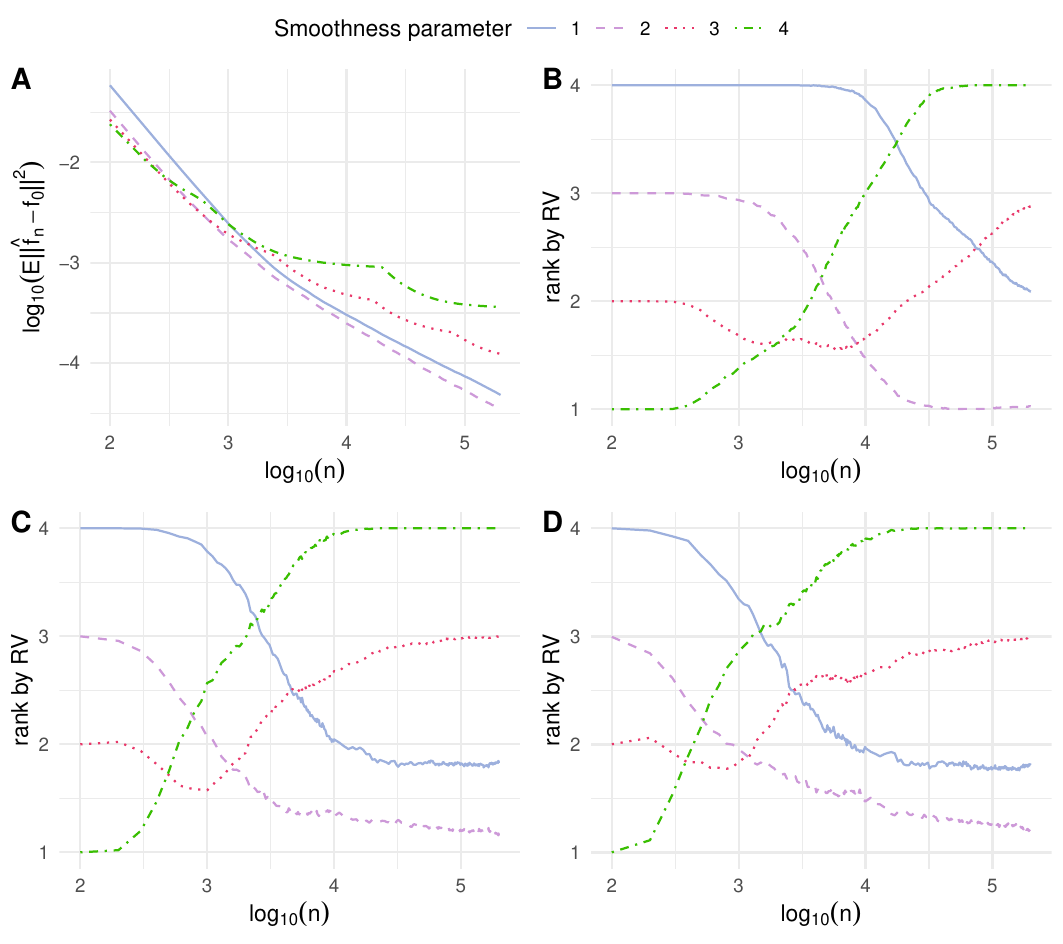}
    \caption{Model selection results of \Cref{exa:1}. (A) the true MSE of the candidate models; (B)-(D), average rankings of the models at different sample sizes over 500 repetitions, according to wRV, with weighting exponent $\xi = 0, 1, 2$, respectively.}
    \label{fig:theoretical simulation}
\end{figure}

In Figure~\ref{fig:theoretical simulation}, we present averaged simulation results from 500 repetitions. Panel A shows the true distance between the estimators and the underlying regression function, which is a piece of information not available in practice. Panels B-D show the average ranking of the four models, based on rolling statistics wRV with different weighting exponent $\xi\in\{0,1,2\}$. A smaller ranking value corresponds to a model preferred by wRV. In \Cref{app: extra numerical}, we also provide supplementary results using larger $\xi \in \{3,4\}$ and extending the sample size to $10^6$ to demonstrate selection consistency for $\xi > 0$. A larger $\xi$ can decrease the detection delay as discussed in \Cref{section: detection delay}. Comparing \Cref{fig:theoretical simulation} and \Cref{fig: theoretical larger xi}, we can see $\xi \in \{2,3,4\}$ gives very similar choices (measured by the location of the trajectory lines' crosses), slightly different from $\xi = 1$. We do not recommend using $\xi = 0$ (\Cref{fig:theoretical simulation}, B) mainly due to its insensitivity to the change of model ranking.

After processing the first $100$ samples ($\log_{10}(n) = 2$), Panel A tells us the best initial estimate corresponds to $s=4$, followed very closely by $s=3$ and $2$. The corresponding wRV statistics in Panels B, C, and D at $\log_{10} n = 2$ also rank $s=4$ as the best. However, as more data comes in, the models $s = 3 $ and $s=4$ start to be outperformed by the other two due to over-smoothing, as reflected in Panel A. They do not increase the number of basis functions fast enough and miss a better trade-off between the estimation and approximation errors. 

Asymptotically, we would expect $s=2$ to give the lowest estimation in this specific simulation setting since it properly specifies the smoothness level. Indeed, it eventually becomes the best-performing model among the four after processing $10^3$ samples (Panel A). The unweighted RV, unfortunately, cannot adjust its outdated choice until $10^4$ samples are collected. But the weighted version has a much faster transition: they start to consistently pick the correct model before $n \sim 2\times 10^3$. Although they are all asymptotically consistent methods, the finite sample performance can be significantly improved by implementing a diverging weight even in cases as simple as univariate regression.

\begin{example}[Multivariate nonparametric regression]\label{exa:2}
We also include a set of simulation results of multivariate nonparametric online estimation. In this setting, we have
\begin{align*}
X_i & \sim \mathcal{U}([0,1]^{10}), \nonumber\\
 \epsilon_i &\sim \mathcal{N}(0,2^2), \nonumber\\
 f_0(X_i)&=\sum_{j \text { is odd }}\left\{ 0.5-\left|X_i[j]-0.5\right|\right\}+\sum_{j \text { is even }} \exp \left(-X_i[j]\right), \nonumber\\
  Y_i&=f_0\left(X_i\right)+\epsilon_i\,.
\end{align*}
\end{example}

\begin{figure}[htbp!]
    \centering
    \includegraphics[width = \textwidth]{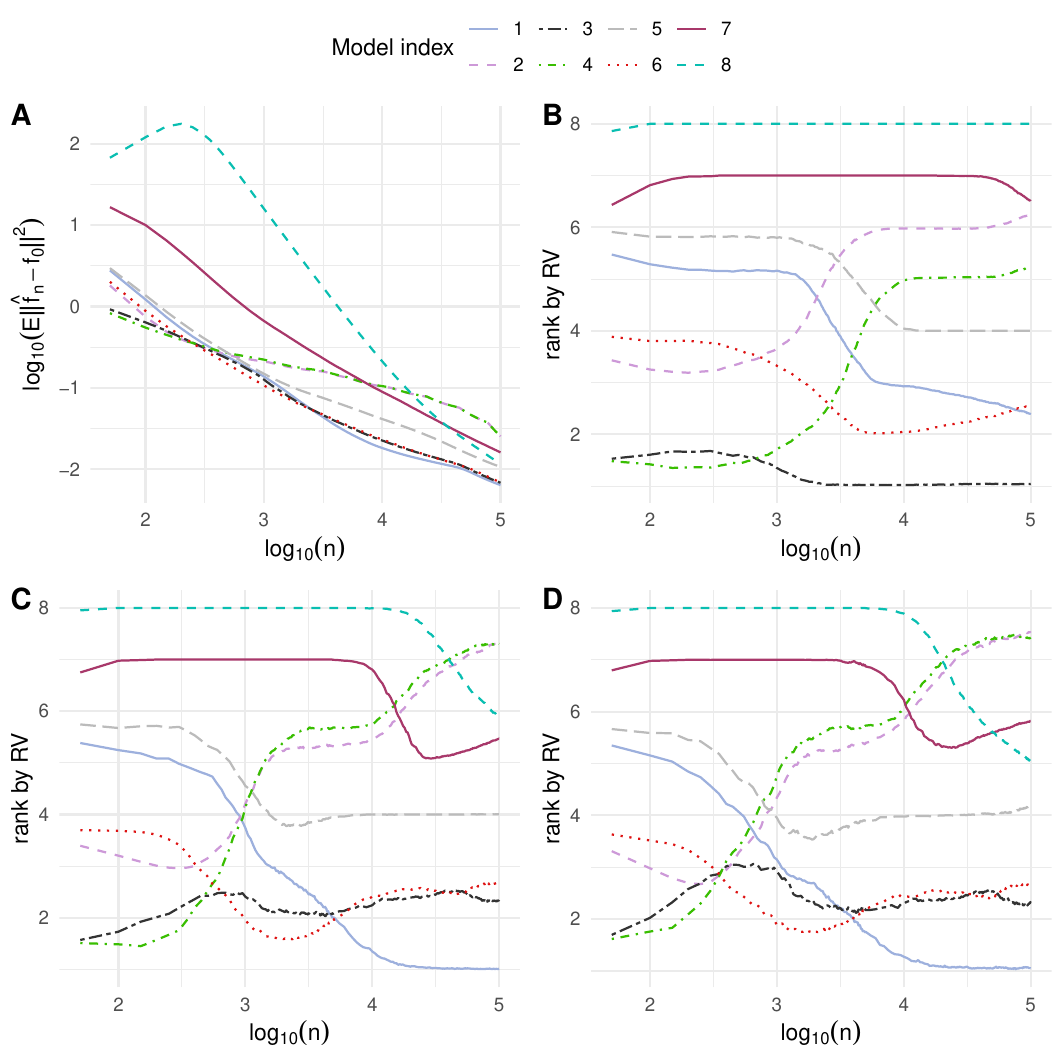}
    \caption{Model selection results of Example 2. (A) the true MSE of $8$ models; (B)-(D), average ranking of the models at different sample size, according to wRV. The RV weighting exponents are $\xi = 0, 1, 2$, respectively. Averaged over 500 repeats.}
    \label{fig:realistic simulation}
\end{figure}

We compare eight different sieve-SGD estimators that take different combinations of (i) assumed smoothness level $s\in\{1,2\}$; (ii) initial step size $A\in\{0.1,1\}$ and (iii) the initial number of basis functions $B\in\{2,8\}$. The definitions of these hyperparameters are given in \eqref{eq: step size and basis number}. The model index and hyperparameter combination correspondence are presented in \Cref{tab:parameter example 2}. 

The basis functions we use are products of univariate cosine functions: $\prod_{m=1}^{10} \cos((l[m]-1)\pi x[m])$ for some $l\in(\mathbb{N}^{+})^{10}$. In the univariate feature case, people almost always use the cosine basis of lower frequency first. However, in the multivariate case, there is no such ordering of the product basis functions. We reorder the multivariate cosine functions based on the product magnitude ($=\prod_{m=1}^{10}l[m]$) of the index vector $l$, in increasing order. This way of reordering multivariate sieve basis functions can lead to rate-optimal estimators in certain tensor product Sobolev spaces \citep{zhang2023regression}.

The true MSE (Panel A) and model selection frequencies (Panels B-D) for this example are presented in Figure~\ref{fig:realistic simulation}. According to Panel A, models 3 and 4 are better when the sample size is less than $250$. As more samples are being processed, models 1, 3, and 6 perform better. 

According to the unweighted RV statistics (Panel B), model 3 is always one of the best-performing ones (second place when $n < 10^{2.7}$), consistent with the MSE in Panel A. The wRV statistics also captured the phenomenon that model 4 is better when the sample size is small but falls behind when the sample size is larger. 

In Panel A, a close observation reveals that Model 1 becomes the best performer when $n>10^{3.5}$. Unfortunately, this phenomenon is not captured by the unweighted RV in Panel B, where Model 1 does not frequently secure the first place. This issue is addressed by wRV in Panels C and D. Model 1 is identified as the best choice by wRV when the sample size exceeds $10^{3.7}$, and this selection remains consistent when $n\ge 10^4$. Models 7 and 8, which use a large initial learning rate and initial basis function number, exhibit larger MSE at the early stages—Model 8 even shows an irregular increasing error. However, after a rapid improvement phase, they eventually surpassed Models 2 and 4, which initially performed well. This dynamic is scarcely captured by the unweighted RV within the sample size range considered here but is better reflected by the wRV. 

\section{Discussion}\label{section: discussion}
\begin{figure}[htbp!]
    \centering
    \includegraphics[width = \textwidth]{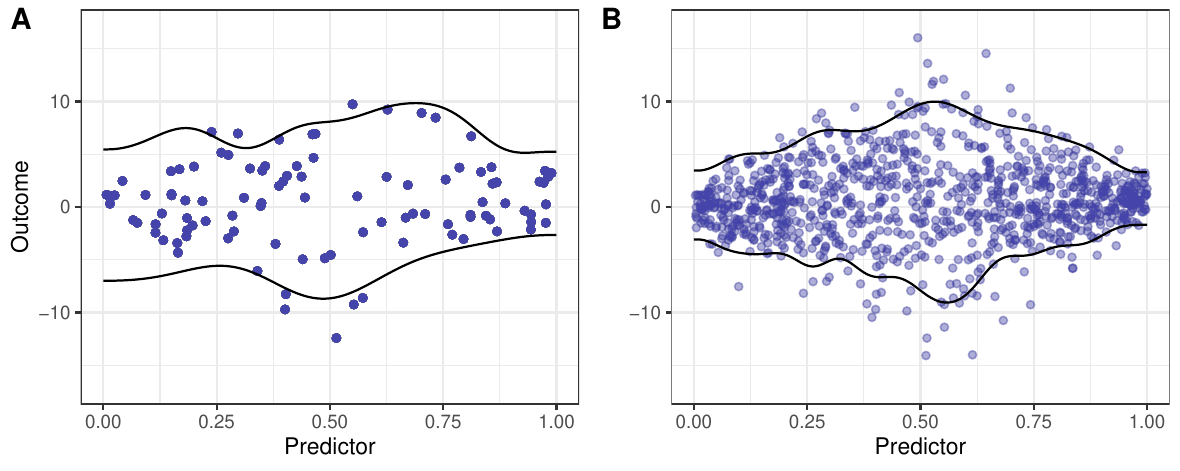}
    \caption{An example of applying wRV to quantile regression. The blue points represent training samples. (A) training sample size = 100. (B) training size increases to $10^3$, the blue points in (A) are a strict subset of those shown in (B). The estimators are trained in an online fashion using the same trajectory. The estimated coverages between the quantiles using a testing sample size $=10^4$ are (A) $91.7\%$ and (B) $90.6\%$.}
    \label{fig:quantile regression}
\end{figure}

The rolling validation method provides a computationally efficient procedure for adaptive online regression, with a theoretically justifiable performance guarantee under stability conditions. While we analyze our method in the regression setting, it can also be extended to other nonparametric estimation tasks. For example, one can replace the squared loss with the ``pinball'' loss 
\begin{equation*}
l(f ; X, Y):= \begin{cases}\alpha(Y-f(X)) & \text { if } Y>f(X) \\ (1-\alpha)(f(X)-Y) & \text { otherwise }\end{cases}
\end{equation*}
with $\alpha\in (0,1)$ to perform conditional quantile regression. In Figure~\ref{fig:quantile regression} we demonstrate the $95\%$ and $5\%$ quantile regression results (using sieve-SGD). The hyperparameter is tuned by wRV with $\xi = 1$. We also report the probability of the outcome lying between the regression curves.

So far, we have focused on the scenario where data points are generated from an unknown but static distribution. However, the analysis framework presented in this paper can be generalized to accommodate time-varying distributions. We are interested in adapting the procedures when $f_0$ experiences abrupt or incremental concept drifts. Some variants of the wRV method can provide valuable insights and practical tools for tracking the time-varying model selection target.

In applications with a large space of hyperparameters or candidate models, the dynamic average method proposed by \cite{fang2023online} may be more efficient in searching the candidate space.  It would be an interesting future direction to combine our framework with the dynamic averaging scheme and develop theoretically justifiable adaptive model searching methods. 

\section*{Acknowledgments}

We thank the Reviewers and Associate Editors for their constructive comments, which helped to improve this work.\\
Part of this work was completed while Tianyu Zhang was at Carnegie Mellon University. Tianyu Zhang's work at Carnegie Mellon University was partially supported by National Institute of Mental Health grant R01MH123184. Jing Lei is partially supported by National Science Foundation grants DMS-2310764 and DMS-2515687.

\bibliographystyle{chicago}
\bibliography{main}

%%%%%%%
\newpage
\appendix

\section{Proof of Theorem~\ref{th: main theorem}}
\label{app:proof of main theorem}
\begin{proof}
We first present the proof of comparing two estimators that are different in convergence rates, that is, under the assumption that $a_g > a_h$. Denote $H_n = \sum_{i=1}^n i^\xi \mathbb{E}[r(h_{i-1})]$. For any $\delta\in (0,1/4)$:
\begin{equation}
\label{eq: overall split}
    \begin{aligned}
& \mathbb{P}(\text { correct selection })\\
%%%%%%%
&=\mathbb{P}\left(\sum_{i=1}^ni^\xi\left\{h_{i-1}\left(X_i\right)-Y_i\right\}^2-i^\xi\epsilon_i^2 \geq \sum_{i=1}^n i^\xi\left\{g_{i-1}\left(X_i\right)-Y_i\right\}^2-i^\xi\epsilon_i^2\right)\\
%%%%%%%
& \geq \mathbb{P}\left(\sum_{i=1}^ni^\xi\left\{h_{i-1}\left(X_i\right)-Y_i\right\}^2-i^\xi\epsilon_i^2 \geq(1-\delta) H_n\right. \text { and }\\
& \left.\sum_{i=1}^ni^\xi\left\{g_{i-1}\left(X_i\right)-Y_i\right\}^2-i^\xi\epsilon_i^2 \leq(1-\delta) H_n\right)\\
& = 1-\mathbb{P}\left(\sum_{i=1}^ni^\xi\left\{h_{i-1}\left(X_i\right)-Y_i\right\}^2-i^\xi\epsilon_i^2 \leq(1-\delta) H_n\right. \text { or }\\
& \quad \left.\sum_{i=1}^ni^\xi\left\{g_{i-1}\left(X_i\right)-Y_i\right\}^2-i^\xi\epsilon_i^2 \geq(1-\delta) H_n\right)\\ 
& \geq 1-\mathbb{P}\left(\sum_{i=1}^ni^\xi\left\{h_{i-1}\left(X_i\right)-Y_i\right\}^2-i^\xi\epsilon_i^2 \leq(1-\delta) H_n\right)\\
& \quad -\mathbb{P}\left(\sum_{i=1}^ni^\xi\left\{g_{i-1}\left(X_i\right)-Y_i\right\}^2-i^\xi\epsilon_i^2 \geq(1-\delta) H_n\right).
    \end{aligned}
\end{equation}
For the first probability, we have:
\begin{equation}
\label{eq: easy first half}
    \begin{aligned}
& \mathbb{P}\left(\sum_{i=1}^n i^{\xi}\left\{h_{i-1}\left(X_i\right)-Y_i\right\}^2-i^{\xi} \epsilon_i^2 \leq(1-\delta) H_n\right) \\
& \leq \mathbb{P}\left(\left|\sum_{i=1}^n i^{\xi}\left\{h_{i-1}\left(X_i\right)-Y_i\right\}^2-i^{\xi} \epsilon_i^2  - H_n \right|\geq \delta H_n\right)\\
&\stackrel{(I)}{\leq} \mathbb{P}\left(\left|\sum_{i=1}^n i^{\xi}\left\{h_{i-1}\left(X_i\right)-Y_i\right\}^2-i^{\xi} \epsilon_i^2-H_n\right| \geq \delta n^{1-a_h+\xi} \tilde{M}_h / 2\right)\\
& \text { (for large } n \text { such that } \inf_{m\geq n}m^{a_h-\xi-1} H_m \geq \tilde M_h/2 \text { ) }
    \end{aligned}
\end{equation}
In step $(I)$ we used the fact that (it can be proved using the definition of limit, see an example in Lemma~\ref{lemma: elementary limit lemma}):
\begin{equation}
\label{eq: limit of sum}
    \begin{aligned}
        \liminf_{i\rightarrow\infty} i^{a_h} \mathbb{E}[r(h_{i-1})]& = M_h\\
        \Rightarrow \liminf_{n\rightarrow\infty} n^{a_h-1-\xi} \sum_{i=1}^n i^\xi \mathbb{E}[r(h_{i-1})] & \geq  M_h/(1+\xi - a_h) =: \tilde M_h\\
    \end{aligned}
\end{equation}

Then we can apply Corollary~\ref{main theorem each sequence} to show the last line of \eqref{eq: easy first half} converging to zero as $n\rightarrow\infty$.

To bound $\mathbb{P}\left(\sum_{i=1}^ni^\xi\left\{g_{i-1}\left(X_i\right)-Y_i\right\}^2-i^\xi\epsilon_i^2 \geq(1-\delta) H_n\right)$, we denote $G_n = \sum_{i=1}^n i^\xi \mathbb{E}[r(g_{i-1})]$. Then we have:
\begin{equation}
\label{eq: harder second}
\begin{aligned}
    & \mathbb{P}\left(\sum_{i=1}^ni^\xi\left\{g_{i-1}\left(X_i\right)-Y_i\right\}^2-i^\xi\epsilon_i^2 \geq(1-\delta) H_n\right)\\
    & \stackrel{(I)}{=} \mathbb{P}\left(\sum_{i=1}^ni^\xi\left\{g_{i-1}\left(X_i\right)-Y_i\right\}^2-i^\xi\epsilon_i^2-G_n \geq(1-\delta) H_n-G_n\right)\\
    &\text{(for large $n$ such that $G_n < (1-\delta)H_n$)}\\ 
    & \stackrel{(II)}{\leq} \mathbb{P}\left(\sum_{i=1}^ni^\xi\left\{g_{i-1}\left(X_i\right)-Y_i\right\}^2-i^\xi\epsilon_i^2-G_n \geq(1-2 \delta) H_n\right)\\
    &\text{(for large $n$ such that $G_n < \min\{\delta H_n, (1-\delta)H_n\} = \delta H_n$, since $\delta < 1/4$)}\\ 
    & \stackrel{(III)}{\leq}  \mathbb{P}\left(\sum_{i=1}^ni^\xi\left\{g_{i-1}\left(X_i\right)-Y_i\right\}^2-i^\xi\epsilon_i^2-G_n \geq (1-2 \delta)n^{1-a_h+\xi}\tilde M_h/2\right)\\
    &\text{(for large $n$ such that $G_n < \delta H_n$ and $|n^{a_h-\xi-1}H_n - \tilde M_h|\leq \tilde M_h/2$)}\\ 
\end{aligned}
\end{equation}
In step $(I),(II)$ we used $H_n$ is a positive sequence diverging in a faster rate than $G_n$ (by our assumption that $a_g > a_h$). In step $(III)$ we used \eqref{eq: limit of sum} again.
Under the conditions on the convergence rate and estimator stability, we can verify that $c_n = (1-2 \delta)n^{1-a_h+\xi}\tilde M_h/2$ satisfies the requirement in \eqref{clean cn scale}. Therefore the above deviation probability converges to $0$ according to Corollary~\ref{main theorem each sequence}.

In the case that $a_g = a_h$ the proof is very similar. The only twist is that we cannot pick any $\delta \in (0,1/4)$ as above. Instead, we should more carefully choose $\delta \in (0, 1-(M_g/M_h)^{1/3})$. We keep the same argument in \eqref{eq: overall split} and \eqref{eq: easy first half}, but modify \eqref{eq: harder second} as:
\begin{equation}
\label{eq: final lines of equal rate}
\begin{aligned}
    & \mathbb{P}\left(\sum_{i=1}^ni^\xi\left\{g_{i-1}\left(X_i\right)-Y_i\right\}^2-i^\xi\epsilon_i^2 \geq(1-\delta) H_n\right)\\
    & \leq \mathbb{P}\left(\sum_{i=1}^ni^\xi\left\{g_{i-1}\left(X_i\right)-Y_i\right\}^2-i^\xi\epsilon_i^2 \geq(1-\delta)^2 n^{1-a_h+\xi} \tilde{M}_h \right)\\
    % &\text{(for large $n$ such that $\left|n^{a_h-\xi-1} H_n-\tilde{M}_h\right| \leq \delta\tilde{M}_h$)}\\ 
    & \text{( for large $n$ such that $\inf _{m \geq n} m^{a_h-\xi-1} H_m \geq (1-\delta)\tilde{M}_h $)}\\
    & \stackrel{(IV)}{\leq} \mathbb{P}\left(\sum_{i=1}^ni^\xi\left\{g_{i-1}\left(X_i\right)-Y_i\right\}^2-i^\xi\epsilon_i^2-G_n \geq (1-\delta)^2 n^{1-a_h+\xi} \tilde{M}_h  - G_n\right)\\
    &\text{(for large $n$ such that $\inf _{m \geq n} m^{a_h-\xi-1} H_m \geq (1-\delta)\tilde{M}_h $ and $G_n \leq (1-\delta)^3 n^{1-a_h+\xi} \tilde{M}_h$)}\\ 
    & \leq  \mathbb{P}\left(\sum_{i=1}^ni^\xi\left\{g_{i-1}\left(X_i\right)-Y_i\right\}^2-i^\xi\epsilon_i^2-G_n \geq \delta(1-\delta)^2 n^{1-a_h+\xi}\tilde M_h\right)\\
\end{aligned}
\end{equation}

In step $(IV)$, it is possible to find large enough $n$ such that $G_n \leq (1-\delta)^3 n^{1-a_h+\xi} \tilde{M}_h$. In fact, similar to $H_n$, we also have (recall $a_g = a_h$):
\begin{equation*}
    \limsup_{n \rightarrow \infty} n^{a_h-1-\xi} \sum_{i=1}^n i^{\xi} \mathbb{E}[r(g_{i-1})] \leq M_g /\left(1+\xi-a_h\right)=: \tilde{M}_g.
\end{equation*}
This means we can find large enough $n$ such that
\begin{equation*}
    G_n \leq (\tilde M_g + \epsilon)n^{1+\xi-a_h} \text{ for any }\epsilon >0.
\end{equation*}
Specifically, we can take $\epsilon = (1-\delta)^3\tilde M_h - \tilde M_g$ (Due to our choice of $\delta \in\left(0,1-\left(M_g / M_h\right)^{1 / 3}\right)$, this is indeed a positive number). The last line of \eqref{eq: final lines of equal rate} can be bounded using Corollary~\ref{main theorem each sequence}.
\end{proof}

\begin{lemma}
\label{lemma: elementary limit lemma}
Suppose for a bounded positive sequence $\{h_i, i \in \mathbb{N}^+\}$ we have
\begin{equation*}
\lim _{i \rightarrow \infty} i^a h_i=M,
\end{equation*}
for some $a\in(0,1), M>0$. Then for any $\xi \geq 0$, we have
\begin{equation*}
\lim_{n\rightarrow\infty} n^{a-\xi-1}\sum_{i=1}^n i^\xi h_i = M/(1+\xi - a) =:\tilde M.
\end{equation*}
\end{lemma}
\begin{proof}
By definition of limit, we have: for any $\epsilon_1 > 0$, there exists some $N_1$ such that for any $i \geq N_1$:
\begin{equation*}
    |i^a h_i - M|\leq \epsilon_1.
\end{equation*}
We want to show that for any $\epsilon_2 > 0$, there exists some $N_2$ such that for any $n \geq N_2$:
\begin{equation*}
\left|n^{a-\xi-1} \sum_{i=1}^n i^{\xi} h_i-\tilde{M}\right| \leq \epsilon_2
\end{equation*}
We show the details of upper bounding $n^{a-\xi-1} \sum_{i=1}^n i^{\xi} h_i-\tilde{M}$:
\begin{equation*}
    \begin{aligned}
        & n^{a-\xi-1} \sum_{i=1}^n i^{\xi} h_i-\tilde{M}\\
        & \leq n^{a-\xi-1}\{\sum_{i=1}^{N_1} i^\xi h_i + \sum_{i = N_1 + 1}^n i^{\xi} (M+\epsilon_1)i^{-a}\} -\tilde{M}\\
        & = n^{a-\xi-1}\sum_{i=1}^{N_1} i^\xi h_i + n^{a-\xi-1}(M+\epsilon_1)\sum_{i = N_1 + 1}^n i^{\xi-a} -\tilde{M}
    \end{aligned}
\end{equation*}
The term $n^{a-\xi-1}\sum_{i=1}^{N_1} i^\xi h_i$ above can be arbitrarily small as $n \rightarrow \infty$ since $h_i$ is uniformly bounded. 

When $\xi - a > 0$ (the other cases can be done similarly):
\begin{equation*}
    \sum_{i = N_1 + 1}^n i^{\xi-a} \leq \int_{N_1+1}^{n+1} x^{\xi - a} dx \leq (1+\xi-a)^{-1}(n+1)^{1+\xi-a},
\end{equation*}
then
\begin{equation}
\label{eq: elementary difference}
    \begin{aligned}
n^{a-\xi-1}\left(M+\epsilon_1\right) \sum_{i=N_1+1}^n i^{\xi-a} \leq \frac{M+\epsilon_1}{1+\xi-a}\left(\frac{n+1}{n}\right)^{1+\xi-a}\\
\Rightarrow n^{a-\xi-1}\left(M+\epsilon_1\right) \sum_{i=N_1+1}^n i^{\xi-a} - \tilde M \leq \frac{M+\epsilon_1}{1+\xi-a}\left(\frac{n+1}{n}\right)^{1+\xi-a} - \tilde M
    \end{aligned}
\end{equation}
Since $\left(\frac{n+1}{n}\right)^{1+\xi-a}$ converges to $1$, \eqref{eq: elementary difference} can also be arbitrarily small.
\end{proof}

\section{Proof of Lemma~\ref{puzzle 1}}
\label{app: proof of puzzle 1}
\begin{proof}
All we need is Chebyshev's inequality:
\begin{align*}
     & \mathbb{P}\left(\left|\sum_{i=1}^n i^\xi\epsilon_i\left\{\hat{f}_{i-1}\left(X_i\right)-f_0\left(X_i\right)\right\}\right| \geq c_n\right) 
     \leq 2c_n^{-2}\text{var}\left(\sum_{i=1}^n i^\xi\epsilon_i\left\{\hat{f}_{i-1}\left(X_i\right)-f_0\left(X_i\right)\right\}\right)\\
    & = 2c_n^{-2}\sum_{i=1}^n \text{var}[i^\xi\epsilon_i \{\hat f_{i-1}(X_i) - f_0(X_i)\}]
    \\
    & \quad + 2c_n^{-2}\sum_{i\neq j} \text{cov}[i^\xi\epsilon_i\{\hat f_{i-1}(X_i) - f_0(X_i)\}, j^\xi\epsilon_j\{\hat f_{j-1}(X_j) - f_0(X_j)\}]\\
    & \stackrel{(I)}{=} 2c_n^{-2} \sum_{i=1}^n i^{2\xi}\mathbb{E}[\epsilon_i^2u_i]
    \leq 2 \sigma^2c_n^{-2}\sum_{i=1}^n i^{2\xi}\mathbb{E}[u_i] \stackrel{(II)}{\rightarrow} 0
\end{align*}
In step $(I)$, we used the noise variables $\epsilon_i$ that are centered and independent from each other. Recall $\sigma^2$ is the bound on the variance of $\epsilon_i$ (conditioned on $X_i$), for step $(II)$ we combined \Cref{ass: sup estimator_quality} and the assumption on sequence $c_n$. Note that $\limsup _{i \rightarrow \infty} i^{2 a} \mathbb{E}\left[\left\{\hat{f}_{i-1}(X)-f_0(X)\right\}^4\right]=C_1$ implies $\limsup _{i \rightarrow \infty} i^{ a} \mathbb{E}\left[\left\{\hat{f}_{i-1}(X)-f_0(X)\right\}^2\right]=C_1^{1/2}$.
\end{proof}

\section{Proof of Lemma~\ref{puzzle 2}}
\label{app: proof of puzzle 2}
\begin{proof}
Denote $Z_i = u_i - \mathbb{E}[u_i\mid F^{i-1}]$, we have $\mathbb{E}[Z_i] = 0 = \mathbb{E}[Z_i \mid F^{j}]$ for any $j\in [i-1]$. Then, applying Chebyshev's inequality:

\begin{equation*}
    \begin{aligned}
      \mathbb{P}\left(\left|\sum_{i=1}^n Z_i\right| 
      \geq c_n\right) &\leq 2 c_n^{-2} \operatorname{var}\left(\sum_{i=1}^n Z_i\right)\\
      & =2 c_n^{-2} \mathbb{E}\left(\sum_{i=1}^n Z_i\right)^2\\
      & = 2c_n^{-2}\left\{ \sum_{i=1}^n \mathbb{E} Z_i^2 + \sum_{i \neq j} \mathbb{E}[Z_i Z_j]\right\}.
    \end{aligned}
\end{equation*}
Note that
\begin{equation*}
    \mathbb{E}[Z_i Z_j] = \mathbb{E}[Z_j \mathbb{E}[Z_i | F^j]] = 0.
\end{equation*}
So we have
\begin{equation*}
    \begin{aligned}
     \mathbb{P}\left(\left|\sum_{i=1}^n Z_i\right| \geq c_n\right) 
     & \leq2 c_n^{-2} \sum_{i=1}^n \mathbb{E} Z_i^2=2 c_n^{-2} \sum_{i=1}^n \mathbb{E}\left(u_i-\mathbb{E}\left[u_i \mid F^{i-1}\right]\right)^2\\
      & = 2c_n^{-2} \sum_{i=1}^n \mathbb{E} u_i^2 + \mathbb{E}[E^2[u_i \mid F^{i-1}]] - 2\mathbb{E}[u_i \mathbb{E}[u_i \mid F^{i-1}]]\\
      & \leq 2c_n^{-2} \sum_{i=1}^n \mathbb{E} u_i^2  = 2 c_n^{-2} \sum_{i=1}^n i^{2 \xi} \mathbb{E}\left[\left\{\hat{f}_{i-1}\left(X_i\right)-f_0\left(X_i\right)\right\}^4\right]\stackrel{(I)}{\rightarrow} 0
    \end{aligned}
\end{equation*}
In step $(I)$ we used \Cref{ass: sup estimator_quality} and the assumption on sequence $c_n$.
\end{proof}

\newpage
\section{Proof of \Cref{puzzle 3} and some Related Examples}
\subsection{Proof of \Cref{puzzle 3}}
\label{app:proof of puzzle 3}
\begin{proof}
We will show that
\begin{equation*}
   c_n^{-1}\left |\sum_{i=1}^n \mathbb{E}[u_i | F^{i-1}] - \mathbb{E}[u_i]\right | \stackrel{P}{\rightarrow} 0 
\end{equation*}
which is implied by
\begin{equation*}
    \mathbb{E}c_n^{-2} \left|\sum_{i=1}^n \mathbb{E}[u_i | F^{i-1}] - \mathbb{E}[u_i]\right|^2 \rightarrow 0\,,
\end{equation*}
or equivalently
\begin{equation}
\label{eq:momentgoesto0}
    c_n^{-1} \left\|\sum_{i=1}^n \mathbb{E}[u_i | F^{i-1}] - \sum_{i=1}^n \mathbb{E}[u_i]\right\|_2 = c_n^{-1}\left(\mathbb{E}\left| \sum_{i=1}^n \mathbb{E}[u_i | F^{i-1}] - \sum_{i=1}^n \mathbb{E}[u_i]\right|^2\right)^{1/2} \rightarrow 0
\end{equation}

We rearrange the second moment term on the left-hand-side of \eqref{eq:momentgoesto0}:
    \begin{align}
    &\left\|\sum_{i=1}^n \mathbb{E}[u_i | F^{i-1}] - \sum_{i=1}^n \mathbb{E}[u_i]\right\|_2^2 \nonumber\\
    =&  \left\|\sum_{i=1}^n \left(\mathbb{E}[u_i | F^{i-1}]  - \mathbb{E}[u_i | F^{i-2}] + \mathbb{E}[u_i | F^{i-2}] - \mathbb{E}[u_i | F^{i-3}] + \cdot\cdot\cdot + \mathbb{E}[u_i | F^1] -  \mathbb{E}[u_i]\right)\right\|_2^2 \nonumber\\
    = &\left\|\sum_{i=1}^n \sum_{j=1}^{i-1} \mathbb{E}\left[u_i \mid F^j\right]-\mathbb{E}\left[u_i \mid F^{j-1}\right]\right\|_2^2 \nonumber\\
    = &\left\|\sum_{j=1}^{n-1} \sum_{i=j+1}^n \mathbb{E}\left[u_i \mid F^j\right]-\mathbb{E}\left[u_i \mid F^{j-1}\right]\right\|_2^2\nonumber\\
    \stackrel{(I)}{=}& \sum_{j=1}^{n-1}\left\|\sum_{i=j+1}^{n} \mathbb{E}[u_i | F^j] - \mathbb{E}[u_i | F^{j-1}]\right\|_2^2\nonumber\\
    \stackrel{(II)}{=}& \sum_{j=1}^{n-1}\left\|\sum_{i=j+1}^{n} \mathbb{E}[\nabla_j u_i|F^j]\right\|_2^2 \stackrel{(III)}{\leq} \sum_{j=1}^{n-1}\left(\sum_{i=j+1}^{n} \|\mathbb{E}[\nabla_j u_i|F^j]\|_2\right)^2 \,.\label{eq:longrearrange}
    \end{align}
In step $(I)$ we used that $M_j = \sum_{i = j+1}^n \mathbb{E}[u_i|F^j] - \mathbb{E}[u_i|F^{j-1}]$ is a martingale-difference sequence w.r.t. to the filtration $F^j$. Therefore $\mathbb{E}[M_lM_k] = 0$ for all $l \neq k$. In step $(II)$ we used the definition of the difference operator $\nabla_j$ and realize that $E(u_i(Z_i)|F^{j-1})=E(u_i(Z_i^j)|F^j)$. Step $(III)$ is triangular inequality.

Now we look at the scale of each term in the summation above. We will need the following decomposition of the stability of the loss function:
\begin{align}
&i^{-\xi}\nabla_j u_i \nonumber\\
=&  \nabla_j \{\hat f_{i-1}(X_i) - f_0(X_i)\}^2\nonumber\\
=&  \{\hat f_{i-1}(X_i;Z_{i-1}) - f_0(X_i)\}^2 - \{\hat f_{i-1}(X_i;Z_{i-1}^j) - f_0(X_i)\}^2 \nonumber\\
=&  \{\hat f_{i-1}(X_i;Z_{i-1}) - \hat f_{i-1}(X_i;Z_{i-1}^j)\}\{\hat f_{i-1}(X_i;Z_{i-1}) - f_0(X_i) + \hat f_{i-1}(X_i;Z_{i-1}^j) - f_0(X_i)\}\nonumber\\
=&\nabla_j \hat f_{i-1}(X_i)\{\hat f_{i-1}(X_i;Z_{i-1}) - f_0(X_i) + \hat f_{i-1}(X_i;Z_{i-1}^j) - f_0(X_i)\}\,.\label{important decomposition for stability}
\end{align}

We use the shorter notation $\hat f_{i-1}^{(j)}(x) = \hat f_{i-1}(x;Z_{i-1}^j)$. Then
\begin{align*}
        &\mathbb{E}[\nabla_j u_i | F^j]\\
     =& i^\xi \mathbb{E}[ \{\nabla_j \hat f_{i-1}(X_i)\}\{\hat f_{i-1}(X_i) - f_0(X_i) + \hat f_{i-1}^{(j)}(X_i) - f_0(X_i)\} | F^j]\\
    \leq & i^\xi\sqrt{\mathbb{E}[\{\nabla_j \hat f_{i-1}(X_i)\}^2| F^j]}  \sqrt{\mathbb{E}[\{\hat f_{i-1}(X_i) - f_0(X_i) + \hat f_{i-1}^{(j)}(X_i) - f_0(X_i)\}^2| F^j]}\,.
\end{align*}

Under \Cref{ass:stability} that $\sqrt{\mathbb{E}[\{\nabla_j \hat f_{i-1}(X_i)\}^2| F^j]}$ can be bounded by $\Xi_i = D_2i^{-b}$ almost surely, we can continue \eqref{eq:longrearrange} as following:
    \begin{align*}
        &\sum_{i = j+1}^n \left\|\mathbb{E}[\nabla_j u_i|F^j]\right\|_2\\
         \leq & \sum_{i = j+1}^n i^\xi\left\|\sqrt{\mathbb{E}[\{\nabla_j \hat f_{i-1}(X_i)\}^2| F^j]}  \sqrt{\mathbb{E}[\{\hat f_{i-1}(X_i) - f_0(X_i) + \hat f_{i-1}^{(j)}(X_i) - f_0(X_i)\}^2| F^j]}\right\|_2\\
        \leq & \sum_{i = j+1}^n i^\xi\Xi_i \left\| \sqrt{\mathbb{E}[\{\hat f_{i-1}(X_i) - f_0(X_i) + \hat f_{i-1}^{(j)}(X_i) - f_0(X_i)\}^2| F^j]}\right\|_2\\
        = & \sum_{i=j+1}^n i^{\xi} \Xi_i\left\|\hat{f}_{i-1}\left(X_i\right)-f_0\left(X_i\right)+\hat{f}_{i-1}^{(j)}\left(X_i\right)-f_0\left(X_i\right)\right\|_2\\
         \lesssim & \sum_{i=j+1}^n i^{\xi} \Xi_i\left\|\hat{f}_{i-1}\left(X_i\right)-f_0\left(X_i\right)\right\|_2.
    \end{align*}

\Cref{ass: sup estimator_quality} implies $\|\hat f_{i-1}(X_i) - f_0(X_i)\|_2 \lesssim Ci^{-a/2}$ for some $a\in(1/2,1]$, we can further simplify the above to be 
\begin{equation*}
\sum_{i = j+1}^n  \left\| \mathbb{E}[\nabla_j u_i|F^j]\right\|_2 \lesssim \sum_{i= j+1}^n i^{-c}, \text{ where }c:= (a/2 + b - \xi)\,.
\end{equation*}

\begin{enumerate}
    \item When $c > 1$,
\begin{equation*}
     \sum_{j=1}^{n-1}\left(\sum_{i=j+1}^{n}  \left\| \mathbb{E}[\nabla_j u_i|F^j]\right\|_2\right)^2 \lesssim \sum_{j=1}^{n-1}(j^{1-c} - n^{1-c})^2 \lesssim n^{3-2c} \vee \log n
\end{equation*}    

Plugging above into \eqref{eq:longrearrange},
\begin{equation*}
    c_n^{-1} \left\|\sum_{i=1}^n \mathbb{E}[u_i | F^{i-1}] - \sum_{i=1}^n \mathbb{E}[u_i]\right\|_2 
    \lesssim c_n^{-1} n^{3/2-c} \vee c_n^{-1}\log^{1/2} n
\end{equation*}
Under the assumptions of $c_n$, the RHS above converges to $0$ as $n\rightarrow\infty$.

\item When $c = 1$, by definition of $c$, we have $a/2 +b -\xi= 1$. Our assumptions on $c_n$ are reduced to $\lim_n c_n^{-1} n^{1/2}\log^{1/2} n = 0$. So in this case we have
\begin{equation*}
\begin{aligned}
    & \sum_{j=1}^{n-1}\left(\sum_{i=j+1}^{n}  \left\| \mathbb{E}[\nabla_j u_i|F^j]\right\|_2\right)^2 \\
    & \lesssim \sum_{j=1}^{n-1}(\log(n) - \log(j))^2\\
    & = (n-1)\log^2(n)-2\log n\sum_{j=1}^n \log j + \sum_{j=1}^n \log^2j\\
    & \stackrel{(I)}{\lesssim}n\log^2n - 2\log n (n\log n - n/2) + n\log^2n\\
    & = n\log n
\end{aligned}
\end{equation*}
In step (I) we used Stirling's approximation treating $\sum_j \log j = \log n!$.
% \begin{equation*}
%      \sum_{j=1}^{n-1}\left(\sum_{i=j+1}^{n}  \left\| \mathbb{E}[\nabla_j u_i|F^j]\right\|_2\right)^2 \lesssim \sum_{j=1}^{n-1}\log^2 n \lesssim n\log^2 n\,,
% \end{equation*} 
% and
% \begin{equation*}
%     c_n^{-1} \left\|\sum_{i=1}^n \mathbb{E}[u_i | F^{i-1}] - \sum_{i=1}^n \mathbb{E}[u_i]\right\|_2 
%     \lesssim c_n^{-1} n^{1/2}\log n \rightarrow 0 \text{ as }n\rightarrow\infty
% \end{equation*}

\item When $c<1$, we have
\begin{equation*}
     \sum_{j=1}^{n-1}\left(\sum_{i=j+1}^{n}  \left\| \mathbb{E}[\nabla_j u_i|F^j]\right\|_2\right)^2\lesssim \sum_{j=1}^{n-1}(n^{1-c} - j^{1-c})^2 \lesssim n^{3-2c}
\end{equation*}  
Similar to case 1, we have $ c_n^{-1}\left\|\sum_{i=1}^n \mathbb{E}[u_i | F^{i-1}] - \sum_{i=1}^n \mathbb{E}[u_i]\right\|_2 \rightarrow 0$ as $n \rightarrow \infty$. \qedhere
\end{enumerate}
\end{proof}

\subsection{\Cref{puzzle 3} and the Variance of Mean Estimation}\label{app: non_martingale_lemma_is_sharp}
In this section, we are going to establish some sharp lower bounds (up to a constant) on the variance for the random variable considered in \Cref{puzzle 3}. Combined with numerical results, we conjecture the deviation presented in \Cref{puzzle 3} can not be further tightened. 

We are taking sequential mean estimation as an example. Suppose we observe a sequence of $Y_i\in\mathbb{R}$, $i\geq 1$, and our goal is to estimate its mean (so the covariates degenerate to constant $1$). The underlying distribution of $Y_i$ is standard normal $\mathcal{N}(0,1^2)$. 

Goal is estimating $\mu = \mathbb{E}[Y]$, our estimator at step $i$ is $\hat \mu_i = i^{-1}\sum_{l=1}^i Y_i$. In this case, $u_i = i^\xi (\hat \mu_i - \mu)^2$. The random quantity of interest is 
$$\mathcal{D}_n = \sum_{i=1}^n \mathbb{E}\left[u_i \mid F^{i-1}\right]-\mathbb{E}\left[u_i\right].$$

We can explicitly calculate the perturb-one version of $u_{i+1}$:
\begin{equation*}
\begin{aligned}
\nabla_j u_{i+1} & =i^{\xi}\left(\hat{\mu}_i-\hat{\mu}_i^j\right)\left(\hat{\mu}_i+\hat{\mu}_i^j-2 \mu\right) \\
& =i^{\xi-2}\left(Y_j-Y_j^{\prime}\right)\left(\sum_{l=1}^i Y_l+\sum_{l=1, l \neq j}^i Y_l+Y_j^{\prime}-2 i \mu\right).
\end{aligned}
\end{equation*}
Take the conditional expectation
\begin{equation}\label{eq: take_conditional_expectation_mean}
\mathbb{E}\left[\nabla_j u_{i+1} \mid F^j\right]=i^{\xi-2} S_{j},
\end{equation}
where
\begin{equation*}
\begin{aligned}
S_j& =2\left(Y_j-\mu\right)\left(\sum_{l=1}^{j-1} Y_l-j \mu\right)+\left(Y_j^2-\alpha^2\right)\\
\alpha^2 & = \mathbb{E}[Y^2].
\end{aligned}
\end{equation*}
The detail is presented in \Cref{lemma: take_conditional_expectation}. Therefore,
\begin{equation}\label{eq: sum_of_nabla}
\sum_{i=j+1}^n \mathbb{E}\left[\nabla_j u_i \mid F^j\right]=S_j\left(\sum_{i=j+1}^n(i-1)^{\xi-2}\right) =: S_j A(j;n).
\end{equation}
Treat LHS of \Cref{eq: sum_of_nabla} as a random variable, and we calculate its second moment
\begin{equation*}
\mathbb{E}\left[\left(\sum_{i=j+1}^n \mathbb{E}\left[\nabla_j u_i \mid F^j\right]\right)^2\right] = A^2(j;n) \mathbb{E}[(S_j)^2].
\end{equation*}
We can directly calculate
\begin{equation*}
\begin{aligned}
    \mathbb{E}[(S_j)^2]  &  = 4 {\rm var}^2(Y) j + (\mathbb{E}[Y^4]- 4\mu \mathbb{E}[Y^3] - 5\alpha^4 - 8\mu^4 + 16\mu^2\alpha^2)\\
  & =: Bj + D.
\end{aligned}
\end{equation*}
According to the proof of \Cref{puzzle 3}, we know 
\begin{equation*}
\begin{aligned}
  &\operatorname{var}\left(\mathcal{D}_n\right)\\
=&\sum_{j=1}^{n-1} \mathbb{E}\left[\left(\sum_{i=j+1}^n \mathbb{E}\left[\nabla_j u_i \mid F^j\right]\right)^2\right]\\
= & \sum_{j=1}^{n-1}A^2(j;n)(Bj+D)
\end{aligned}
\end{equation*}
When $\xi = 0$, $A(j;n) \geq (j+1/2)/j^2 - (n+1)/n^2$ (properties of Polygamma functions). So we have
\begin{equation*}
\begin{aligned}
\sum_{j=1}^{n-1} A^2(j ; n)Bj 
&\geq B\sum_{j=1}^{n-1}((j+1/2)/j^2 - (n+1)/n^2)^2j\\
& = B\sum_{j=1}^{n-1} 
\frac{1}{4 j^3}+\frac{1}{j^2}+\frac{1}{j}+\frac{j}{n^4}+\frac{2 j}{n^3}-\frac{2}{n^2}-\frac{1}{j n^2}+\frac{j}{n^2}-\frac{2}{n}-\frac{1}{j n}\\
& \geq B^\prime \log n \quad \text{ for large }n. 
\end{aligned}
\end{equation*}
Similarly, we can show that the $\sum_{j=1}^{n-1} A(j ; n) D$ term is greater than a negative constant. This concludes that $\operatorname{var}\left(\mathcal{D}_n\right)$ is greater than a constant times $\log n$ for large enough $n$. This corresponds to the $a = b = 1$, $\xi = 0$ case. Compared with the bound in \Cref{puzzle 3}, we know the order of $c_n$ cannot be improved in general: We required $c_n^2$ to be slightly larger than $\log n$, and the latter matches the lower bound of ${\rm var}(\mathcal{D}_n)$.

When $\xi = 2$, $A(j;n) = n-j$. We have 
\begin{equation*}
    \begin{aligned}
    \sum_{j=1}^{n-1} A^2(j ; n) B j
    & = B\sum_{j=1}^{n-1}(n-j)^2j\\
    & = n^2 \sum_{j=1}^{n-1}j - 2n \sum_{j=1}^{n-1}j^2 + \sum_{j=1}^{n-1}j^3\\
    & \geq n^4/12 +o(n^4).
    \end{aligned}
\end{equation*}

This corresponds to the $a = b = 1$, $\xi = 2$ case, and we can compare the bound in \Cref{puzzle 3}.

\begin{figure}[t]
    \centering
\includegraphics[width = \textwidth]{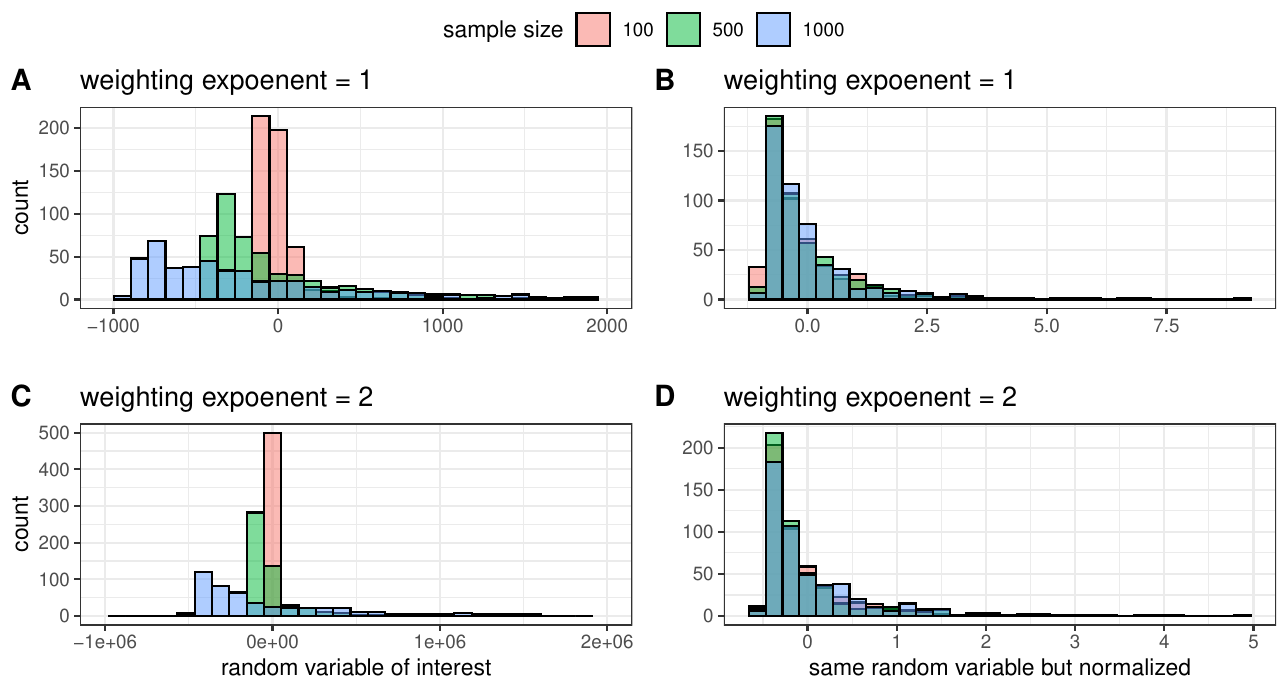}
    \caption{Distribution of the random variable studied in \Cref{puzzle 3}. \textbf{(A)} and \textbf{(C)}, histograms of $\mathcal{D}_n$ at different sample size $n$; \textbf{(A)} and \textbf{(B)}, $\xi = 1$; \textbf{(C)} and \textbf{(D)}, $\xi = 2$; \textbf{(B)}, divide the random variable in \textbf{(A)} by its predicted SD $n$; \textbf{(D)}, divide the random variable in \textbf{(C)} by its predicted SD $n^2$. The histograms are created using $500$ repeats.}
    \label{fig: non_martingale_lemma_sharp}
\end{figure}

In \Cref{fig: non_martingale_lemma_sharp}, we present the histograms of $\mathcal{D}_n$ for the mean estimation example with $\xi = 1$ and $2$. The variances can be shown to be larger than $n^2$ and $n^4$, respectively, using the argument above. After dividing them by their standard deviation (SD)---more accurately, a sequence of the same order as their SD lower bounds---we observe their normalized distributions do not visually vary between different sample sizes (and are close to a tamed continuous distribution). Given the above evidence, we conjecture $c_n$ in \Cref{puzzle 3} needs to be of order larger than $\sqrt{{\rm var}(\mathcal{D}_n)}$ to have the probability in \eqref{eq: non_martingale_concentration} converging to zero.

\begin{lemma}\label{lemma: take_conditional_expectation}
Using the notation in \Cref{eq: take_conditional_expectation_mean},
\begin{equation*}
    \mathbb{E}\left[\nabla_j u_{i+1} \mid F^j\right]=i^{\xi-2} S_j
\end{equation*}
where
\begin{equation*}
    S_j=2\left(Y_j-\mu\right)\left(\sum_{l=1}^{j-1} Y_l-j \mu\right)+\left(Y_j^2-\alpha^2\right).
\end{equation*}
\end{lemma}
\begin{proof}
\begin{equation*}
    \begin{aligned}
       & \mathbb{E}\left[\nabla_j u_{i+1} \mid F^j\right]\\
       & =i^{\xi - 2} \mathbb{E}\left[\left(Y_j-Y_j^{\prime}\right)\left(2\sum_{l=1, l \neq j}^i Y_l+Y_j + Y_j^{\prime}-2 i \mu\right) \mid F^j\right]\\
       & = i^{\xi - 2} Y_j \mathbb{E}\left[2 \sum_{l=1, l \neq j}^i Y_l+Y_j+Y_j^{\prime}-2 i \mu \mid F^j\right]\\
       &\quad -i^{\xi - 2} \mathbb{E}\left[Y_j^{\prime}\left(2 \sum_{l=1, l \neq j}^i Y_l+Y_j+Y_j^{\prime}-2 i \mu\right) \mid F^j\right]\\
       & = i^{\xi-2} Y_j\left(2 \sum_{l=1}^{j-1} Y_l+2(i-j) \mu+Y_j+\mu-2 i \mu\right)\\
       &\quad - i^{\xi-2}\left(2 \mu \sum_{l=1}^{j-1} Y_l+2 \mu(i-j) \mu+Y_j \mu+\alpha^2-2 i \mu^2\right)\\
       & = i^{\xi - 2}\left(2\left(Y_j-\mu\right) \sum_{l=1}^{j-1} Y_l+2\left(Y_j-\mu\right)(i-j) \mu+Y_j^2-\alpha^2-2 i \mu\left(Y_j-\mu\right)\right)\\
       & = i^{\xi - 2}2\left(Y_j-\mu\right)\left( \sum_{l=1}^{j-1} Y_l- j \mu\right)+i^{\xi - 2}\left(Y_j^2-\alpha^2\right). \end{aligned}
\end{equation*}
\end{proof}

\newpage
\section{Consistency under Diverging Candidate Cardinality}\label{app: infinite model}

We present the proof of Theorem~\ref{th: diverging cardinality} below.
\begin{proof}
We define $H_n = \sum_{i=1}i^\xi \min_{k\in K_i} \mathbb{E}[r(h_{i-1}^{(k)})]$ as the accumulated error of the best alternative models.

\begin{equation}\label{eq: final bound, decomposition}
\begin{aligned}
    & \mathbb{P}(\text { correct selection })\\
= & \mathbb{P}\left(\inf_{k \in K_n}\sum_{i=1}^n i^{\xi}\left\{h^{(k)}_{i-1}\left(X_i\right)-Y_i\right\}^2-i^{\xi} \epsilon_i^2 \geq \sum_{i=1}^n i^{\xi}\left\{g_{i-1}\left(X_i\right)-Y_i\right\}^2-i^{\xi} \epsilon_i^2\right)\\
\geq & \mathbb{P}\left(\inf_k \sum_{i=1}^n i^{\xi}\left\{h^{(k)}_{i-1}\left(X_i\right)-Y_i\right\}^2-i^{\xi} \epsilon_i^2 \geq(1-\delta) H_n \text{ and }\right. \\
& \quad \left.\sum_{i=1}^n i^{\xi}\left\{g_{i-1}\left(X_i\right)-Y_i\right\}^2-i^{\xi} \epsilon_i^2 \leq(1-\delta) H_n\right)\\
= & 1-\mathbb{P}\left(\inf_{k}\sum_{i=1}^n i^{\xi}\left\{h_{i-1}^{(k)}\left(X_i\right)-Y_i\right\}^2-i^{\xi} \epsilon_i^2 \leq(1-\delta) H_n\right. \text { or } \\
& \left.\sum_{i=1}^n i^{\xi}\left\{g_{i-1}\left(X_i\right)-Y_i\right\}^2-i^{\xi} \epsilon_i^2 \geq(1-\delta) H_n\right)\\
\geq &  1-\mathbb{P}\left(\inf_k \sum_{i=1}^n i^{\xi}\left\{h_{i-1}^{(k)}\left(X_i\right)-Y_i\right\}^2-i^{\xi} \epsilon_i^2 \leq(1-\delta) H_n\right) \\
& \quad-\mathbb{P}\left(\sum_{i=1}^n i^{\xi}\left\{g_{i-1}\left(X_i\right)-Y_i\right\}^2-i^{\xi} \epsilon_i^2 \geq(1-\delta) H_n\right)\\
\geq &  1-\sum_{k\in [K_n]}\mathbb{P}\left(\sum_{i=1}^n i^{\xi}\left\{h_{i-1}^{(k)}\left(X_i\right)-Y_i\right\}^2-i^{\xi} \epsilon_i^2 \leq(1-\delta) H_n\right)\\
& \quad-\mathbb{P}\left(\sum_{i=1}^n i^{\xi}\left\{g_{i-1}\left(X_i\right)-Y_i\right\}^2-i^{\xi} \epsilon_i^2 \geq(1-\delta) H_n\right).
\end{aligned}
\end{equation}
Define $H_n^{(k)}=\sum_{i=1}^n i^{\xi} \mathbb{E}[r(h_{i-1}^{(k)})]$, for each term in the $\sum_{k\in [K_n]}$ summation:
\begin{equation}\label{eq: decomposition into three terms infinite}
    \begin{aligned}
    & \mathbb{P}\left(\sum_{i=1}^n i^{\xi}\left\{h_{i-1}^{(k)}\left(X_i\right)-Y_i\right\}^2-i^{\xi} \epsilon_i^2 \leq(1-\delta) H_n\right)\\
    & \leq \mathbb{P}\left(\sum_{i=1}^n i^{\xi}\left\{h_{i-1}^{(k)}\left(X_i\right)-Y_i\right\}^2-i^{\xi} \epsilon_i^2 \leq(1-\delta) H_n^{(k)}\right)\\
    & \leq \mathbb{P}\left(\left|\sum_{i=1}^n i^{\xi}\left\{h_{i-1}^{(k)}\left(X_i\right)-Y_i\right\}^2-i^{\xi} \epsilon_i^2-H_n^{(k)}\right| \geq \delta H_n^{(k)}\right)\\
    & \stackrel{(I)}{\leq} \mathbb{P}\left(\left|\sum_{i=1}^n i^{\xi}\left\{h_{i-1}^{(k)}\left(X_i\right)-Y_i\right\}^2-i^{\xi} \epsilon_i^2-H_n^{(k)}\right| \geq \delta \sum_{i=1}^n i^{-a_k + \xi} M_k\right)\\
    & \leq \mathbb{P}\left(\left|\sum_{i=1}^n i^{\xi}\left\{h^{(k)}_{i-1}\left(X_i\right)-Y_i\right\}^2-i^{\xi} \epsilon_i^2-H_n^{(k)}\right| \geq \delta n^{1-a_k+\xi} \tilde{M}_k/ 2\right)\\
    & \quad \text{ for large }n\text{ not depending on }k\text{,and } \tilde M_k  = M_k/(1 + \xi - a_k)\\
    & \stackrel{(II)}{\leq} \mathbb{P}\left(\left|\sum_{i=1}^n i^{\xi} \epsilon_i\left\{h_{i-1}^{(k)}\left(X_i\right)-f_0\left(X_i\right)\right\}\right| \geq c_n\right) +\\
    & \quad \mathbb{P}\left(\left|\sum_{i=1}^n u_i^{(k)}-\mathbb{E}\left[u_i^{(k)} \mid F^{i-1}\right]\right| \geq c_n\right)+ \\
    & \quad \mathbb{P}\left(\left|\sum_{i=1}^n \mathbb{E}\left[u_i^{(k)} \mid F^{i-1}\right]- H_n^{(k)}\right| \geq c_n\right).
    \end{aligned}
\end{equation}
Here $c_n = \delta n^{1-a_k+\xi} \tilde{M}_k / 6$ and recall $u_i^{(k)} = i^{\xi}\left\{h^{(k)}_{i-1}\left(X_i\right)-f_0\left(X_i\right)\right\}^2$. In step (I) we used Assumption~\ref{ass: subW estimator quality}. In step (II) we applied the decomposition \eqref{concentration splitting}.

For each of the probabilities in the last line of \eqref{eq: decomposition into three terms infinite}, we can derive an exponential tail concentration inequality. The details are presented in Lemma~\ref{lemma: infinite concentration 1}-\ref{lemma: infinite concentration 3}. Then we have 
\begin{equation}\label{eq: bound depending on k}
\mathbb{P}\left(\sum_{i=1}^n i^{\xi}\left\{h_{i-1}^{(k)}\left(X_i\right)-Y_i\right\}^2-i^{\xi} \epsilon_i^2 \leq(1-\delta) H_n\right) \lesssim  \exp(-K^\prime (\mathcal{C}^{(k)}_n)^{\theta^\prime}),
\end{equation}
where $K^\prime, \theta^\prime$ are constants depending on $K_\epsilon, \theta_\epsilon, K_{\rm est}, \theta_{\rm est}, K_{\rm stab}, \theta_{\rm stab}$. The sequence 
\begin{equation*}
    \mathcal{C}^{(k)}_n = n^{(1-a_k)/2}\wedge \log ^{-1} n\left(n^{b_k-a_k / 2-1 / 2} \wedge n^{1-a_k+\xi}\right).
\end{equation*}
Now we applied the assumed bound (U1). Denote $\mathcal{C}_n := \min_{k\in [K_n]} (b_k-a_k / 2-1 / 2) \wedge (1-a_k)$, we have a uniform lower bound for $\mathcal{C}^{(k)}_n$ not depending on $k$:
\begin{equation*}
    \log^{-1} n \cdot n^{\mathcal{C}_n/2} \lesssim \mathcal{C}^{(k)}_n.
\end{equation*}
Go back to \eqref{eq: bound depending on k}:
\begin{equation*}
\mathbb{P}\left(\sum_{i=1}^n i^{\xi}\left\{h_{i-1}^{(k)}\left(X_i\right)-Y_i\right\}^2-i^{\xi} \epsilon_i^2 \leq(1-\delta) H_n\right)
\lesssim \exp(-K^\prime (\log ^{-1} n \cdot n^{\mathcal{C}_n / 2})^{\theta^\prime}).
\end{equation*} 
Then go back to \eqref{eq: final bound, decomposition}:
\begin{equation*}
    \begin{aligned}
        &\mathbb{P}(\text { correct selection }) \\
     \geq  & 1 - K_n \exp \left(-K^{\prime}\left(\log ^{-1} n \cdot n^{\mathcal{C}_n / 2}\right)^{\theta^{\prime}}\right) \\
     &\quad -\mathbb{P}\left(\sum_{i=1}^n i^{\xi}\left\{g_{i-1}\left(X_i\right)-Y_i\right\}^2-i^{\xi} \epsilon_i^2 \geq(1-\delta) H_n\right)\\
      = & 1 - K_n \exp \left(-K^{\prime}\left(\log ^{-1} n \cdot n^{\mathcal{C}_n / 2}\right)^{\theta^{\prime}}\right) + o(1).
    \end{aligned}
\end{equation*}
Establishing the last equality above is identical to the proof of Theorem~\ref{th: main theorem}, leveraging assumption (U2) that all the alternative models are worse than $\{g_i\}$. So long as
\begin{equation}\label{eq: general bound Kn}
K^{\prime}\left(\log ^{-1} n \cdot n^{\mathcal{C}_n / 2}\right)^{\theta^{\prime}} \gg \log K_n
\end{equation}
as $n\rightarrow\infty$, we can obtain consistency with a diverging number of models under comparison. We can verify that the conditions
\begin{equation*}
    \begin{aligned}
        (\log \log n) / \log n & = o(\mathcal{C}_n)\\
       \log K_n & \lesssim \log n,
    \end{aligned}
\end{equation*}
are sufficient for \eqref{eq: general bound Kn} to hold. 
\end{proof}

\subsection{Technical Results for Exponential Concentration}

In this section, we will present the details of concentration results under the sub-Weibull assumptions. We will repeatedly apply the following result to derive exponential tails for martingale difference sequences. This result can be seen as a Hoeffding-type inequality for sub-Weibull random variables. We replicate its proof as in \cite[Lemma~A.3]{kissel2022high}, which applies a useful moment martingale inequality of \cite{rio2009moment}.

\begin{lemma}\label{lemma: hoeffding sub-W}
    Let $D=\sum_{i=1}^n D_i$ where the sequence $\left(D_i, i\in [n]\right)$ satisfies 
    \begin{enumerate}
        \item martingale property: $\mathbb{E}\left(D_i \mid D_j: j\in [i-1]\right)=0$ for all $2 \leq i \leq n$, and $\mathbb{E} D_1=0$.
        \item sub-Weibull tail: $\left\|D_i\right\|_q \leq c K_i q^{1 / \alpha_i}$ for some $c, \alpha_i>0$ and all $q \geq 1$.
    \end{enumerate}
Then we have, for $\alpha^{\prime}=2 \alpha^*/(2+\alpha^*)$ with $\alpha^*=\min _{j \leq n} \alpha_j$ and a positive constant $c^{\prime}$,
\begin{equation*}
\|D\|_q \leq c^{\prime}\left(\sum_{i=1}^n K_i^2\right)^{1 / 2} q^{1 / \alpha^{\prime}}, \quad \forall q \geq 1.
\end{equation*}
If $K_i=K$ for all $i \in[n]$, then
\begin{equation*}
\|D\|_q \leq c^{\prime} \sqrt{n} K q^{1 / \alpha^{\prime}}, \quad \forall q \geq 1.
\end{equation*}
\end{lemma}
\begin{proof}
    By Theorem~2.1 of \cite{rio2009moment}, we have for any $q \geq 2$
    \begin{equation*}
        \|D\|_q \leq\left[(q-1) \sum_{i=1}^n\left\|D_i\right\|_q^2\right]^{1 / 2} \leq\left[C(q-1) q^{2 / \alpha^*} \sum_{i=1}^n K_i^2\right]^{1 / 2} \leq C^{1 / 2} q^{\frac{2+\alpha^*}{2 \alpha^*}}\left(\sum_{i=1}^n K_i^2\right)^{1 / 2}
    \end{equation*}
where $C$ is a constant depending only on $c$, and the second inequality follows from the assumption $\left\|D_i\right\|_q \leq c K_i q^{1 / \alpha_i}$. For any $q < 2$, we apply $\|D\|_q \leq \|D\|_2$.
\end{proof}

Now we present the three bounds applied in \eqref{eq: decomposition into three terms infinite}.

\begin{lemma}\label{lemma: infinite concentration 1}
Under Assumptions~\ref{ass: subW noise}, and \ref{ass: subW estimator quality}, for any $c_n >0$, we have for each $\hat{f}_{i-1}^{(k)}$:
\begin{equation*}
\mathbb{P}\left(\left|\sum_{i=1}^n i^{\xi} \epsilon_i\left\{\hat{f}_{i-1}^{(k)}\left(X_i\right)-f_0\left(X_i\right)\right\}\right| \geq c_n\right) \leq \exp \left(-K'\left(c_n / \mathbb{c}_n^{(k)}\right)^{\theta^{\prime}}\right),
\end{equation*}
where $\mathbb{c}_n^{(k)}=n^{\xi-a_k / 2+1 / 2}$, $K' = \left(K_\epsilon K_{\rm est}\right)^{-\theta^{\prime}}$, and $\theta^{\prime} = 2 \theta_\epsilon \theta_{\rm est} /\left(2 \theta_\epsilon+2 \theta_{\rm est}+\theta_\epsilon \theta_{\rm est}\right)$. Specifically, when $c_n = n^{1-a_k+\xi}$:
\begin{equation*}
\mathbb{P}\left(\left|\sum_{i=1}^n i^{\xi} \epsilon_i\left\{\hat{f}_{i-1}^{(k)}\left(X_i\right)-f_0\left(X_i\right)\right\}\right| \geq c_n\right) \leq 
\exp \left(-K^{\prime}n^{\theta^\prime(1-a_k)/2} \right).
\end{equation*}
\end{lemma}
\begin{proof}

We drop the super-(sub-)script $k$ for notation simplicity. Denote $D_i = i^{\xi} \epsilon_i\left\{\hat{f}_{i-1}\left(X_i\right)-f_0\left(X_i\right)\right\}$, it is direct to check that, 
$$\mathbb{E}\left[D_i \mid D_j : j\in[i-1]\right]=0$$ 
for all $2 \leq i \leq n$, and $\mathbb{E} D_1 =0$ ($\epsilon_i$ is centered conditioning on $X_i$). 

We can show that under the tail assumptions, the random variable $D_i $ also has an exponential tail. To this end, it is sufficient to derive bounds for all the moments of $D_i$:
\begin{equation*}
    \begin{aligned}
        & \mathbb{E}[(D_i )^q ] \\
        & = i^{\xi q} \mathbb{E}\left[\left(\epsilon_i\right)^q\left\{\hat{f}_{i-1} \left(X_i\right)-f_0\left(X_i\right)\right\}^q\right]\\
        &\leq i^{\xi q} i^{-a  q / 2} \sqrt{\mathbb{E}\left[\epsilon_i^{2 q}\right]} \sqrt{\mathbb{E}\left[i^{a q}\left\{\hat{f}_{i-1} \left(X_i\right)-f_0\left(X_i\right)\right\}^{2q} \right]}\\
        & \stackrel{(I)}{\leq} i^{(\xi-a  / 2) q}\left(c K_\epsilon q^{1 / \theta_\epsilon}\right)^q\left(c K_{\rm est} q^{1 / \theta_{\rm est}}\right)^q\\
        & =\left(i^{\xi-a  / 2} c^2 K_\epsilon K_{\rm est} q^{\theta_\epsilon^{-1}+\theta_{\rm est}^{-1}}\right)^q\\
        \Rightarrow (\mathbb{E}[(D_i)^q])^{1/q} & \leq i^{\xi-a / 2} c^2 K_\epsilon K_{\rm est} q^{\theta_\epsilon^{-1}+\theta_{\rm est}^{-1}}
    \end{aligned}
\end{equation*}
In step (I), we used the moment bounds for sub-Weibull random variables. So $D_i $ is $(ci^{\xi-a  / 2}  K_\epsilon K_{\rm est}, (\theta_\epsilon^{-1} + \theta_{\rm est}^{-1})^{-1})$-sub-Weibull. Given the martingale structure and sub-Weibull conditions for each $D_i$, we can employ Lemma~\ref{lemma: hoeffding sub-W} to derive sub-Weibull properties for $\sum_{i=1}^n D_i $:
\begin{equation*}
\begin{aligned}
\left\|\sum_{i=1}^n D _i\right\|_q & \leq c\left(\sum_{i=1}^n\left(K_\epsilon K_{\rm est}\right)^2 i^{2 \xi-a }\right)^{1 / 2} q^{1 / \theta^{\prime}} \text { for all } q \geq 1 \\
& \leq c K_\epsilon K_{\rm est} \mathbb{c}_n q^{1 / \theta^{\prime}} \text{ where }\mathbb{c} _n = n^{\xi-a  / 2+1 / 2} 
\end{aligned}
\end{equation*}
and $\theta^\prime = 2\theta_\epsilon \theta_{\rm est}/(2\theta_\epsilon + 2\theta_{\rm est} + \theta_\epsilon \theta_{\rm est})$. Then we know the deviation of interest has an exponential tail:
\begin{equation*}
\begin{aligned}
& \mathbb{P}\left(\left|\sum_{i=1}^n i^{\xi} \epsilon_i\left\{\hat{f}_{i-1}\left(X_i\right)-f_0\left(X_i\right)\right\}\right| \geq c_n\right) \\
& =\mathbb{P}\left(\frac{\left|\sum_{i=1}^n D_i\right|}{K_\epsilon K_{\rm est} \mathbb{c}_n} \geq \frac{c_n}{K_\epsilon K_{\rm est} \mathbb{c} _n}\right) \\
& \leq c \exp \left(-(K_{\epsilon}K_{\rm est})^{-\theta^\prime}\left(c_n /\mathbb{c}_n\right)^{\theta^{\prime}}\right).
\end{aligned}
\end{equation*}
\end{proof}

We can similarly derive a bound for $\sum_{i=1}^n u_i^{(k)}-\mathbb{E}\left[u_i^{(k)} \mid F^{i-1}\right]$.

\begin{lemma}\label{lemma: infinite concentration 2}
    Under Assumptions~\ref{ass: subW noise}, and \ref{ass: subW estimator quality}, for any $c_n >0$, we have for each $\hat{f}_{i-1}^{(k)}$:
    \begin{equation*}
\mathbb{P}\left(\left|\sum_{i=1}^n u_i^{(k)}-\mathbb{E}\left[u_i^{(k)} \mid F^{i-1}\right]\right| \geq c_n\right) \leq \exp \left(-K^{\prime}\left(c_n n^{-\xi+a_k-1 / 2}\right)^{\theta^\prime}\right)
\end{equation*}
    where $K^{\prime}=\left(2^{\left(\theta_{\rm {est }}+2\right) / \theta_{\rm {est }}} K_{\rm {est }}\right)^{-\theta_{\rm {est }} / 2}$ and $\theta^\prime = \theta_{\rm est}/2$. Specifically, when $c_n=n^{1-a_k+\xi}$:
    \begin{equation*}
\mathbb{P}\left(\left|\sum_{i=1}^n u_i^{(k)}-\mathbb{E}\left[u_i^{(k)} \mid F^{i-1}\right]\right| \geq c_n\right) \leq \exp(-K^{\prime}n^{\theta^\prime/2} ).
    \end{equation*}
\end{lemma}

\begin{proof}
    We drop the $k$ in the superscript for notation simplicity. For each $k$, denote $D_i = u_i -\mathbb{E}\left[u_i  \mid F^{i-1}\right]$, we need to derive some moments bounds on $D_i $.
    \begin{equation*}
    \begin{aligned}
       \|D_i \|_q 
        &\leq \|u_i \|_q + \|\mathbb{E}[u_i  \mid F^{i-1}]\|_q\\
        & \leq 2\|u_i \|_q\\
        & =2\left(\mathbb{E}\left[\left(i^{\xi}\left\{\hat{f}_{i-1} \left(X_i\right)-f_0\left(X_i\right)\right\}^2\right)^q\right]\right)^{1 / q}\\
        & = 2i^\xi \|\hat{f}_{i-1} \left(X\right)-f_0\left(X\right)\|_{2q}^2\\
        & =2 i^{\xi-a}\left\|i^{a/2}\left\{\hat{f}_{i-1}(X)-f_0(X)\right\}\right\|_{2 q}^2\\
        & \leq c (2^{(\theta_{\rm est} + 2)/\theta_{\rm est}}K_{\rm est})i^{\xi - a} q^{2/\theta_{\rm est}}.
    \end{aligned}
    \end{equation*}
So $D_i$ is $(\left(2^{\left(\theta_{\rm est}+2\right) / \theta_{\rm est}} K_{\rm est}\right) i^{\xi-a}, \theta_{\rm est}/2)$-sub-Weibull. This implies 
\begin{equation*}
\begin{aligned}
   \left\|\sum_{i=1}^n D_i\right\|_q & \leq c\left(\sum_{i=1}^n C\left(K_{\text {est }}, \theta_{\mathrm{est}}\right) i^{2 \xi-2 a}\right)^{1 / 2} q^{2 / \theta_{\text {est }}}\\
& \leq  cC\left(K_{\text {est }}, \theta_{\text {est }}\right) n^{\xi - a+1/2} q^{2/\theta_{\rm est}}.
\end{aligned}
\end{equation*}
Therefore
\begin{equation*}
\begin{aligned}
    & \mathbb{P}\left(\left|\sum_{i=1}^n u_i -\mathbb{E}\left[u_i  \mid F^{i-1}\right]\right| \geq c_n\right) \\
&\leq \exp(-K'(c_n n^{-\xi+a -1/2})^{\theta_{\rm est}/2}).
\end{aligned}
\end{equation*}
\end{proof}

The last deviation can be bounded as follows:
\begin{lemma}\label{lemma: infinite concentration 3}
    Under Assumptions~\ref{ass: subW noise}-\ref{ass: subW stability}, for any $c_n > 0$, we have for each $\hat{f}_{i-1}^{(k)}$:
    \begin{equation*}
        \mathbb{P}\left(\left|\sum_{i=1}^n \mathbb{E}\left[u_i^{(k)} \mid F^{i-1}\right]-\mathbb{E}\left[u_i^{(k)}\right]\right| \geq c_n\right) \leq \exp \left[-K^\prime\left\{c_n\left(n^{-\mathbb{c}-3 / 2} \log ^{-1} n \wedge \log ^{-1} n\right)\right\}^{\theta^{\prime}}\right],
    \end{equation*}
    where $\mathbb{c} = \xi - b_k - a_k/2$. And $K^\prime$, $\theta^{\prime}$ are constants depending on $K_{\rm stab}, \theta_{\rm stab}, K_{\rm est}$ and $\theta_{\rm est}$. Specifically, when $c_n=n^{1-a_k+\xi}$:
    \begin{equation*}
        \mathbb{P}\left(\left|\sum_{i=1}^n \mathbb{E}\left[u_i^{(k)} \mid F^{i-1}\right]-\mathbb{E}\left[u_i^{(k)}\right]\right| \geq c_n\right) \leq \exp \left[-K^\prime\left\{\log ^{-1} n\left(n^{b_k-a_k/2-1/ 2}  \wedge n^{1-a_k+\xi}\right)\right\}^{\theta^{\prime}}\right].
    \end{equation*}
\end{lemma}
\begin{proof} We will drop the $k$ in the superscripts. We first rearrange the terms and relate them to the stability quantity:
    \begin{equation*}
\begin{aligned}
& \sum_{i=1}^n \mathbb{E}\left[u_i \mid F^{i-1}\right]-\sum_{i=1}^n \mathbb{E}\left[u_i\right]\\
= & \sum_{i=1}^n\left(\mathbb{E}\left[u_i \mid F^{i-1}\right]-\mathbb{E}\left[u_i \mid F^{i-2}\right]+\mathbb{E}\left[u_i \mid F^{i-2}\right]-\mathbb{E}\left[u_i \mid F^{i-3}\right]+\cdots+\mathbb{E}\left[u_i \mid F^1\right]-\mathbb{E}\left[u_i\right]\right)\\
= & \sum_{i=1}^n \sum_{j=1}^{i-1} \mathbb{E}\left[u_i \mid F^j\right]-\mathbb{E}\left[u_i \mid F^{j-1}\right] \\
= & \sum_{j=1}^{n-1} \sum_{i=j+1}^n \mathbb{E}\left[u_i \mid F^j\right]-\mathbb{E}\left[u_i \mid F^{j-1}\right] =: \sum_{j=1}^{n-1} D_j,
\end{aligned}
\end{equation*}
where $D_j:= \sum_{i=j+1}^n \mathbb{E}\left[u_i \mid F^j\right]-\mathbb{E}\left[u_i \mid F^{j-1}\right]$. Note that $D_j$ is a martingale difference sequence $$\mathbb{E}[D_j \mid D_k: k\in[j-1]] =\mathbb{E}[E\{D_j \mid D_k: k\in[j-1], F^{j-1}\}] = 0.$$
Given Lemma~\ref{lemma: hoeffding sub-W}, if we can establish bounds on the moments of each $D_j$, we can control the tail of the quantity of interest.
\begin{equation}\label{eq: decompose Dj}
\begin{aligned}
& \left\|D_j\right\|_q \\
= & \left\|\sum_{i=j+1}^n \mathbb{E}\left[u_i \mid F^j\right]-\mathbb{E}\left[u_i \mid F^{j-1}\right]\right\|_q \\
\leq & \sum_{i=j+1}^n\left\|\mathbb{E}\left[u_i \mid F^j\right]-\mathbb{E}\left[u_i \mid F^{j-1}\right]\right\|_q \\
= & \sum_{i=j+1}^n\left\|\mathbb{E}\left[\nabla_j u_i \mid F^j\right]\right\|_q \\
\leq & \sum_{i=j+1}^n\left\|\nabla_j u_i\right\|_q .
\end{aligned}
\end{equation}
In equation \eqref{important decomposition for stability}, we showed that 
\begin{equation*}
    \nabla_j u_i = i^\xi \nabla_j \hat{f}_{i-1}\left(X_i\right)\left\{\hat{f}_{i-1}\left(X_i ; Z_{i-1}\right)-f_0\left(X_i\right)+\hat{f}_{i-1}\left(X_i ; Z_{i-1}^j\right)-f_0\left(X_i\right)\right\}.
\end{equation*}
Therefore 
\begin{equation*}
\begin{aligned}
& \left\|\nabla_j u_i\right\|_q \\
\leq & \left\|i^{\xi} \nabla_j \hat{f}_{i-1}\left(X_i\right)\left\{\hat{f}_{i-1}\left(X_i ; Z_{i-1}\right)-f_0\left(X_i\right)\right\}\right\|_q+ \\
& \left\|i^{\xi} \nabla_j \hat{f}_{i-1}\left(X_i\right)\left\{\hat{f}_{i-1}\left(X_i ; Z_{i-1}^j\right)-f_0\left(X_i\right)\right\}\right\|_q.
\end{aligned}
\end{equation*}
Consider the first term (the second one can be analyzed identically):
\begin{equation*}
\begin{aligned}
    & \left\|i^{\xi} \nabla_j \hat{f}_{i-1}\left(X_i\right)\left\{\hat{f}_{i-1}\left(X_i ; Z_{i-1}\right)-f_0\left(X_i\right)\right\}\right\|_q^q\\
    & \leq i^{\xi q} \sqrt{\mathbb{E}\left[\left(\nabla_j \hat{f}_{i-1}\left(X_i\right)\right)^{2 q}\right]} \sqrt{\mathbb{E}\left[\left\{\hat{f}_{i-1}\left(X_i ; Z_{i-1}\right)-f_0\left(X_i\right)\right\}^{2 q}\right]}\\
     & =  i^{(\xi - b - a/2) q} \sqrt{\mathbb{E}\left[i^{2qb}\left(\nabla_j \hat{f}_{i-1}\left(X_i\right)\right)^{2 q}\right]} \sqrt{\mathbb{E}\left[i^{aq}\left\{\hat{f}_{i-1}\left(X_i ; Z_{i-1}\right)-f_0\left(X_i\right)\right\}^{2 q}\right]}\\
     & \leq i^{\mathbb{c} q} (cK_{\rm stab}(2q)^{1/\theta_{\rm stab}})^{q} (cK_{\rm est}(2q)^{1/\theta_{\rm est}})^{q}\text{ where }\mathbb{c} = \xi-b-a / 2\\
     \Rightarrow & \left\|i^{\xi} \nabla_j \hat{f}_{i-1}\left(X_i\right)\left\{\hat{f}_{i-1}\left(X_i ; Z_{i-1}\right)-f_0\left(X_i\right)\right\}\right\|_q \\
     & \leq  i^{\mathbb{c}} c^2 (K_{\rm stab}K_{\rm est})(2q)^{\theta_{\rm stab}^{-1} + \theta_{\rm est}^{-1}}.
\end{aligned}
\end{equation*}
So we know under the assumptions $\nabla_j u_i$ is $(C(K,\theta)i^{\mathbb{c}}, (\theta_{\rm stab}^{-1} + \theta_{\rm est}^{-1})^{-1})$-sub-Weibull. Go back to \eqref{eq: decompose Dj}, we know $D_j$ is $(C(K,\theta)\sum_{i = j+1}^n i^{\mathbb{c}}, (\theta_{\rm stab}^{-1} + \theta_{\rm est}^{-1})^{-1})$-sub-Weibull as well. The summation of $D_j$, which is the quantity of interest, is $(C(K,\theta)\sqrt{\sum_{j=1}^{n-1}(\sum_{i=j+1}^n i^{\mathbb{c}}})^2, (\theta_{\rm stab}^{-1} + \theta_{\rm est}^{-1})^{-1})$-sub-Weibull, thanks to the martingale structure of $D_j$.

Apply the same argument as in the proof of Lemma~\ref{puzzle 3}, we know 
\begin{equation*}
    \mathbb{c}_n = \sqrt{\sum_{j=1}^{n-1}\left(\sum_{i=j+1}^n i^{\mathbb{c}}\right)^2} \leq c(n^{3 / 2+\mathbb{c}}\log n \vee \log n  ).
\end{equation*}
So overall, we have 
\begin{equation*}
    \begin{aligned}
        & \mathbb{P}\left(\left|\sum_{i=1}^n \mathbb{E}\left[u_i \mid F^{i-1}\right]-\mathbb{E}\left[u_i\right]\right| \geq c_n\right)\\
        = & \mathbb{P}\left(\left|\sum_{i=1}^n \mathbb{E}\left[u_i \mid F^{i-1}\right]-\mathbb{E}\left[u_i\right]\right|/\mathbb{c}_n \geq c_n/\mathbb{c}_n\right)\\
        \leq & \exp(-C'(K,\theta)(c_n/\mathbb{c}_n)^{\theta^\prime})\\
        \leq & \exp[-C'(K,\theta)\{c_n(n^{-\mathbb{c}- 3/2}\log^{-1} n \wedge \log^{-1} n)\}^{\theta^\prime}]
    \end{aligned}
\end{equation*}
\end{proof}

\newpage
\section{Technical Details for Choosing Weighting Exponents}\label{app: delay ratio}
This section offers more technical details for the discussion in \Cref{section: detection delay}. We will use the notation
$$
T(a,b,\xi) = \left(\frac{\xi+1-a}{\xi+1-b}\right)^{1/(b-a)}.
$$

\begin{lemma}\label{lemma: delay_ratio}

Consider two sequences $\alpha_i = A i^{-a}$ and $\beta_i = B i^{-b}$ for $i \geq 1$. Suppose $0 < A < B$ and $0 \leq a < b < 1$. Then $\beta_i \leq \alpha_i$ when
\[
i \geq \left(\frac{B}{A}\right)^{1/(b-a)}.
\]

Moreover, for any $\xi \geq 0$, we have
\[
\sum_{i=1}^n i^{\xi} \beta_i \leq \sum_{i=1}^n i^{\xi} \alpha_i
\]
whenever $n$ satisfies
\[
n \geq \left(\frac{B}{A}\right)^{1/(b-a)} T(a,b,\xi) + r_n,
\]
and
\[
\frac{B(\xi - a + 1)}{A n} + \frac{1}{n^{\xi - b + 1}} \leq \frac{B}{A} T^{b-a}(a, b, \xi),
\]
where
\begin{align*}
r_n &= \frac{1}{b - a} \cdot 2^{1/(b-a) - 1} \left( \frac{B}{A} \right)^{1/(b-a) - 1} 
T^{1 - (b-a)} \left( \frac{B(\xi - a + 1)}{A n} + \frac{1}{n^{\xi - b + 1}} \right) \\
&= O\left( \frac{1}{n} + \frac{1}{n^{\xi - b + 1}} \right).
\end{align*}
\end{lemma}

\begin{proof}
The first claim is straightforward to verify. We now detail the derivation of the second one. By assumption, $\xi - b > -1$. We will use the monotone–rectangle argument:
\begin{equation*}
\min_{t \in [i, i+1]} f(t) \leq \int_i^{i+1} f(t)\,dt \leq \max_{t \in [i, i+1]} f(t).
\end{equation*}

For $\beta_i$, when $\xi - b \leq 0$,
\begin{equation*}
\sum_{i=1}^{n-1} B(i+1)^{\xi - b} \leq \int_1^n B t^{\xi - b} \, dt.
\end{equation*}
So,
\begin{align*}
\sum_{i=1}^n i^{\xi} \beta_i 
&= \beta_1 + \sum_{i=2}^n B i^{\xi-b} \\
&\leq B + B \int_1^n t^{\xi-b} \, dt \\
&= B\left(1 + \frac{n^{\xi - b + 1} - 1}{\xi - b + 1}\right) \\
&\leq B\left(\frac{n^{\xi - b + 1}}{\xi - b + 1} + n^{\xi - b}\right).
\end{align*}

When $\xi - b > 0$,
\begin{align*}
\sum_{i=1}^n i^{\xi} \beta_i 
&= \sum_{i=1}^{n-1} i^{\xi} \beta_i + B n^{\xi - b} \\
&\leq B \int_1^n t^{\xi - b} \, dt + B n^{\xi - b} \\
&= B\left(\frac{n^{\xi - b + 1} - 1}{\xi - b + 1} + n^{\xi - b}\right) \\
&\leq B\left(\frac{n^{\xi - b + 1}}{\xi - b + 1} + n^{\xi - b}\right).
\end{align*}

Similarly, for $\alpha_i$, when $-1 < \xi - a \leq 0$,
\begin{align*}
\sum_{i=1}^n i^{\xi} \alpha_i 
&= \sum_{i=1}^{n-1} A i^{\xi - a} + A n^{\xi - a} \\
&\geq \int_1^n A t^{\xi - a} \, dt + A n^{\xi - a} \\
&= A\left(\frac{n^{\xi - a + 1} - 1}{\xi - a + 1} + n^{\xi - a}\right) \\
&> A \frac{n^{\xi - a + 1} - 1}{\xi - a + 1}.
\end{align*}

When $\xi - a > 0$,
\begin{align*}
\sum_{i=1}^n i^{\xi} \alpha_i 
&= A + \sum_{i=2}^n A i^{\xi - a} \\
&\geq A + A \int_1^n t^{\xi - a} \, dt \\
&= A\left(1 + \frac{n^{\xi - a + 1} - 1}{\xi - a + 1}\right) \\
&> A \frac{n^{\xi - a + 1} - 1}{\xi - a + 1}.
\end{align*}

Hence, to ensure $\sum_{i=1}^n i^{\xi} \beta_i \leq \sum_{i=1}^n i^{\xi} \alpha_i$ we just need to find a sufficient condition for
\begin{equation*}
A \frac{n^{\xi - a + 1} - 1}{\xi - a + 1} \geq B\left(\frac{n^{\xi - b + 1}}{\xi - b + 1} + n^{\xi - b}\right).
\end{equation*}
Solve this inequality:
\begin{align*}
n \geq \left(
\frac{B}{A} \cdot \frac{\xi + 1 - a}{\xi + 1 - b} 
+ \frac{B(\xi - a + 1)}{A n}
+ \frac{1}{n^{\xi - b + 1}}
\right)^{1 / (b - a)}.
\end{align*}

We will use the inequality $(x + e)^p \leq x^p + p 2^{p-1} x^{p-1} e$ for $x > 0$, $0 < e \leq x$, and $p > 1$. Define
\[
e_n := \frac{B(\xi - a + 1)}{A n} + \frac{1}{n^{\xi - b + 1}}.
\]
When $e_n < \frac{B}{A} \cdot \frac{\xi + 1 - a}{\xi + 1 - b}$, a further sufficient condition is
\begin{equation*}
n \geq \left( \frac{B}{A} \cdot \frac{\xi + 1 - a}{\xi + 1 - b} \right)^{1 / (b - a)} + r_n.
\end{equation*}
\end{proof}

%     \begin{equation*}
%             \sum_{i=1}^{n-1} i^{\xi} \beta_i   = \sum_{i=1}^{n-1} Bi^{\xi-b}\leq \int_1^{n} Bt^{\xi - b}   dt = \frac{B}{\xi - b +1}n^{\xi - b+1}.      
%     \end{equation*}
%     Meanwhile
%     \begin{equation*}
%         \sum_{i=0}^n i^\xi \alpha_i = \sum_{i=0}^n Ai^{\xi - a} \geq \int_0^{n+1} A(t-1)^{\xi - a} dt = \frac{A}{\xi - a + 1} n^{\xi - a + 1}.
%     \end{equation*}
% We just need to solve the condition for
% \begin{equation*}
%     \frac{A}{\xi-a+1} n^{\xi-a+1} \geq \frac{B}{\xi-b+1}(n+1)^{\xi-b+1} + Bn^{\xi-b}.
% \end{equation*}
% It is equivalent to
% \begin{equation*}
%     n \geq (B / A)^{(b-a)^{-1}} \left(\frac{\xi+1-a}{\xi+1-b}\right)^{1 /(b-a)} \left(1+\frac{\xi-b+1}{n}\right)^{1 /(b-a)}
% \end{equation*}

\begin{figure}[t]
    \centering
    \includegraphics[width = 0.8\textwidth]{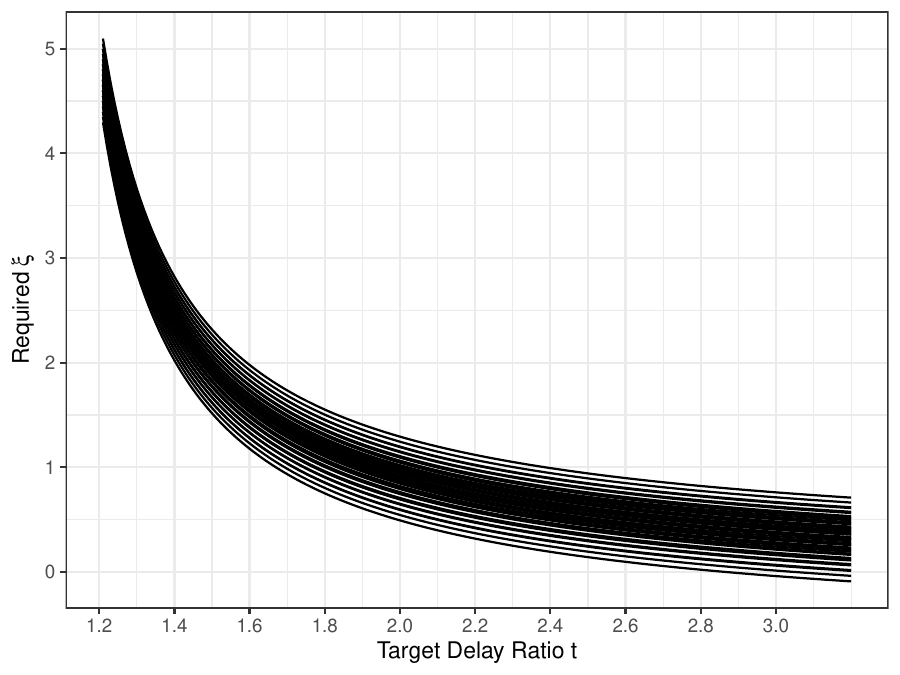}
    \caption{Function plots of $\xi(t)$, defined in the right-hand-side of \eqref{eq: threshold xi}. The $(a,b)$'s take value in $\{0,0.1,...,0.9\}^2$, subjecting to the restriction $a<b$. Each combination of $(a,b)$ corresponds to one curve in the plot.}
    \label{fig: threshold xi}
\end{figure}
The ratio between the thresholds
$$
(B / A)^{1 /(b-a)}
$$
and
$$
(B / A)^{1 /(b-a)} T(a, b, \xi)+r_n
$$
in \Cref{lemma: delay_ratio} is essentially $T(a,b,\xi)$ up to an $O(r_n)$ term. We will analyze this main component that does not vanish with an increasing $n$. The solution of $T(a,b,\xi)\leq t$ for some $t > 1$ is 
\begin{equation}\label{eq: threshold xi}
    \xi \geq \frac{1 - a - (1-b)t^{b-a}}{t^{b-a} - 1} =: \xi(t).
\end{equation}

For each fixed pair of $(a,b)$, we can treat $\xi$ as a function of the target delay ratio $t$. This function with some typical choices of the $a,b$-convergence rates is plotted in Figure~\ref{fig: threshold xi}. We can observe that in order to achieve a smaller delay ratio $t$, we need to set the diverging weighting exponent $\xi$ to be a larger value. Specifically, if we take $\xi = 1$, we can control the ratio to be less than $2.4$; if it is set to be a larger number such as $2$ or $3$, the delay can be further shortened to $1.6$ or $1.4$. 

In fact, we can also analytically analyze the above inequality of interest and get the following result:
\begin{lemma}\label{lemma: delay ratio}
    Let $0 \leq a < b < 1$ and $t > 1$ be three fixed numbers. Then we know when $\xi \geq (\log t)^{-1}$, we have 
        \begin{equation}\label{eq: eq to control for delay}
\left(\frac{\xi+1-a}{\xi+1-b}\right)^{1 /(b-a)} \leq t.
    \end{equation}
\end{lemma}
\begin{proof}
    Define $\delta = b-a$, and we rewrite the LHS as
\begin{equation*}
\begin{aligned}
& \left(\frac{\xi+1+\delta-b}{\xi+1-b}\right)^{\delta^{-1}} \\
& =\left(1+\frac{\delta}{\xi+1-b}\right)^{\delta^{-1}} \\
& \leq\left(1+\frac{\delta}{\xi}\right)^{\delta^{-1}}
\end{aligned}
\end{equation*}
Applying $(1+x^{-1})^{x} \leq e$ for $x \geq 0$. We have
\begin{equation*}
\left(\frac{\xi+1+\delta-b}{\xi+1-b}\right)^{\delta^{-1}} \leq \exp \left(\xi^{-1}\right)
\end{equation*}
Solving $\exp(\xi^{-1}) \leq t$ gives a sufficient condition for \eqref{eq: eq to control for delay}.
\end{proof}

Plugging in $t = 1.4 , 1.6$ and $2$ into $(\log t)^{-1}$ would give us $3.0$, $2.1$ and $1.4$, which roughly corresponds to the numerical values presented in Figure~\ref{fig: threshold xi}.

\subsection{Efficiency loss at $I_F$}
The impact of $\xi$ can also be examined in the axis of efficiency loss at time $I_F$. At $i = I_F = (B / A)^{1/(b-a)} \left(\frac{\xi+1-a}{\xi+1-b}\right)^{1 /(b-a)}$,
\begin{equation*}
    \mathbb{E}[r_{i,1}] = AI_F^{-a}\text{ and }\mathbb{E}[r_{i,2}] = BI_F^{-b}.
\end{equation*}
The ratio
\begin{equation*}
\begin{aligned}
    \mathbb{E}\left[r_{i,1}\right] / \mathbb{E}\left[r_{i,2}\right]
    & =\frac{A I_F^{-a}}{B I_F^{-b}}=\frac{\xi+1-a}{\xi+1-b}\\
    & = 1+\frac{b-a}{\xi+1-b}.
\end{aligned}
\end{equation*}

\newpage
\section{Stability of Estimator Examples}
\subsection{Regarding Batch Sieve Estimators}\label{app: rigorous_discussion}

The excess risk of the batch sieve estimator under $P_X$-orthonormal design can be decomposed as follows:
\begin{equation*}
    \begin{aligned}
        & \mathbb{E}\{\hat f(X) - f(X)\}^2\\
        = & \mathbb{E}\left\{\sum_{j=1}^J \hat{\beta}_j \phi_j(X)-\sum_{j=1}^J \beta_j \phi_j(X)-\sum_{j=J+1}^{\infty} \beta_j \phi_j(X)\right\}^2\\
        = & \mathbb{E}\left\{\sum_{j=1}^J (\hat{\beta}_j - \beta_j) \phi_j(X)-\sum_{j=J+1}^{\infty} \beta_j \phi_j(X)\right\}^2\\
        = & E\sum_{j=1}^J (\hat \beta_j - \beta_j)^2 + \sum_{j= J+1}^\infty \beta_j^2
    \end{aligned}
\end{equation*}

For each $j \in [J]$, the estimation error is
\begin{equation*}
    \mathbb{E}\left(\hat{\beta}_j-\beta_j\right)^2=\mathbb{E}\left(n^{-1} \sum_{i=1}^n Y_i \phi_j\left(X_i\right)-\mathbb{E}\left[Y \phi_j(X)\right]\right)^2=n^{-1} \operatorname{Var}\left(Y \phi_j(X)\right) .
\end{equation*}
Under very mild conditions on the noise distribution and basis functions, each $\mathbb{E}\left(\hat{\beta}_j-\beta_j\right)^2$ can be bounded by $n^{-1}$ times a constant. For the approximation error:
\begin{equation*}
    \sum_{j=J+1}^{\infty} \beta_j^2 \leq  J^{-2s}\sum_{j = J+1}^\infty j^{2s}\beta_k^2 = J^{-2s}Q^2,
\end{equation*}
where the constants $s, Q$ are the same as in the Sobolev ellipsoid assumptions.

In \Cref{section: batch sieve stability}, we establish the selection consistency for $\alpha < (2s+1)^{-1}$, which satisfies the stability conditions, and our convergence guarantee can cover it. Now we also demonstrate numerically that wRV can consistently rule out overfitting models ($\alpha > (2s+1)^{-1}$) as well. 

\begin{figure}[!htbp]
\centering
\includegraphics[width = \textwidth]{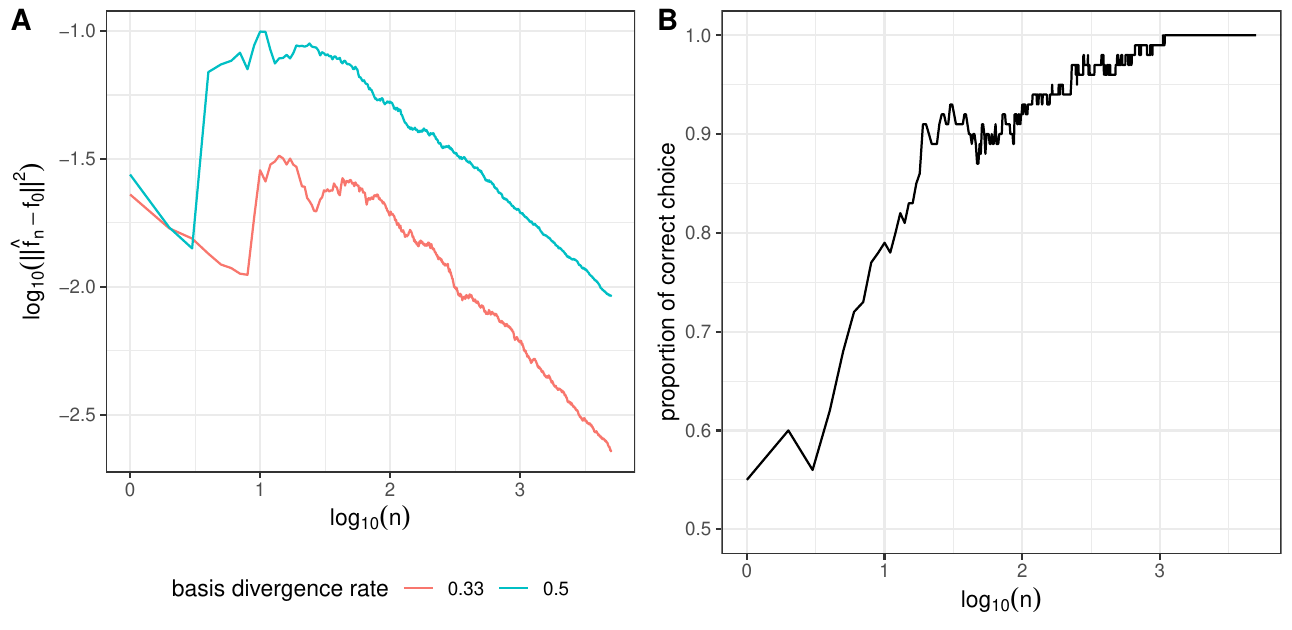}
\caption{wRV can rule out overfitting models. \textbf{(A)} The true average MSE for two sequences, each implements a different divergence rate for the number of basis functions. \textbf{(B)} The estimated probability of wRV choosing the correct sequence (i.e. favor basis divergence rate = $0.33$). Repetition $=100$ and weighting exponent $\xi = 1$.}
\label{fig: batch_sieve_consistency}
\end{figure}

We conduct a numerical experiment under the following setting:
\begin{equation*}
    \begin{aligned}
        X_i & \sim \mathcal{U}([0,1]),\nonumber\\
       \epsilon_i & \sim \mathcal{N}(0, 0.1^2),\nonumber\\
       Y_i & = \mathcal{K}(0.75, X_i)+\epsilon_i,
    \end{aligned}
\end{equation*}
where the function
\begin{equation*}
    \begin{aligned}
        \mathcal{K}(a,b) & = \sinh(1)^{-1}\cdot \cosh(a\wedge b)) \cdot \cosh(1 - a \vee b)\\
        & = \sum_{j=1}^\infty \lambda_j \psi_j(a)\psi_j(b).
    \end{aligned}
\end{equation*}
The eigenvalues are $\lambda_1 = 1$, $\lambda_j=(1+(j-1)^2 \pi^2)^{-1}$ for $j \geq 2$, and eigenfunctions $\psi_j(x)$ are the cosine basis:
\begin{equation*}
    \psi_1(x) = 1,\quad \psi_j(x) = \sqrt{2}\cos((j-1) 2\pi x).
\end{equation*}
The estimators are also constructed with the same basis functions. We consider two $J_n = n^{0.33}$ (optimal choice) and $J_n = n^{0.5}$. The latter overfits the data (\Cref{fig: batch_sieve_consistency}, A) and has worse stability. However, wRV can still choose the superior sequence with probability converges to $1$ (\Cref{fig: batch_sieve_consistency}, B), which is more optimistic than the current theoretical guarantees.

\subsection{Proof of Theorem~\ref{th: stability of parametric sgd}}
\label{app:stability of SGD}

\begin{proof}
Recall the iteration formula of the model parameter:
\begin{equation*}
   \begin{aligned}
 \beta_0&=0 \in \mathbb{R}^p \\
 \beta_i&=\beta_{i-1}+\gamma\left(Y_i-X_i^{\top} \beta_{i-1}\right) X_i\\
& = \gamma Y_i X_i + (I - \gamma X_i X_i^\top) \beta_{i-1}
\end{aligned} 
\end{equation*}
In our stability notation, $\beta_i$ can be written as $\beta_i(Z_i)$, meaning it is calculated using the true sample $Z_i = \{(X_1,Y_1),...,(X_i,Y_i)\}$. And we use $\beta'_i = \beta'_i(Z_i^j)$ to denote the model parameter trained with sample $Z_i^j$, whose $j$-th sample is a IID copy of $(X_j,Y_j)$ in $Z_i$. We use $\eta_i$ to denote the difference between $\beta_i$ and $\beta'_i$
\begin{equation*}
    \eta_i = \beta_i - \beta'_i.
\end{equation*}
We are going to show $\eta_i$'s are vectors of small norm. It is trivial that $\eta_i = 0$ when $j>i$ (the parameters are calculated using the same samples).
When $i = j$, 
\begin{equation*}
    \eta_j = \gamma Y_i X_i+\left(I-\gamma X_i X_i^{\top}\right) \beta_{i-1} - \gamma Y'_i X'_i+\left(I-\gamma X'_i X_i^{\prime\top}\right) \beta_{i-1}
\end{equation*}

For $i \geq j+1$, we have the iteration formula:
\begin{align*}
    \eta_i
    & = (I - \gamma \|X_i\|^2 H_i)\eta_{i-1}\\
    & = H_i^{\perp} \eta_{i-1} + (1 - \gamma \|X_i\|^2) H_i \eta_{i-1},
\end{align*}
where we used $H_i = X_i(\|X_i\|^2)^{-1}X_i^{\top}$ to denote the projection matrix onto the direction of $X_i$ and $H_i^{\perp} = I - H_i$ is the projection onto its orthogonal complement.

Then we have
\begin{align*}
    \|\eta_i\|^2 
    & = \|H_i^{\perp} \eta_{i-1}\|^2 + (1 - \gamma \|X_i\|^2)^2 \|H_i \eta_{i-1}\|^2\\
    &\stackrel{(I)}{\leq} \|H_i^{\perp} \eta_{i-1}\|^2 + (1 - \gamma \|X_i\|^2) \|H_i \eta_{i-1}\|^2\\
    &\leq \|\eta_{i-1}\|^2 - \gamma \|X_i\|^2\|H_i \eta_{i-1}\|^2
\end{align*}
In step $(I)$ we used $\gamma \|X_i\|^2 \leq 1$. Therefore,
\begin{equation*}
    \frac{\|\eta_i\|^2}{\|\eta_{i-1}\|^2} 
    \leq 1 - \gamma \|X_i\|^2 \frac{\|H_i \eta_{i-1}\|^2}{\|\eta_{i-1}\|^2} 
    = 1 - \gamma \frac{\eta_{i-1}^{\top} X_i X_i^{\top} \eta_{i-1}} {\|\eta_{i-1}\|^2}
\end{equation*}

Take conditional expectations on both sides:
\begin{equation*}
\begin{aligned}
        \mathbb{E}\left[\left.\frac{\|\eta_i\|^2}{\|\eta_{i-1}\|^2} \right| F^{i-1}\right] 
    & \leq 1 - \gamma \|\eta_{i-1}\|^{-2} \eta_{i-1}^{\top} \mathbb{E}[X_i X_i^{\top}] \eta_{i-1}\\
    & \stackrel{(II)}{\leq} 1 - \gamma \underline{\lambda} \leq \exp(-\gamma \underline{\lambda}).
\end{aligned}
\end{equation*}

In step $(II)$ we used $\lambda_{min} (\mathbb{E}[X_iX_i^{\top}]) \geq \underline{\lambda} > 0$. Now we have established some ``exponential contraction" properties for the sequence of vectors $\{\eta_i\}$. Then we can apply Lemma~\ref{summable lemma using gamma function} (with parameter $a$ therein equals to $0$) to conclude that the averaged coefficient vector $\hat\eta_i = i^{-1}\sum_{k=1}^i\eta_k$ has magnitude:
\begin{equation}
\label{eq: parametric magnitude}
     \mathbb{E}[\|\hat\eta_i\|^2 | F^j] \lesssim i^{-2}\|\eta_j\|^2.
\end{equation}

Recall that our goal is to control the stability of the estimator. Use the bound on $\hat\eta_i$ we have:
\begin{equation*}
\begin{aligned}
    \mathbb{E}\left[\left\{\nabla_j \hat{f}_{i-1}\left(X_i\right)\right\}^2 \mid F^j\right]
    & = \mathbb{E}\left[\left\{\hat{f}_{i-1}\left(X_i ; Z_i\right)-\hat{f}_{i-1}\left(X_i ;  Z_i^j\right)\right\}^2 \mid F^j\right]\\
    & = \mathbb{E}\left[\left\{X_i^{\top} \hat{\beta}_{i-1}-X_i^{\top} \hat{\beta}_{i-1}^{\prime}\right\}^2 \mid F^j\right]\\
    & \text{ (where }\hat\beta_i = i^{-1}\sum_{k=1}^i \beta_k,\ \hat\beta'_i = i^{-1}\sum_{k=1}^i \beta'_k) \\
    & = \mathbb{E}\left[\left. \{ X_i^{\top} \hat\eta_{i-1}\}^2\right| F^j\right]\\
    & \leq \mathbb{E}\left[\left. \|X_i\|^2 \|\hat\eta_{i-1}\|^2\right| F^j\right],
\end{aligned}
\end{equation*}

Under the assumption that $\mathbb{E}[\|X_i\|^2]\leq R^2$,
\begin{equation}\label{eq:martingaleoneta}
\begin{aligned}
    \mathbb{E}\left[\left\{\nabla_j \hat{f}_{i-1}\left(X_i\right)\right\}^2 \mid F^j\right] &\leq R^2 \mathbb{E}\left[\left\|\hat{\eta}_{i-1}\right\|^2 \mid F^j\right]\\
    & \stackrel{(III)}{\lesssim} i^{-2}\|\eta_j\|^2
\end{aligned}
\end{equation}
In step $(III)$ we used \eqref{eq: parametric magnitude}. Under the uniform assumptions: $\|X_j\|^2 \leq R^2$ and $|Y_j - X_j^{\top}\beta_{j-1}|\leq M$, we have
\begin{equation}
\begin{aligned}
\label{eq:boundetaj}
   \left\|\eta_j\right\| & \leq\left\|\gamma X_j\left(Y_j-X_j^{\top} \beta_{j-1}\right)\right\|+\left\|\gamma X_j^{\prime}\left(Y_j^{\prime}-X_j^{\prime \top} \beta_{j-1}\right)\right\|\\
    & \leq 2\gamma M R\\
    \Rightarrow 
    \|\eta_j\|^2 
    &\leq 4\gamma^2 M^2 R^2.
\end{aligned}
\end{equation}

Combine \eqref{eq:martingaleoneta} and \eqref{eq:boundetaj}, we know
\begin{equation*}
    \mathbb{E}\left[\left\{\nabla_j \hat{f}_{i-1}\left(X_i\right)\right\}^2 \mid F^j\right] \lesssim i^{-2},
\end{equation*}
which concludes our proof.
\end{proof}

\subsection{Proof of Theorem~\ref{th: stability sieve with shrinkage}}
\label{app:stability of Sieve SGD}
%%%%%%%%%%%%

\begin{proof}
Let $\beta_i$ be the coefficient vectors trained with sample $Z_i$ and $\beta_i^{\prime}$ from $Z_i^j$. At the $j$-th step, $\beta_i$ is updated using sample $\left(X_j, Y_j\right)$:
\begin{equation*}
\beta_j=\beta_{j-1}^{\wedge}+\gamma_j\left(Y_j-\boldsymbol{\phi}_j^{\top} \beta_{j-1}^{\wedge}\right) D_j \phi_j,
\end{equation*}
And we have a similar update rule for $\beta_j^{\prime}$ using sample $\left(X_j^{\prime}, Y_j^{\prime}\right)$ :
\begin{equation*}
\begin{aligned}
\beta_j^{\prime} & =\beta_{j-1}^{\prime \wedge}+\gamma_j\left(Y_j-\boldsymbol{\phi}_j^{\prime \top} \beta_{j-1}^{\prime \wedge}\right) D_j \boldsymbol{\phi}_j^{\prime} \\
& =\beta_{j-1}^{\wedge}+\gamma_j\left(Y_j-\boldsymbol{\phi}_j^{\prime \top} \beta_{j-1}^{\wedge}\right) D_j \boldsymbol{\phi}_j^{\prime}
\end{aligned}
\end{equation*}
Here $\boldsymbol{\phi}_j^{\prime}=\left(\phi_1\left(X_j^{\prime}\right), \ldots, \phi_{J_i}\left(X_j^{\prime}\right)\right)^{\top}$.

Let $\eta_i=\beta_i-\beta_i^{\prime}$ be the difference between the two vectors. We have $\eta_i=0$ for $i<j$. And
\begin{equation*}
\eta_j=\gamma_j\left(Y_j-\boldsymbol{\phi}_j^{\top} \beta_{j-1}^{\wedge}\right) D_j \boldsymbol{\phi}_j-\gamma_j\left(Y_j-\boldsymbol{\phi}_j^{\prime \top} \beta_{j-1}^{\prime \wedge}\right) D_j \boldsymbol{\phi}_j^{\prime} .
\end{equation*}

From the recursive formula of $\beta_i,\beta'_i$ we can also derive one for $\eta_i$:
\begin{equation*}
\eta_i=\eta_{i-1}^{\wedge}-\gamma_i D_i \boldsymbol{\phi}_i \boldsymbol{\phi}_i^{\top} \eta_{i-1}^{\wedge}
\end{equation*}

Multiply $D_i^{-1/2}$ on both sides:
\begin{equation*}
D_i^{-1 / 2} \eta_i =D_i^{-1 / 2} \eta_{i-1}^{\wedge}-\gamma_i D_i^{1 / 2} \boldsymbol{\phi}_i\left(D_i^{1 / 2} \boldsymbol{\phi}_i\right)^{\top} D_i^{-1 / 2} \eta_{i-1}^{\wedge}
\end{equation*}

Denote $\theta_i=D_i^{-1/2} \eta_i, W_i=D_i^{1/2} \boldsymbol{\phi}_i$, we rewrite the above equation as:
\begin{equation}
\label{eq:theta recursive}
\theta_i=\theta^\wedge_{i-1}-\gamma_i W_i W_i^{\top} \theta^\wedge_{i-1}
\end{equation}

For the rest of the proof, we aim to derive some bounds on $\|\theta_i\|^2$ (and their averaged version $\hat\theta_i = i^{-1}\sum_{k=1}^i \theta_{\rm est}$).

For $\theta_j$ we have the following bound:
\begin{equation}
\label{eq: bound on theta j}
\begin{aligned}
\left\|\theta_j\right\| & =\left\|\gamma_j\left(Y_j-\boldsymbol{\phi}_j^{\top} \beta_{j-1}^{\wedge}\right) D_j^{1 / 2} \boldsymbol{\phi}_j-\gamma_j\left(Y_j-\boldsymbol{\phi}_j^{\prime \top} \beta_{j-1}^{\prime \wedge}\right) D_j^{1 / 2} \boldsymbol{\phi}_j^{\prime}\right\| \\
& \leq M \gamma_j\left(\left\|D_j^{1 / 2} \boldsymbol{\phi}_j\right\|+\left\|D_j^{1 / 2} \boldsymbol{\phi}_j^{\prime}\right\|\right) \stackrel{(I)}{\lesssim} \gamma_j
\end{aligned}
\end{equation}

In step $(I)$ above we used the norm of $D_j^{1/2} \boldsymbol{\phi}_j$ is uniformly bounded:
\begin{equation*}
    \|D_j^{1/2} \boldsymbol{\phi}_j\|^2 = \sum_{k=1}^{J_j} (k^{-\omega}\phi_k(X_j))^2 \leq C\sum_{k=1}^\infty k^{-2\omega} < \infty,
\end{equation*}
recall that $\omega > 1/2$. Now we are going to use the iteration formula to derive bounds for $\|\theta_i\|, i > j$. Take inner product on both sides of \eqref{eq:theta recursive}:
\begin{equation*}
\begin{aligned}
    \left\|\theta_i\right\|^2
    & =\left\|\theta_{i-1}^{\wedge}\right\|^2+\gamma_i^2\left(\theta_{i-1}^{\wedge}\right)^{\top} W_i W_i^{\top} W_i W_i^{\top} \theta_{i-1}^{\wedge}-2 \gamma_i\left(\theta_{i-1}^{\wedge}\right)^{\top} W_i W_i^{\top} \theta_{i-1}^{\wedge}\\
    & \stackrel{(II)}{\leq} \left\|\theta^\wedge_{i-1}\right\|^2-\gamma_i \theta_{i-1}^{\wedge\top} W_i W_i^{\top} \theta^\wedge_{i-1}.
\end{aligned}
\end{equation*}
In step $(II)$ above we used $\gamma_i\|W_i\|^2\leq 1$. Take the conditional expectation
\begin{equation*}
    \begin{aligned}
    \mathbb{E}[\|\theta_i\|^2 | F^{i-1}] 
    & \leq \|\theta^\wedge_{i-1}\|^2 - \gamma_i \theta_{i-1}^{^\wedge\top} \mathbb{E}[W_i W_i^{\top}] \theta^\wedge_{i-1}\\
    & = \|\theta_{i-1}\|^2 - \gamma_i \theta_{i-1}^{^\wedge\top} D_i^{1/2}\mathbb{E}[\boldsymbol{\phi}_i \boldsymbol{\phi}_i^\top]D_i^{1/2} \theta^\wedge_{i-1}
    \end{aligned}
\end{equation*}
Divide $\|\theta_{i-1}\|^2 =\|\theta^\wedge_{i-1}\|^2 $ on both sides:
\begin{equation*}
    \mathbb{E}[\|\theta_i\|^2/\|\theta_{i-1}\|^2 | F^{i-1}] \leq 1 - \gamma_i \frac{\theta_{i-1}^{^\wedge\top} D_i^{1/2}\mathbb{E}[\boldsymbol{\phi}_i \boldsymbol{\phi}_i^\top]D_i^{1/2} \theta^\wedge_{i-1}}{\|\theta^\wedge_{i-1}\|^2}
\end{equation*}

By our assumption, $\lambda_{min}(\mathbb{E}[\boldsymbol{\phi}_i \boldsymbol{\phi}_i^\top])$ is strictly larger than a constant (we postpone the proof after equation \eqref{eq:finilizing sieve sgd}). Also note that $D_i^{1/2} = {\rm Diag}(1^{-\omega},...,J_i^{-\omega})$ has eigenvalues greater than $J_i^{-\omega} = i^{-\zeta\omega}$. So overall we have
\begin{equation*}
    \frac{\theta_{i-1}^{^\wedge\top} D_i^{1/2}\mathbb{E}[\boldsymbol{\phi}_i \boldsymbol{\phi}_i^\top]D_i^{1/2} \theta^\wedge_{i-1}}{\|\theta^\wedge_{i-1}\|^2} \geq i^{-2\zeta\omega}.
\end{equation*}

Then we have an exponential contraction recursion:
\begin{equation*}
    \mathbb{E}[\|\theta_i\|^2/\|\theta_{i-1}\|^2 | F^{i-1}] \leq 1 - \gamma_i i^{-2\zeta\omega} = 1 -  i^{-2\zeta\omega-\zeta} \leq \exp(-i^{-\alpha}),
\end{equation*}
for ease of notation, we defined $\alpha = 2\zeta\omega+\zeta$. Now we are ready to apply Lemma~\ref{summable lemma using gamma function} to bound the magnitude of $\hat\theta_i = i^{-1}\sum_{k=1}^i \theta_{\rm est}^\wedge$ as:
\begin{equation*}
    \mathbb{E}[\|\hat\theta_i\|^2|F^j] \lesssim i^{-2}j^{2\alpha}\|\theta_j\|^2.
\end{equation*}

Plug this into \eqref{eq: bound on theta j}, we get 
\begin{equation*}
    \mathbb{E}[\|\hat\theta_i\|^2|F^j] \lesssim i^{-2}j^{2\alpha-2\zeta} = i^{-2}j^{4\omega\zeta}
\end{equation*}

This implies the stability property of the sieve SGD estimator, let $\hat\beta_i = i^{-1}\sum_{k=1}^i \beta_k^\wedge$ denote the averaged coefficient vector ($\hat{\beta}_i^{\prime}, \hat{\eta}_i$ are similarly defined):
\begin{equation}
\label{eq:finilizing sieve sgd}
\begin{aligned}
\mathbb{E}\left[\left\{\nabla_j \hat{f}_{i-1}\left(X_i\right)\right\}^2 \mid F^j\right] & =\mathbb{E}\left[\left\{\hat{f}_{i-1}\left(X_i ; Z_i\right)-\hat{f}_{i-1}^{\prime}\left(X_i ; Z_i^j\right)\right\}^2 \mid F^j\right] \\
& =\mathbb{E}\left[\left\{\left\langle\hat{\beta}_{i-1}^{\wedge}, \boldsymbol{\phi}_i\right\rangle-\left\langle\hat{\beta}_{i-1}^{\prime \wedge}, \boldsymbol{\phi}_i\right\rangle\right\}^2 \mid F^j\right] \\
& =\mathbb{E}\left[\left(\hat{\eta}_{i-1}^{\wedge \top} \boldsymbol{\phi}_i\right)^2 \mid F^j\right] \\
& =\mathbb{E}\left[\left(\hat{\theta}_{i-1}^{\wedge \top} W_i\right)^2 \mid F^j\right] \\
& =\mathbb{E}\left[\hat{\theta}_{i-1}^{\top} \mathbb{E}\left[W_i W_i^{\top}\right] \hat{\theta}_{i-1} \mid F^j\right] \\
& \lesssim \lambda_{\max }\left(\mathbb{E}\left[W_i W_i^{\top}\right]\right) i^{4 \omega \zeta-2}.
\end{aligned}
\end{equation}

To finalize the proof, we just need to bound the minimal (maximal) eigenvalue of matrix $\mathbb{E}[\boldsymbol{\phi}_i \boldsymbol{\phi}_i^\top]$ ($\mathbb{E}[W_iW_i^\top]$). We denote the density of $X$ (with respect to the measure $\nu$) as $p_X$. By our assumption there exists some positive constants $\bar p_X, \underline{p}_X$ such that $0 < \underline{p}_X \leq p_X(x) \leq \bar p_X < \infty$. Let $v$ denote any eigenvector of $\mathbb{E}[\boldsymbol{\phi}_i \boldsymbol{\phi}_i^\top]$, we have:
\begin{equation*}
      v^\top \mathbb{E}[\boldsymbol{\phi}_i \boldsymbol{\phi}_i^\top]v 
      = 
    \mathbb{E}\left[\left(\sum_{k=1}^{J_i}\phi_k(X_i)v_k\right)^2\right] \geq \underline{p}_X \int \left(\sum_{k=1}^{J_i}\phi_k(X_i)v_k\right)^2 d\nu = \underline{p}_X \|v\|^2.
\end{equation*}
So we know $\lambda_{min}(\mathbb{E}[\boldsymbol{\phi}_i \boldsymbol{\phi}_i^\top]) \geq \underline{p}_X$. Similarly, let $u$ denote any eigenvector of $\mathbb{E}[W_iW_i^\top]$:
\begin{equation*}
    u^\top \mathbb{E}[W_iW_i^\top]u = \mathbb{E}\left[\left(\sum_{k=1}^{J_i}k^{-\omega}\phi_k(X_i)u_k\right)^2\right] \lesssim \bar p_X \|u\|^2.
\end{equation*}
So we know $\lambda_{max}(\mathbb{E}[W_iW_i^\top]) \lesssim \bar p_X $. Combining these eigenvalue results with \eqref{eq:finilizing sieve sgd}, we conclude that 
$$ \mathbb{E}\left[\left\{\nabla_j \hat f_{i-1}\left(X_i\right)\right\}^2 \mid F^j\right] \lesssim i^{4\omega\zeta - 2}$$.
\end{proof}

\subsection{A Technical Lemma}
The following lemma was used in the proof of the stability properties of SGD estimators.

\begin{lemma}
\label{summable lemma using gamma function}
Let $\{\eta_i, i\in \mathbb{N}^+\}$ be a sequence of real vectors. Assume that there exists some $j\in \mathbb{N}^+$ and some $a \in [0,1), A>0$ such that:
\begin{itemize}
    \item $\eta_i = 0$ when $i \in [j-1] $;
    \item For all $i \geq j$, the sequence satisfies the following step-wise ``exponential contraction" condition:
\begin{equation}
\label{eq:exponential contraction for technical lemma}
\begin{aligned}
     \mathbb{E}[\|\eta_{i+1}\|^2/\|\eta_i\|^2 | F^i] &
     \leq \exp(-Ai^{-a}).
\end{aligned}
\end{equation}
\end{itemize}

Then for $\hat\eta_n = n^{-1}\sum_{i=1}^n \eta^\wedge_i$, we have
\begin{equation*}
    n^2\mathbb{E}[\|\hat\eta_n\|^2 | F^j] \leq C(a,A)j^{2a}\|\eta_j\|^2,
\end{equation*}
with some constant $C(a,A)$ depending on $a,A$.
\end{lemma}

\begin{remark}
   the $\eta_i$ vectors in Lemma~\ref{summable lemma using gamma function} may belong to the same real-vector space $\mathbb{R}^J$ for a fixed, positive integer $J$. It is also possible that they belong to vector spaces of different dimensions, that is, $\eta_i\in\mathbb{R}^{J_i}$ where $J_i$ is a sequence of positive integers. The latter case is especially important when we prove the stability of Sieve-type SGD estimators (Theorem~\ref{th: stability sieve with shrinkage}).
\end{remark}

\begin{proof}
    Our readers can take $a = 0$ for the first pass which alleviates the technical complexity. For the ease of notation, we use $\mathbb{E}_j$ to denote the conditional expectation $\mathbb{E}[\cdot | F^j]$.
    
    By our assumption \eqref{eq:exponential contraction for technical lemma}, we have that for any $k > i \geq j$:
    \begin{equation*}
    \begin{aligned}
    \mathbb{E}_{k-1}\left[\|\eta_k\|^2/\|\eta_i\|^2\right]
    & = \mathbb{E}_{k-1}\left[ \frac{\|\eta_{k-1}\|^2}{\|\eta_i\|^2} \frac{\|\eta_{k}\|^2}{\|\eta_{k-1}\|^2}\right]
    = \frac{\|\eta_{k-1}\|^2}{\|\eta_{i}\|^2} \mathbb{E}_{k-1}\left[ \frac{\|\eta_{k}\|^2}{\|\eta_{k-1}\|^2} 
    \right]\\
    & \leq \exp(-A(k-1)^{-a}) \|\eta_{k-1}\|^2/\|\eta_{i}\|^2.
    \end{aligned}
    \end{equation*}
    
    Iterate the above argument we have:
    \begin{equation*}
    \mathbb{E}_i[\|\eta_k\|^2/\|\eta_i\|^2] 
     \leq \exp\left(-A\sum_{l = i}^{k-1} l^{-a}\right).
    \end{equation*}
    
    Apply Jensen's inequality, we also have:
    \begin{equation*}
    \begin{aligned}
     \mathbb{E}_i[\|\eta_k\|/\|\eta_i\|] 
     & \leq \sqrt{\mathbb{E}_i[\|\eta_k\|^2/\|\eta_i\|^2]} \leq \exp\left(-A\sum_{l = i}^{k-1} l^{-a}/2\right)\\
    \Rightarrow \mathbb{E}_i[\|\eta_k\|] & \leq \exp \left(-A \sum_{l=i}^{k-1} l^{-a} / 2\right) \|\eta_i\|.
    \end{aligned}
    \end{equation*}
    
    Expand the quantity of interest and plug in the above iteration equations we can derive our results:
    \begin{equation}
    \label{skeleton of sieve sgd stability}
        \begin{aligned}
            n^2 \mathbb{E}_j[\|\hat \eta_n \|^2]
            & = \mathbb{E}_j\left[\left\|\sum_{i=j}^n \eta_i\right\|^2\right]\quad \text{ ($\eta_i = 0$ for $i < j$)}\\
            & \leq 2\sum_{i=j}^n \sum_{k=i}^n \mathbb{E}_j[\eta_k^\top \eta_i]\\
            & \leq 2\sum_{i=j}^n \sum_{k=i}^n \mathbb{E}_j[\|\eta_k\|\|\eta_i\|] \\
            & \leq 2\sum_{i=j}^n \sum_{k=i}^n \exp\left(-A\sum_{l = i}^{k-1} l^{-a}/2\right)\mathbb{E}_j[\|\eta_i\|^2] \\
            & \stackrel{(I)}{\leq} C(a,A)\sum_{i = j}^n i^a \mathbb{E}_j[\|\eta_i\|^2] \\
            & \leq C(a,A)\sum_{i = j}^n i^a \exp\left(-\sum_{l = j}^{i-1} l^{-a}\right)\|\eta_j\|^2\\
            & \stackrel{(II)}{\leq} C(a,A)j^{2a}\|\eta_j\|^2.
        \end{aligned}
    \end{equation}

%%%%%

Steps $(I)$ and $(II)$ are technical; we will present them below. 
% Under the assumption that $Y_i - X_i^\top \beta_{i-1}$ is bounded almost surely, we have $\|\eta_j\|^2 \leq \gamma_j^2 \|X_j\|^2 \leq j^{-2a}j^a = j^{-a}$. Recall that our choice of $\gamma_i = i^{-a}$. 

For step $(I)$, we need to show $\sum_{k=i}^n \exp(-\sum_{l = i}^{k-1} A l^{-a}) \leq C(a,A)i^a$. Each term in the summation has the following bound:
\begin{equation*}
        \exp\left(-\sum_{l=i}^{k-1} Al^{-a}\right) 
        \leq \exp(-A(1-a)^{-1}\{(k-1)^{1-a} - i^{1-a}\}).
\end{equation*}
Then we have
\begin{equation}
\label{first incomplete gamma}
    \begin{aligned}
        \sum_{k = i}^n \exp\left(-\sum_{l=i}^{k-1} Al^{-a}\right) 
        &\leq \exp(A(1-a)^{-1}i^{1-a}) \sum_{k=i}^n \exp(-A(1-a)^{-1}(k-1)^{1-a}) \\
        &\leq  \exp(A(1-a)^{-1}i^{1-a}) \int_i^{n+1} \exp(-A(1-a)^{-1}(k-1)^{1-a}) dk
    \end{aligned}
\end{equation}

%%%%
Define $u = A(1-a)^{-1}(k-1)^{1-a}$, we have $du = A(k-1)^{-a}dk$ and $(k-1)^{-a} = \{A^{-1}(1-a)u\}^{-a/(1-a)}$.
Therefore,
\begin{equation*}
    \begin{aligned}
    du & = A\{A^{-1}(1-a)u\}^{-a/(1-a)} dk\\
    \Rightarrow dk & = A^{-1}\{A^{-1}(1-a)u\}^{a/(1-a)} du
    \end{aligned}
\end{equation*}
Make the change of variable in \eqref{first incomplete gamma}, we simplify the integral as
\begin{equation}
\label{my upper incomplete gamma}
    \int_{A(1- a)^{-1}(i-1)^{1-a}}^{A(1- a)^{-1}n^{1-a}} \exp(-u) A^{-1}\{A^{-1}(1-a)u\}^{a/(1-a)} du
\end{equation}
Note that the above integral is essentially the upper incomplete gamma function. We have the following upper bound of its magnitude:

\begin{lemma}(Theorem~4.4.3 of \cite{gabcke1979neue})
\label{lemma: incomplete gamma}
Let $\Gamma(m,x)$ denote the upper incomplete gamma function:
\begin{equation*}
    \Gamma(m,x)=\int_x^\infty \exp(-u) u^{m-1} du,
\end{equation*}
with $x > m \geq 1$, then we have the upper bound
\begin{equation*}
    \Gamma(m,x) \leq m\exp(-x) x^{m-1}.
\end{equation*}
\end{lemma}

%%%%%

Applying Lemma~\ref{lemma: incomplete gamma}, we can continue \eqref{my upper incomplete gamma} as:
\begin{equation*}
    \begin{aligned}
    \label{good term}
    & \int_i^{n+1} \exp(-A(1-a)^{-1}k^{1-a}) dk\\
    & \leq A^{-1-a/(1-a)}(1-a)^{a/(1-a)}\Gamma((1-a)^{-1}, A(1- a)^{-1}(i-1)^{1-a}) \\
    & \leq A^{-1-a/(1-a)}(1-a)^{a/(1-a)-1} \exp(-A(1- a)^{-1}(i-1)^{1-a}) (A(1- a)^{-1}(i-1)^{1-a})^{a/(1-a)}\\
    & \leq A^{-1}(1-a)^{-1}\exp(-A(1- a)^{-1}(i-1)^{1-a}) i^a.
    \end{aligned}
\end{equation*}
Combining \eqref{first incomplete gamma} and \eqref{good term}, we conclude the proof of step $(I)$ in \eqref{skeleton of sieve sgd stability}.

The step $(II)$ of \eqref{skeleton of sieve sgd stability} can be shown using a similar argument.
\begin{equation}
\label{eq: wher the second integral shows up}
    \begin{aligned}
        \sum_{i=j}^n i^a \exp \left(-\sum_{l=j}^{i-1} Al^{-a}\right) 
        & \leq \exp(A(1-a)^{-1}j^{1-a}) \sum_{i=j}^n \exp(-A(1-a)^{-1}(i-1)^{1-a}) i^a \\
        & \leq \exp(A(1-a)^{-1}j^{1-a}) \int_j^{n+1} \exp(-A(1-a)^{-1}(i-1)^{1-a}) i^a di
    \end{aligned}
\end{equation}
Change the variable in the integral: $u = A(1-a)^{-1}(i-1)^{1-a}$, $du = A(i-1)^{-a} di$, $di = A^{-1}\{A^{-1}(1-a)u\}^{a/(1-a)}du$, so we can bound the integral in \eqref{eq: wher the second integral shows up} by
\begin{equation*}
\begin{aligned}
   &\int_{A(1-a)^{-1}(j-1)^{1-a}}^{\infty} \exp(-u) i^a \{(1-a) u \}^{a/(1-a)} du\\
 & \leq 2^a A^{-1-2a/(1-a)}(1-a)^{2a/(1-a)}\int_{A(1-a)^{-1}(j-1)^{1-a}}^{\infty}\exp(-u) u^{2a/(1-a)} du\\
 & = 2^aA^{-1-2a/(1-a)}(1-a)^{2a/(1-a)}\Gamma((1+a)/(1-a), A(1-a)^{-1}(j-1)^{1-a}) \\
 & \leq 2^aA^{-1-2a/(1-a)}(1-a)^{2a/(1-a)} (1+a)/(1-a)\\
 &\qquad \exp(-A(1-a)^{-1}(j-1)^{1-a}) \{A(1-a)^{-1}(j-1)^{1-a}\}^{2a/(1-a)}\\
 & \leq 2^aA^{-1}(1+a)(1-a)^{-1} \exp(-A(1-a)^{-1}(j-1)^{1-a}) j^{2a}.
\end{aligned}
\end{equation*}
\end{proof}
\newpage
\section{Supplementary Numerical Studies and More Detail}
\subsection{Choice of Hyperparameters}\label{app:details_numerica}
In Section~\ref{section: numerical}, \Cref{exa:2}, we described performing model selection among $8$ sieve-SGD estimators. Here we list their hyperparameters. Recall that $\gamma_i=A i^{-(2 s+1)^{-1}}$ is the learning rate and $J_i=B i^{(2 s+1)^{-1}}$ is the number of basis function, as defined in \eqref{eq: step size and basis number}. For all eight methods, $\omega$ in \eqref{eq: sieve sgd update} is 0.51.
\begin{table}[ht]
\centering
\begin{tabularx}{0.7\textwidth}{XXXX}
\hline
Model Index & $s$ & $A$ & $B$ \\
\hline
1 & 1 & 0.1 & 2 \\
\hline
2 & 2 & 0.1 & 2 \\
\hline
3 & 1 & 1 & 2 \\
\hline
4 & 2 & 1 & 2 \\
\hline
5 & 1 & 0.1 & 8 \\
\hline
6 & 2 & 0.1 & 8 \\
\hline
7 & 1 & 1 & 8 \\
\hline
8 & 2 & 1 & 8 \\
\hline
\end{tabularx}
\caption{Hyperparameter for the estimators of \Cref{exa:2}.}
\label{tab:parameter example 2}
\end{table}

\subsection{Numerical Comparison using Larger Exponents $\xi$}\label{app: extra numerical}

\begin{figure}[!t]
    \centering
    \includegraphics[width = 0.8\textwidth]{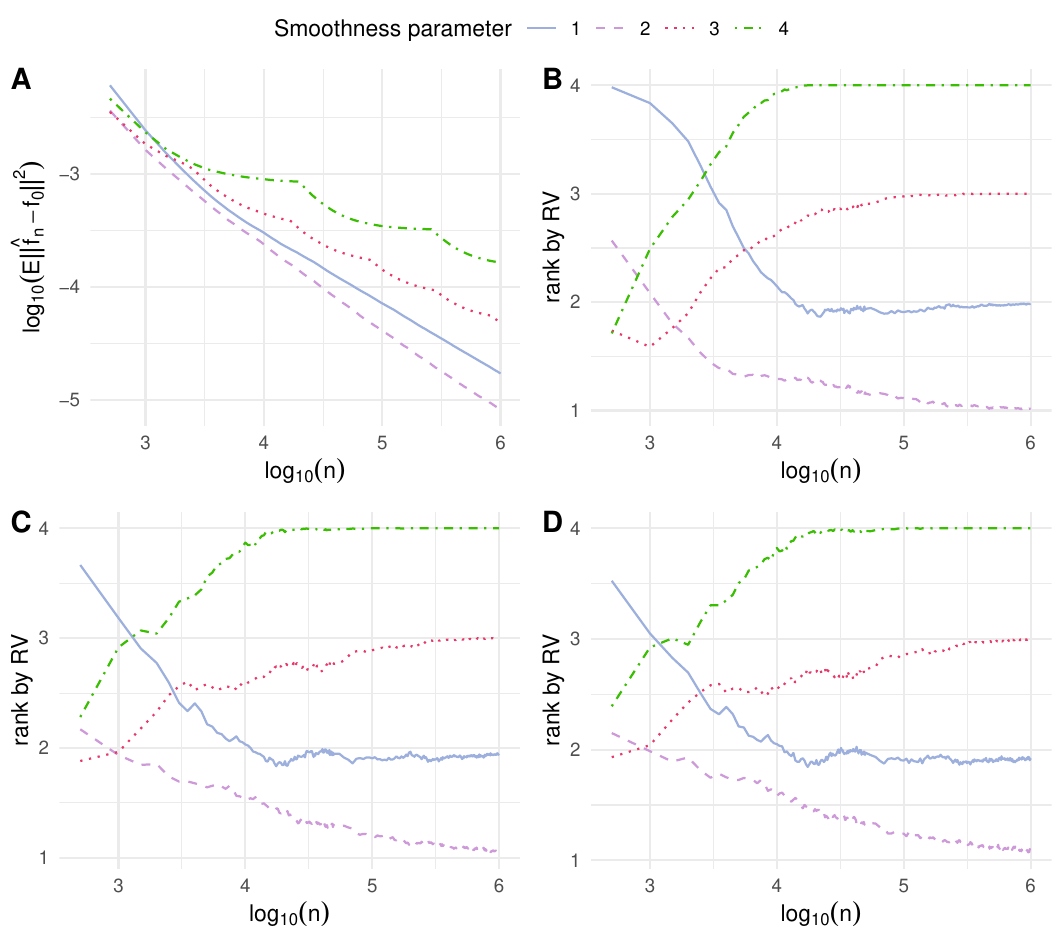}
    \caption{Model selection results of \Cref{exa:1}. (A) the true MSE of the candidate models; (B)-(D), average rankings of the models at different sample sizes over 500 repetitions, according to wRV, with weighting exponent $\xi = 1, 3,4$, respectively.}
    \label{fig: theoretical larger xi}
\end{figure}

In the main text, we examined the numerical performance of wRV for $\xi = 0, 1,2$. The same simulation is performed under $\xi = 3,4$ and the supplementary results are presented in \Cref{fig: theoretical larger xi}. Larger $\xi$ is more sensitive to the rank change of the model sequences and has a shorter delayed time frame. 

\end{document}